\documentclass[11pt,letterpaper,reqno,oneside]{amsart}
\makeatletter
\def\@settitle{\begin{center}\baselineskip14\p@\relax\normalfont\uppercasenonmath\@title\@title\end{center}}
\makeatother
\usepackage{scalerel,stackengine,relsize}
\stackMath
\usepackage{amssymb}
\usepackage{ifthen}
\usepackage{tikz}
\usepackage{tikz-3dplot}
\usepackage{bm}
\usepackage[hyphens]{url} 
\usepackage{amsmath,amssymb,amsfonts,amsthm,amsopn,dsfont,mathrsfs}
\usepackage{esint}
\usepackage{graphicx}
\usepackage{enumitem}
\usepackage{bigints}
\usepackage{fancyhdr}
\usepackage[T1]{fontenc}
\usepackage{mathtools}
\usepackage[utf8]{inputenc}
\usepackage{verbatim} 
\usepackage[backend=bibtex,style=numeric,bibencoding=utf8,language=auto,autolang=other,maxnames=10]{biblatex}
\usepackage[a4paper, total={6in, 9in}]{geometry}
\allowdisplaybreaks

\theoremstyle{plain}
\newtheorem{thm}{Theorem}[section]
\newtheorem*{thm*}{Theorem}
\newtheorem{prop}{Proposition}[section]
\newtheorem{lem}{Lemma}[section]
\newtheorem{cor}{Corollary}[section]
\theoremstyle{definition}
\newtheorem{dfn}{Definition}[section]
\theoremstyle{remark}
\newtheorem{rem}{Remark}[section]
\newtheorem{exm}{Example}[section]
\numberwithin{equation}{section}

\renewcommand{\r}{\mathbb{R}}
\renewcommand{\c}{\mathbb{C}}
\newcommand{\cp}{\mathbb{CP}}

\newcommand{\D}{\mathcal{D}}
\newcommand{\opwedge}{\mathlarger{\mathlarger{\wedge}}}
\newcommand{\n}{\mathbb{N}}
\newcommand{\m}{\mathbb{M}}
\renewcommand{\le}{\leqslant}
\renewcommand{\ge}{\geqslant}

\renewcommand{\H}{\mathscr{H}}
\renewcommand{\L}{\mathcal{L}}
\newcommand{\F}{\mathscr{F}}
\DeclareRobustCommand{\rchi}{{\mathpalette\irchi\relax}}
\newcommand{\irchi}[2]{\raisebox{\depth}{$#1\chi$}}

\DeclareMathOperator{\vol}{vol}

\DeclareMathOperator{\Id}{Id}

\DeclareMathOperator{\dist}{dist}
\DeclareMathOperator{\Lip}{Lip}
\DeclareMathOperator{\spt}{spt}

\DeclareMathOperator{\reg}{RegVal}

\DeclareMathOperator{\sing}{Sing}
\newcommand{\res}{\mathbin{\vrule height 1.6ex depth 0pt width
		0.13ex\vrule height 0.13ex depth 0pt width 1.3ex}}
\newcommand\reallywidehat[1]{%
	\savestack{\tmpbox}{\stretchto{%
			\scaleto{%
				\scalerel*[\widthof{\ensuremath{#1}}]{\kern-.6pt\bigwedge\kern-.6pt}%
				{\rule[-\textheight/2]{1ex}{\textheight}}
			}{\textheight}%
		}{0.5ex}}%
	\stackon[1pt]{#1}{\tmpbox}%
}

\title[]{The Unique Tangent Cone Property for Weakly Holomorphic Maps into Projective Algebraic Varieties}
\author[R. Caniato]{Riccardo Caniato}  
\address{Department of Mathematics \\ ETH Zürich\\ Zürich\\ Switzerland}
\email{riccardo.caniato@math.ethz.ch}
\author[T. Rivière]{Tristan Rivière}  
\address{Department of Mathematics \\ ETH Zürich\\ Zürich\\ Switzerland}
\email{tristan.riviere@math.ethz.ch}

\begin{document}
\begin{abstract}
In the present paper, we establish the uniqueness of tangent maps for general weakly holomorphic and locally approximable maps from an arbitrary almost complex manifold into projective algebraic varieties. As a byproduct of the approach and the techniques developed we also obtain the {\it unique tangent cone property} for a special class of non-rectifiable positive pseudo-holomorphic cycles. This approach gives also a new proof of the main result by C.Bellettini in \cite{bellettini-p} on the uniqueness of tangent cones for positive integral $(p,p)$-cycles in arbitrary almost complex manifolds. 
\end{abstract}
\maketitle
\tableofcontents


\section{Introduction}
\subsection{Weakly holomorphic and locally approximable maps}
For the purposes of this introduction, we always denote by $M$ a connected smooth manifold without boundary and we will need to endow $M$ with an arbitrarily chosen reference metric $g$. Since our results are local, all our discussion will be totally independent or, in any case, not essentially effected by such arbitrary choice.  
Moreover, throughout the paper we always stick to the following conventions:
\begin{enumerate}
    \item for  $K\subset M$ compact, $W^{1,2}(K)$ is the closure of the space $C^{\infty}(K)$ with respect to the strong topology induced by the $W^{1,2}$-norm, given by
    \begin{align*}
        ||u||_{W^{1,2}(K)}:=\bigg(\int_K|u|^2\, d\vol_g\bigg)^{1/2}+\bigg(\int_K|du|_g^2\, d\vol_g\bigg)^{1/2},
    \end{align*}
    for every $u\in C^{\infty}(K)$;
    \item $W_{loc}^{1,2}(M):=\{u\in W^{1,2}(K) \mbox{ for every }\,K\subset M\mbox{ compact}\}$;
    \item for any $k\in\n$, the space $W_{loc}^{1,2}\big(M,\r^k\big)$ is the real vector space of functions $u=(u_1,...,u_k)$ such that $u_j\in W_{loc}^{1,2}(M)$, $j=1,...,k$;
    \item given a closed (i.e. connected, compact and without boundary) smooth manifold $N$ and a smooth isometric embedding $N\hookrightarrow\r^k$, for some $k\in\n$ large enough, we let
    \begin{align*}
        W_{loc}^{1,2}(M,N):=\big\{u\in W_{loc}^{1,2}\big(M,\r^k\big) \mbox{ s.t. } u(x)\in N, \mbox{ for } \vol_g\mbox{-a.e. } x\in M\big\}.        
    \end{align*}
\end{enumerate}
\begin{dfn}
\label{definition of weakly holomorphic map}
Let $M$ be any even-dimensional smooth manifold without boundary and let $J$ be a Lipschitz almost complex structure on $M$. Let $(N,J_N)$ be any closed smooth almost complex manifold. We say that a map $u\in W_{loc}^{1,2}(M,N)$ is \textbf{weakly} $(J,J_N)$-\textbf{holomorphic} if 
\begin{align*}
    du_x(JX)=J_Ndu_x(X), \qquad\mbox{ for } \vol_g\mbox{-a.e. } x\in M,\, \forall\, X\in T_xM.
\end{align*}
\end{dfn}
Whenever we don't need to specify which couple of complex structures is involved in the previous definition, we simply say that the map $u$ is weakly pseudo-holomorphic or even just weakly holomorphic to lighten the notation.  

Assume that $M$ is any even-dimensional smooth manifold without boundary and let $J$ be a Lipschitz almost complex structure on $M$ which admits a compatible symplectic form $\Omega$, meaning that the bilinear form $(X,Y)\mapsto \Omega(X,JY)$ defines a Lipschitz Riemannian metric on $M$. We will show (Lemma \ref{weaklyharmonic}) that in this particular framework any weakly $(J,J_N)$-holo-morphic map taking values into a closed smooth almost K\"ahler manifold $N$ is weakly harmonic, i.e. 
\begin{align*}
    \left.\ \frac{d}{dt}\int_M\left|d\big(\pi_{N}\circ(\Phi\circ u+tX)\big)\right|_g^2\, d\vol_g\right|_{t=0}=0, \qquad \forall\, X\in C^{\infty}\big(M,\r^k\big),
\end{align*}
where $\pi_N:W\rightarrow N$ is the nearest-point projection into $N$, defined on a suitable tubular neighbourhood $W$ of $N$, and $\Phi:N\hookrightarrow\r^k$ denotes a smooth, isometric embedding of $N$ into $\r^k$.
Nevertheless, it is well-known that no regularity is ensured for weakly harmonic maps when the dimension of the domain is larger than 2 (see \cite{tristan}). Thus, we will need to prescribe some additional condition in order to get that the map $u$ is at least stationary harmonic, i.e.
\begin{align*}
    \left.\ \frac{d}{dt}\int_M\left|d(u\circ\Phi_t)\right|_g^2\, d\vol_g\right|_{t=0}=0, 
\end{align*}
for any smooth one-parameter family of diffeomorphisms $\Phi_t$ of $M$ with compact support. We will show (see Lemma \ref{stationaryharmonic}) that imposing the following local, strong approximability property with respect to the $W^{1,2}$-norm sufficies to our purposes. 
\begin{dfn}
\label{definition of locally approxiamble map}
Let $M$ be a smooth manifold without boundary and $N$ be a closed smooth manifold. We say that a map $u\in W_{loc}^{1,2}(M,N)$ is \textbf{locally (strongly) approximable} with respect to the $W^{1,2}$-norm if for every open set $U\subset M$ such that $U$ is diffeomorphic to some euclidean ball there exists a sequence of smooth maps $\{u_j\}_{j\in\n}\subset C^{\infty}(U,N)$ such that $u_j\rightarrow u$ as $j\rightarrow+\infty$, strongly in $W^{1,2}(U,N)$. 
\end{dfn}
If a map $u$ is locally approximable, then the following cohomological condition follows easily:
\begin{align}
\label{cohomological condition}
    d(u^*\omega)=0, \qquad \mbox{distributionally on } M,
\end{align}
for every closed $2$-form $\omega\in\Omega^2(N)$. We refer the reader to \cite{bethuel-coron-demengel} for further reading concerning the deep link between local approximability and \eqref{cohomological condition}.

In the most general case that we will address, i.e. when $J$ doesn't admit a compatible symplectic form (even locally), weakly holomorphic and locally approximable maps are not stationary harmonic. Nevertheless, they are \textbf{almost stationary harmonic}, in the sense that they satisfy a perturbed version of the harmonic map equation and that there exists $C>0$ such that
\begin{align*}
    \left.\ \frac{d}{dt}\int_M\left|d(u\circ\Phi_t)\right|_g^2\, d\vol_g\right|_{t=0}\ge-C\int_{M}\lvert X\rvert\lvert du\rvert^2, 
\end{align*}
with $X:=\partial_t\Phi_t|_{t=0}$ for any smooth one-parameter family of diffeomorphisms $\Phi_t$ of $M$ with compact support. We underline that such maps were also studied before by C. Bellettini and G. Tian in \cite{bellettini-tian}.
\subsection{Statement of the main results and previous literature} 
Given any $\rho\in(0,+\infty)$, we denote by $B_{\rho}\subset\r^{2m}$ the open unit ball in $\r^{2m}$ centred at the origin and having radius $\rho$. When we simply write $B$, we always mean the open unit ball $B_1\subset\r^{2m}$. From now on, for every $\rho\in(0,1)$ we let $\Phi_{\rho}:B_{\rho}\rightarrow B$ be given by $\Phi_{\rho}(x):=\rho^{-1}x$, for every $x\in B_{\rho}$. 

Let $M$ be a smooth manifold without boundary and $N$ be a closed smooth manifold. Consider any map $u\in W_{loc}^{1,2}(M,N)$. Given a point $x_0\in M$, pick any coordinate chart $\varphi:U\subset M\rightarrow B$ with relatively compact domain $U$ at $x_0$, i.e. such that $x_0\in U$ and $\varphi(x_0)=0$. The family of the \textbf{blow-ups} of $u$ at the point $x_0$, denoted by $\{u_{\rho}\}_{\rho\in (0,1)}\subset W^{1,2}(U,N)$, is given by $u_{\rho}:=u\circ\varphi^{-1}\circ\Phi_{\rho}^{-1}\circ\varphi$, for every $\rho\in (0,1)$. If such family is bounded in $W^{1,2}(U,N)$, by standard compactness arguments it follows that for every sequence $\rho_k\rightarrow 0^+$ as $k\rightarrow+\infty$ there exists a subsequence $\{\rho_{k_j}\}_{j\in\n}$ such that $u_{\rho_{k_j}}\rightharpoonup u_{\infty}\in W^{1,2}(U,N)$, weakly in $W^{1,2}(U,N)$. We say that $u_{\infty}$ is a \textbf{tangent map} for $u$ at the point $x_0$. Any tangent map at $x_0$ is meant to represent a picture of the map $u$ when one gets closer and closer to $x_0$. Such limiting configuration may very well depend on the sequence $\{\rho_{k_j}\}_{j\in\n}$ that we have chosen to approach $x_0$. If this is not the case, we say that the map $u$ \textbf{has a unique tangent map} at the point $x_0$. 

In the present paper, we aim to give a complete and self-contained proof of the following theorem.
\begin{thm}
\label{main theorem introduction}
Let $m,n\in\n_0$ be such that $m\ge 2$. Let $M$ be a smooth $2m$-dimensional manifold without boundary and let $J$ be any Lipschitz almost complex structure on $M$. Let $N\subset\cp^n$ be a projective algebraic variety. Assume that $u\in W_{loc}^{1,2}(M,N)$ is weakly $(J,j_n)$-holomorphic and locally approximable, where $j_n$ stands for the standard complex structure on $\cp^n$ (restricted to $N$).

Then, $u$ has a unique tangent map at every point. 
\end{thm}
As we have seen before, if $J$ admits a compatible symplectic form then weakly holomorphic and locally approximable maps are a special subclass in the much wider family of stationary harmonic maps (see previous subsection 1.1). Hence, our main result relates with the whole set of well-know facts concerning stationary harmonic maps between manifolds (see e.g. \cite{evans-harmonic-maps}, \cite{bethuel-stationary-harmonic}, \cite{riviere-struwe}). In particular, both the existence of tangent maps at every point and Theorem \ref{epsilonreg} are immediate consequences of the general theory of stationary harmonic maps. However, nothing can be said a priori about uniqueness of tangent maps to general stationary harmonic maps, since B. White (see \cite{white-tangent-maps}) has shown that such property might fail even for energy-minimizing maps at their singular points. Nevertheless, whenever the target manifold is analytic, uniqueness of tangent maps was proved to hold for energy-minimizing harmonic maps by L. Simon in \cite{simon-asymptotics}. Hence, since any projective algebraic variety is analytic, if weakly holomorphic and locally approximable maps were energy-minimizing then Theorem 1.1 would be a direct consequence of Simon's result. Unfortunately, it's not hard to build sequences of even holomorphic maps that converge weakly but \textbf{not} strongly in $W_{loc}^{1,2}(\r^{2m},N)$ (see Example \ref{Example: weak different from strong}). Since for energy-minimizing harmonic maps weak convergence implies strong convergence (see \cite{schoen-uhlenbeck}), this sufficies to convince ourselves that weakly holomorphic and locally approximable maps are \textbf{not} energy-minimizing harmonic maps in general.
\begin{exm}\label{Example: weak different from strong}
Let $N=\mathbb{S}^2$, equipped with the standard K\"ahler structure. Let $S\in\mathbb{S}^2$ be the south-pole in and let $p_S:\mathbb{S}^2\rightarrow\r^2$ be the stereographic projection from the south-pole. For $\lambda>0$, we define the map $u_{\lambda}:\r^2\rightarrow\mathbb{S}^2$ as follows:
\begin{align*}
    u_{\lambda}(x):=p_S^{-1}(\lambda x)\, \qquad\forall\,x\in\r^2.
\end{align*}
For every $\lambda>0$, $u_{\lambda}$ is a finite-energy, orientation-preserving conformal map from $\r^2$ to $\mathbb{S}^2$. Hence, $\{u_{\lambda}\}_{\lambda>0}$ is a family of holomorphic maps in $W^{1,2}(\r^2,\mathbb{S}^2)$. An easy computation shows that $u_{\lambda}\rightharpoonup u_{\infty}\equiv S$ weakly in $W^{1,2}$ as $\lambda\rightarrow+\infty$. Nevertheless, the convergence cannot be strong because 
\begin{align*}
    \int_{\r^2}\lvert du_{\lambda}\rvert^2\, d\L^2\rightarrow 8\pi\neq 0=\int_{\r^2}\lvert du_{\infty}\rvert^2\, d\L^2 \qquad (\lambda\rightarrow+\infty).
\end{align*}
\end{exm}
Theorem \ref{main theorem introduction} was already proven when the almost complex structure $J$ on the domain is integrable by S. Sun and X. Chen in \cite{sun-chen-1}, thanks also to the previous contributions \cite{king} and \cite{harvey-shiffman} who established the optimal bound for the Hausdorff measure of the singular set, namely
\begin{align}
\label{optimal size singular set}
\H^{2m-4}\big(\sing(u)\cap K\big)<+\infty, \qquad\forall\, K\subset M \mbox{ compact.}
\end{align}
Nevertheless, the proof provided by S. Sun and X. Chen in the integrable case makes an extensive use of complex holomorphic coordinates on the base manifold and of several algebraic tools that are not available in case we work in the almost complex framework.

As far as we know, the only available result concerning the non-integrable case that can be found in literature was achieved by the second author and G. Tian in \cite{riviere-tian}. In such paper, the case of a $4$-dimensional domain $M$ is completely solved, providing also the optimal size \eqref{optimal size singular set} for the singular set. Unfortunately, the proof that is given there strongly relies on positive intersection arguments that cannot be reproduced when $m>2$. 
\subsection{Key ideas to face the non-integrable case in higher dimensions}
In view of what we have seen in subsection 1.2, we need to think of a completely new analytic approach in order to prove Theorem \ref{main theorem introduction} in its full generality. From now on, we will denote by $J_0$ the standard complex structure on $\r^{2m}\cong\c^m$.

Let $M$, $N$ and $u$ satisfy the hypotheses of Theorem \ref{main theorem introduction}. Given any point $x_0\in M$ and a local chart $\varphi:U\subset M\rightarrow B_2$ with compact domain $U$ at $x_0$ such that $J(0)=J_0$, it's clear that $u$ has a unique tangent map at $x_0$ if and only if the local representative $\tilde u:=u\circ\varphi^{-1}\in W^{1,2}(B_2,N)$ of $u$ has a unique tangent map at the origin. Notice that $\tilde u$ is weakly $(\tilde J,j_n)$-holomorphic on $B_2$, where $\tilde J:=d\varphi\circ J\circ d\varphi^{-1}$ is a Lipschitz almost complex structure on $B_2$. Moreover, a straightforward computation allows to conclude that $\tilde u$ is locally approximable on $B_2$. Hence, Theorem \ref{main theorem introduction} will be proved if we just manage to show that the statement holds in case $M=B_2\subset\r^{2m}$, $x_0=0$, $J(0)=J_0$ and $u\in W^{1,2}(B_2,N)$.

The fact that we have reduced to prove the statement on some open ball leads to a key advantage. Since $B_2$ is contractible, we can find a Lipschitz almost symplectic form $\Omega$ on $B_2$ which is compatible with $\tilde J$, i.e. the symmetric bilinear form $(X,Y)\mapsto g(X,Y):=\Omega(X,JY)$ defines a Lipschitz metric $g$ on $B_2$. Here and throughout the whole paper, by "almost symplectic form" we mean any non-degenerate $2$-form, even if not necessarily closed (the reader should be aware that the term "almost Hermitian structure" can also be found in literature to refer to the triple $(B_2,J,\omega)$). From now on, we will assume that the domain $B_2$ is endowed with such special metric $g$ and all the computations involving scalar products will be referring to this specific choice. Notice that $\Omega^k/k!$ is a \textbf{semicalibration} on $B_2$ with respect to the metric $g$, for every $k=1,....,m$. This simply amounts to the fact that the comass norm of $\Omega^k/k!$ with respect to the metric $g$, given by
\begin{align*}
    \left|\left|\frac{\Omega^k}{k!}\right|\right|_{\ast}:=\sup\bigg\{\left<\frac{\Omega_x^k}{k!},\xi\right> \mbox{ s.t. } x\in B_2 \mbox{, } \xi\in \opwedge_k\r^{2m} \mbox{ unit simple $k$-vector}\bigg\},
\end{align*}
is equal to $1$. In case $\Omega^k/k!$ were also a closed form, we would say that it is a \textbf{calibration} on $B_2$. The notion of calibration has a long and rich history. The paper which gave its name to the corresponding general mathematical notion is the famous work of Harvey and Lawson \cite{calibrated-geometries} but complex analytic submanifolds and calibrated subvarieties had been introduced before. We invite the reader to consult the works of F. Morgan \cite{morgan1} and \cite{morgan2} for a more complete introduction to this important object of geometry.

Before giving the following, fundamental definition, we recall that a normal $k$-current on $B_2$ is $k$-dimensional current $T\in\D_k(B_2)$ such that
\begin{align*}
    \m(T)\,,\m(\partial T)<+\infty.
\end{align*}
\begin{dfn}[Semicalibrated currents]
    A normal $k$-current $T\in\D_k(B_2)$ for some $k\in\n$ is said to be \textbf{semicalibrated} by a given semicalibration $\omega$ on $B_2$ if one of the following equivalent conditions hold:
    \begin{enumerate}
        \item the measure theoretic orientation $\vec T$ of $T$ is a convex linear combination of $k$-vectors semicalibrated by $\omega$ (i.e. such that their duality with $\omega$ is unitary), $||T||$-a.e. on $B_2$ where $||T||$ stands for the total variation of $T$;
        \item $\m(T)=\left<T,\omega\right>$.
    \end{enumerate}
    If $\omega$ is a calibration (i.e. if it's closed), we say that $T$ is \textbf{calibrated} by $\omega$.
\end{dfn}

One of the reasons why calibrated and semicalibrated currents have been very much studied is that calibrated $k$-cycles are homologically mass-minimizing. Indeed, assume that $T\in\D_k(B_2)$ is a cycle, calibrated by a calibration $\omega$. If we pick any other cycle $T'\in\D_k(B_2)$ in the same homology class of $T$, i.e. such that $T-T'=\partial S$ for some $S\in\D_{k+1}(B_2)$, we immediately get
\begin{align}
\label{calibrated mass-minimizing}
    \m(T)=\left<T,\omega\right>=\left<T'+\partial S,\omega\right>=\left<T',\omega\right>\le\m(T'),
\end{align}
where we used Stokes theorem, $d\omega=0$ and $\lVert\omega\rVert_{\ast}\le 1$. 
Unfortunately, it happens very often that the closeness requirement in the definition of calibration is too strong to suit certain problems, such as the one we are interested in. Therefore, as initiated in \cite{pumberger-riviere}, it is natural to study semicalibrations and semicalibrated cycles. Substantial work on the uniqueness of tangent cones was carried out for special classes of \textbf{integral} semicalibrated cycles. In particular, such result was already obtained by C. Bellettini together with the second author in \cite{bellettini-riviere} for special legendrian integral cycles in $S^5$ and C. Bellettini in \cite{bellettini-p} has proved uniqueness of tangent cone for positive integral $(p,p)$-cycles in arbitrary almost complex manifolds. The case of positive integral $(1,1)$-cycles in arbitrary almost K\"ahler manifolds was previously covered by the main regularity result obtained by the second author and G. Tian in \cite{tristan-tian-annals}. Analogous results were obtained also by C. De Lellis, E. Spadaro and L. Spolaor in \cite{de_lellis-spadaro-spolaor}, by exploiting the fact that any integral semicalibrated $k$-cycle $T\in\D_k(B_2)$ is "almost" homologically mass-minimizing, i.e. for every $x_0\in B_2$ there are constants $C_0,r_0,\alpha_0>0$ such that
\begin{align*}
    \m\big(T\res B_{\rho}(x_0)\big)\le\m\big((T+\partial S)\res B_{\rho}(x_0)\big)+C_0\rho^{k+\alpha_0}, \qquad\forall\, 0<\rho<\rho_0
\end{align*}
and for all $S\in\D_{k+1}(B_2)$ such that $\spt(S)\subset B_{\rho}(x_0)$ (compare with the stronger property \eqref{calibrated mass-minimizing} that holds for calibrated cycles).

Nevertheless, very little is known so far concerning uniqueness of tangent cone when the rectifiability hypothesis is dropped. In general this is not true, counter-examples have been given initially by C. O. Kiselman in \cite{kiselman} and then generalized in \cite{BDM}. On the other hand, a positive result on this matter was obtained by C. Bellettini in \cite{bellettini-1}, where the author proves that the tangent cone to normal positive $(1,1)$-cycles is unique at any point where the density does not have a jump with respect to all of its values in a neighborhood.

As we will see below, the proof of Theorem \ref{main theorem introduction} will be further reduced to the problem of showing uniqueness of tangent cones for a special class of non-rectifiable, semicalibrated $(2m-2)$-cycles on $B_2$. Thus, the present paper is meant to be a contribution to this so far still fairly open class of problems. 

Let $u$ satisfy the assumptions of Theorem \ref{main theorem introduction} with $M=B_2$. We can associate to the map $u$ the $(2m-2)$-current $T_u\in\D_{2m-2}(B)$ given by
\begin{align*}
    \left<T_u,\alpha\right>:=\int_Mu^*\omega_{\cp^n}\wedge\alpha, \qquad\forall\,\alpha\in\D^{2m-2}(B).
\end{align*}
We can show (see section 2) that $T_u$ satisfies the following properties:
\begin{enumerate}
    \item $T_u$ is a cycle, i.e. $\partial T_u=0$ in the sense of currents. 
    \item $T_u$ is normal, with
    \begin{align*}
        \m(T_u)=\int_{B}|u^*\omega_{\cp^n}|_{\ast}\, d\vol_g=\frac{1}{2}\int_{B}|du|_g^2\, d\vol_g<+\infty.
    \end{align*}
    \item $T_u$ is semicalibrated by $\displaystyle{\frac{\Omega^{m-1}}{(m-1)!}}$.
\end{enumerate}
For the reader's convenience, we recall at this point that we say that a current $T$ on $B$ has a \textbf{unique tangent cone} at the origin if given any two sequences $\{\rho_k\}_{k\in\n}\subset (0,1)$ and $\{\rho_k'\}_{k\in\n}\subset (0,1)$ such that 
\begin{enumerate}
    \item $\rho_k\rightarrow 0^+$ and $\rho_k'\rightarrow 0^+$, as $k\rightarrow +\infty$,
    \item $(\Phi_{\rho_k})_{\ast}T\rightharpoonup C_{\infty}$ and $(\Phi_{\rho_k'})_{\ast}T\rightharpoonup C_{\infty}'$, as $k\rightarrow +\infty$,
\end{enumerate}
it follows that $C_{\infty}=C_{\infty}'$.

In section 6.3 we will show that uniqueness of tangent cone for the $(2m-2)$-dimensional cycle $T_u$ and for its "localizations" in the target manifold can be used in order to achieve a full proof of the uniqueness of the tangent map for $u$ at the origin. Therefore, most of our efforts will be devoted to the proof of the following statement. 
\begin{thm}
\label{main theorem}
Let $m,n\in\n_0$ be such that $m\ge 2$. Let $B\subset\r^{2m}$ be the open unit ball in $\r^{2m}$ and let $J$ be a Lipschitz almost complex structure on $B_2:=2 B$ such that $J(0)=J_0$. Assume that $u\in W^{1,2}(B,\cp^n)$ is weakly $(J,j_n)$-holomorphic and locally approximable, where $j_n$ stands for the standard complex structure on $\cp^n$. 

Then, the $(2m-2)$-cycle $T_u\in\D_{2m-2}(B)$ has a unique tangent cone at the origin.
\end{thm}

This last paragraph is dedicated to explain the main new ideas that we have introduced in order to prove Theorem \ref{main theorem}. The whole proof is based on the fact that the level sets of any weakly $(J,j_n)$-holomorphic and locally approximable map are rectifiable, $J$-holomorphic cycles. This fact is proved in section 5. By applying a slicing procedure on the right-hand-side of the monotonicity formula \eqref{monotonicity formula below} (see appendix A), we get a foliation of the region of integration into rectifiable, almost $J$-holomorphic curves (see Definition \ref{definition of almost pseudoholomorphic curve} and Remark \ref{remark about almost pseudoholomorphic curves}). By localizing properly in the target, integrating what we call the fundamental Morrey type estimate for almost $J$-holomorphic curves (see section 4) and passing then to the limit as the localization sets invade the codomain, we finally get uniqueness of tangent cone for the $2$-dimensional current $\big(T_u\res\pi^*\omega_{\cp^{m-1}}^{m-2}\big)/(m-2)!$, where $\pi:\c^m\smallsetminus\{0\}\rightarrow\cp^{m-1}$ is the standard projection map (see the first paragraph of section 2). Then, the statement of Theorem \ref{main theorem} follows as shown in the conclusion of section 6.2. 

\subsection{Final comments and open problems}
We would like to stress that our approach could also give an alternative proof of the uniqueness of tangent cone for integral $(p,p)$-cycles on almost complex manifolds, which was previously obtained by C. Bellettini in \cite{bellettini-p}. This could be achieved by considering maps $u$ that are more and more concentrated on just one rectifiable pseudo-holomorphic set in the domain (see Remark \ref{remark bellettini}).

We also believe that the method that we have developed in this work could be useful in order to proceed further in the analysis of the singular set of weakly holomorphic and locally approximable maps. In particular, we conjecture that the optimal bound \eqref{optimal size singular set} on the size of the singular set of such maps could be achieved as a further development, also exploiting Theorem \ref{main theorem introduction}. Furthermore, an interesting open problem concerns the generalization of our result to arbitrary almost K\"ahler target manifolds.

Taking a wider look and abandoning the framework of weakly holomorphic maps, there are plenty of other related problems that would deserve to be studied more deeply in view of recent developments in the field. In particular, the aim is to invent new analytic techniques that are robust enough to survive the non-availability of holomorphic coordinates in the almost complex non-integrable setting. Among all these problems, for sake of brevity we just mention uniqueness of tangent cone for triholomorphic maps in hyper-K\"ahler geometry (see e.g. \cite{bellettini-tian}) and for Hermitian Yang-Mills connections (see \cite{sun-chen-2}, \cite{sun-chen-3}).


\section{Almost monotonicity formula and tangent cones}

First, let us fix the notation that we will use throughout the present paper, whenever it won't be differently specified.
We denote by $B\subset\r^{2m}$ be the open unit ball in $\r^{2m}$, with $m\ge 2$. We let $J$ be a Lipschitz almost complex structure on $B_2:=2B$ such that $J(0)=J_0$. We let $\Omega$ be a Lipschitz almost symplectic form on $B_2$ which is compatible with $J$ and such that $\Omega(0)=\Omega_0$, where $\Omega_0$ stands for the standard symplectic form on $\r^{2m}\cong\c^m$. We denote by $g$ the Lipschitz metric on $B_2$ given by $g_x(v,w):=\Omega_x(v,Jw)$, for every $x\in B$ and $v,w\in\r^{2m}$. We indicate by $|\cdot|_g$ the norm induced by $g$ and by $|\cdot|$ the standard euclidean norm. Whenever we use the musical isomorphisms $"\sharp"$ and $"\flat"$ or the Hodge $\ast$-operator, we always take as a reference metric $g$. Finally, $\pi:B\smallsetminus\{0\}\rightarrow\cp^{m-1}$ denotes the standard projection map given by
\begin{align*}
    \pi(x_1,y_1,...,x_m,y_m):=[x_1+iy_1:...:x_m+iy_m]
\end{align*}
for every $(x_1,y_1,...,x_m,y_m)\in B\smallsetminus\{0\}$.
\begin{rem}
\label{remark compact=independence reference metric}
Notice that the fact that $g$ is Lipschitz on $B_2$ guarantees that $|\cdot|_g$ is equivalent to the euclidean norm. Consider the function $f:\overline{B}\times S^{2m-1}\rightarrow(0,+\infty)$ given by $f(x,v):=|v|_{g(x)}^2$, where $S^{2m-1}\subset\r^{2m}$ is the unit sphere in $\r^{2m}$ with respect to the euclidean norm. Since $f$ is continuous on the compact set $\overline{B}\times S^{2m-1}$, by Weierstrass theorem and by definition of $S^{2m-1}$, there exists a constant $G>0$ such that 
\begin{align*}
    \frac{1}{G}|v|^2\le f(x,v)=|v|_{g(x)}^2\le G|v|^2, \qquad\forall\, x\in\overline{B},\,\forall\, v\in S^{2m-1}.
\end{align*}
By 2-homogeneity of the squared norm, it follows that
\begin{align}
\label{equivalence of norms}
    \frac{1}{G}|v|^2\le|v|_{g(x)}^2\le G|v|^2, \qquad\forall\, x\in\overline{B},\,\forall\, v\in\r^{2m}
\end{align}
and our claim follows.

Since $|\cdot|_g$ and $|\cdot|$ are equivalent, when we refer to the Sobolev spaces on $B$ we don't need to specify which of these two norms we use in order to define them. In fact, we will use both of them according to what suits better in the context. 
\end{rem}
\begin{lem}
\label{wedgemodlem}
Let $V$ be a $2m$-dimensional real vector space and $J$ a linear complex structure on $V$. Let $\Omega$ be a symplectic form on $V$ which is compatible with $J$.

Then
\begin{align*}
    \frac{\Omega^{m-1}}{(m-1)!}\wedge\xi\wedge J\xi=|\xi|_g^2\frac{\Omega^m}{m!}, \qquad \mbox{ for every } \xi\in V^*,
\end{align*}
where $J\xi$ is given by $(J\xi)(v):=-\xi(Jv)$, for every $v\in V$. 
\begin{proof}
Fix an $g$-orthonormal basis $\{e_j,Je_j\}_{j=1,...,m}$ of $V$, so that
\begin{align*}
    \Omega=\sum_{j=1}^me_j^*\wedge Je_j^*.
\end{align*}
First of all, notice that
\begin{align*}
    \frac{\Omega^m}{m!}&=\frac{1}{m!}\sum_{j_1,...,j_m=1}^m(e_{j_1}^*\wedge Je_{j_1}^*)\wedge...\wedge(e_{j_m}^*\wedge Je_{j_m}^*)\\
    &=e_1^*\wedge Je_1^*\wedge...\wedge e_m^*\wedge Je_m^*.
\end{align*}
Fix any $\xi\in V^*$ and decompose it along the $g$-orthonormal basis denoted by $\{e_j^*,Je_j^*\}_{j=1,...,m}$ as
\begin{align*}
    \xi=\sum_{j=1}^m(\xi_{j1}e_j^*+\xi_{j2}Je_j^*).
\end{align*}
This in turn implies that
\begin{align*}
    J\xi=\sum_{j=1}^m(\xi_{j1}Je_j^*-\xi_{j2}e_j^*)
\end{align*}
and then
\begin{align*}
    \xi\wedge J\xi&=\sum_{j,k=1}^m(\xi_{j1}\xi_{k1}e_j^*\wedge J e_k^*-\xi_{j1}\xi_{k2}e_j^*\wedge e_k^*\\
    &\quad+\xi_{j2}\xi_{k1}Je_j^*\wedge Je_k^*-\xi_{j2}\xi_{k2}Je_j^*\wedge e_k^*).
\end{align*}
Since
\begin{align*}
    \frac{\Omega^{m-1}}{(m-1)!}=\sum_{j_1,...,j_{m-1}=1}^m(e_{j_1}^*\wedge Je_{j_1}^*)\wedge...\wedge(e_{j_{m-1}}^*\wedge Je_{j_{m-1}}^*),
\end{align*}
it's easy to see that
\begin{align*}
    \frac{\Omega^{m-1}}{(m-1)!}\wedge e_j^*\wedge e_k^*=\frac{\Omega^{m-1}}{(m-1)!}\wedge Je_j^*\wedge Je_k^*=0, 
\end{align*}
for every $j,k\in\{1,...,m\}$. Moreover, we notice that
\begin{align*}
    \frac{\Omega^{m-1}}{(m-1)!}\wedge\sum_{j,k=1}^m\xi_{j1}\xi_{k1}e_j^*\wedge J e_k^*=\Bigg(\sum_{j=1}^m\xi_{j1}^2\Bigg)e_1^*\wedge Je_1^*\wedge...\wedge e_m^*\wedge Je_m^*
\end{align*}
and 
\begin{align*} 
    -\frac{\Omega^{m-1}}{(m-1)!}\wedge\sum_{j,k=1}^m\xi_{j2}\xi_{k2}J e_j^*\wedge e_k^*=\Bigg(\sum_{j=1}^m\xi_{j2}^2\Bigg)e_1^*\wedge Je_1^*\wedge...\wedge e_m^*\wedge Je_m^*.
\end{align*}
By adding all the contributions, we get
\begin{align*}
    \frac{\Omega^{m-1}}{(m-1)!}\wedge\xi\wedge J\xi&=\Bigg(\sum_{j=1}^m(\xi_{j1}^2+\xi_{j2}^2)\Bigg)e_1^*\wedge Je_1^*\wedge...\wedge e_m^*\wedge Je_m^*\\
    &=|\xi|_g^2e_1^*\wedge Je_1^*\wedge...\wedge e_m^*\wedge Je_m^*
\end{align*}
and the statement follows.
\end{proof}
\end{lem}

\begin{cor}
\label{weakholcor}
Let $n\in\n_0$ and let $(N^{2n},J_N,\omega_N)$ be a compact almost Kähler manifold. Assume that $u\in W^{1,2}(B,N)$ is weakly $(J,J_N)$-holomorphic. Then
\begin{align}
\label{weakholeq}
    u^*\omega_{N}\wedge\frac{\Omega^{m-1}}{(m-1)!}=\frac{|du|_g^2}{2}\frac{\Omega^m}{m!}, \qquad \L^{2m}\mbox{-a.e. on } B.
\end{align}
\begin{proof}
Let $E\subset B$ be the set of Lebesgue points of $du$. If $x\in E$ is such that $du_x=0$, then equation \eqref{weakholeq} holds trivially. Assume then that $du(x)\neq0$. Fix an $\omega_{N}$-orthonormal basis $\{\xi_i,J_N\xi_i\}_{i=1}^n$ of $T_{u(x)}N$, so that
\begin{align*}
    (\omega_{N})_{u(x)}=\sum_{i=1}^n \xi_i^*\wedge J_N\xi_i^*.
\end{align*}
Notice that, since $u$ is weakly $(J,J_N)$-holomorphic, it holds that
\begin{align*} 
    (u^*\omega_{N})_x=\sum_{i=1}^nu^*\xi_i^*\wedge u^*J_N\xi_i^*=\sum_{i=1}^nu^*\xi_i^*\wedge Ju^*\xi_i^*,
\end{align*}
since $\big(u^*(J_N\xi^*)\big)(v)=\big(J(u^*\xi^*)\big)(v)$ for every $v\in T_xB\cong\r^{2m}$ follows from the definition of $J_N\xi^*$ and $J(u^*\xi^**)$ (see Lemma \ref{wedgemodlem}).  
Thus in particular,
\begin{align*}
    (u^*\omega_{N})_x\wedge\Omega_x^{m-1}=\sum_{i=1}^nu^*\xi_i^*\wedge Ju^*\xi_i^*\wedge\Omega_x^{m-1}.
\end{align*}
By applying Lemma \ref{wedgemodlem} with $\Omega=\Omega_x$ and $\xi=u^*\xi_i$ for every $i=1,...,n$, we get that
\begin{align*}
    (u^*\omega_{N})_x\wedge\frac{\Omega_x^{m-1}}{(m-1)!}=\Bigg(\sum_{i=1}^n|u^*\xi_i^*|_g^2\Bigg)\frac{\Omega_x^m}{m!}=\frac{|du_x|_g^2}{2}\frac{\Omega_x^m}{m!}.
\end{align*}
The statement follows immediately.
\end{proof}
\end{cor}

\begin{lem}
\label{calibratedlem}
Let $n\in\n_0$ and let $(N^{2n},J_N,\omega_N)$ be a closed almost Kähler smooth manifold. Assume that $u\in W^{1,2}(B,N)$ is weakly $(J,J_N)$-holomorphic and locally approximable.

Then, $T_u$ is a normal $(2m-2)$-cycle on $B$ semicalibrated by $\displaystyle{\frac{\Omega^{m-1}}{(m-1)!}}$. 
\begin{proof}
First, we claim that $T_u$ is a cycle, i.e. that $\partial T_u=0$. Indeed, by Stokes theorem and since $d(u^*\omega_{N})=0$ holds distributionally on $B$ by local approximability of $u$, for every fixed $\alpha\in\D^{2m-3}(B)$ we get
\begin{align*}
    \left<\partial T_u,\alpha\right>&=\left<T_u,d\alpha\right>=\int_Bu^*\omega_{N}\wedge d\alpha=0.
\end{align*}
In order to conclude, we just need to show that
\begin{align*}
    \left<T_u,\frac{\Omega^{m-1}}{(m-1)!}\right>=\m(T)<+\infty.
\end{align*}
Notice that, by Corollary \ref{weakholcor}, it holds that
\begin{align*}
    \left<T_u,\frac{\Omega^{m-1}}{(m-1)!}\right>=\int_Bu^*\omega_{N}\wedge\frac{\Omega^{m-1}}{(m-1)!}=\frac{1}{2}\int_B|du|_g^2\, d\vol_g<+\infty,
\end{align*}
since $du\in L^2(B;\r^{2m}\otimes u^*TN)$. We claim that
\begin{align*}
    \m(T_u)=\frac{1}{2}\int_B|du|_g^2\, d\vol_g.
\end{align*}
Indeed, fix any Lebesgue point $x\in B$ for $du$ and let $\{e_1,Je_1,...,e_m,Je_m\}$, $\{\xi_1^*,j\xi_1^*,...,\xi_n^*,j\xi_n^*\}$ be orthonormal bases of $T_xB$ and $T_{u(x)}N$ respectively. Then, we have 
\begin{align*}
    \langle (u^*\omega_N)_x,e_k\wedge Je_h\rangle&=\sum_{i=1}^{n}(u^*\xi_i^*\wedge Ju^*\xi_i^*)(e_k,Je_h)\\
    &=\sum_{i=1}^{n}(u^*\xi_i^*)(e_k)(Ju^*\xi_i^*)(Je_h)-(u^*\xi_i^*)(Je_h)(Ju^*\xi_i^*)(e_k)\\
    &=\sum_{i=1}^{n}(u^*\xi_i^*)(e_k)(u^*\xi_i^*)(e_h)+(u^*\xi_i^*)(Je_h)(u^*\xi_i^*)(Je_k)\\
    &\le\frac{1}{2}\sum_{i=1}^{n}\big(\lvert u^*\xi_i^*(e_k)\rvert^2+\lvert u^*\xi_i^*(e_h)\rvert^2\\
    &\quad+\lvert u^*\xi_i^*(Je_k)\rvert^2+\lvert u^*\xi_i^*(Je_h)\rvert^2\big)\\
    &\le\frac{1}{2}\sum_{i=1}^{n}\sum_{j=1}^m\big(\lvert u^*\xi_i^*(e_j)\rvert^2+\lvert u^*\xi_i^*(Je_j)\rvert^2\big)=\frac{\lvert du_x\rvert_g^2}{2}, 
\end{align*}
for every $k,h=1,...,m$. Moreover, by exactly the same computation we get
\begin{align*}
    \langle (u^*\omega_N)_x,e_k\wedge e_h\rangle=\langle (u^*\omega_N)_x,Je_k\wedge Je_h\rangle=0, \qquad\forall\,k,h=1,...,m.
\end{align*}
Thus, for every unit and simple $2$-vector $v_1\wedge v_2$ with $v_1,v_2\in T_xB$ we have  
\begin{align*}
    \langle (u^*\omega_N)_x,v_1\wedge v_2\rangle\le\frac{\lvert du_x\rvert_g^2}{2}.
\end{align*}
Consider the unit vector 
\begin{align*}
    v:=\frac{1}{\sqrt{m}}\sum_{\substack{i=1,...,m\\ i\text{ odd}}}e_i+\frac{1}{\sqrt{m}}\sum_{\substack{i=1,...,m\\ i\text{ even}}}Je_i\in T_xB
\end{align*}
and notice that
\begin{align*}
    \langle (u^*\omega_N)_x,v\wedge Jv\rangle=\bigg\langle (u^*\omega_N)_x,\frac{\Omega_x}{m}\bigg\rangle=\ast\bigg( (u^*\omega_N)_x,\wedge\frac{\Omega_x^{m-1}}{(m-1)!}\bigg)=\frac{\lvert du_x\rvert_g^2}{2}, 
\end{align*}
By definition of comass norm and since $x\in B$ was any arbitrary Lebesgue point of $du$, we conclude that
\begin{align*}
    \left|u^*\omega_{N}\right|_{*}=\frac{|du|_g^2}{2}, \qquad \vol_g\mbox{-a.e. on } B.
\end{align*}
Moreover, since $T_u$ is the integration current induced by $u^*\omega_{N}$, it holds that 
\begin{align*}
    \m(T_u)=\int_B|u^*\omega_{N}|_{*}\, d\vol_g, \qquad \mbox{ for every open set } U\subset\subset B.
\end{align*}
The statement then follows. 
\end{proof}
\end{lem}
Before stating the following fundamental proposition, we recall the following notation. Given any $x_0\in B$, we define 
\begin{align*}
    \Omega_{t,x_0}:=(dR_{x_0}\wedge\Omega)\res\nu_{x_0}
\end{align*}
where $R_{x_0}:=\lvert\cdot-x_0\rvert$, $\nu_{x_0}:=(dR_{x_0})^{\sharp}$ and with $"\res"$ we denote the interior product. We call $\Omega_{t,x_0}$ the \textbf{tangential part of $\Omega$ with respect to $x_0$}. Such notation is analogous to the one used in \cite{pumberger-riviere}. 
\begin{prop}[Almost monotonicity formula]
\label{monotonicity formula}
Let $n\in\n_0$ and let the triple $(N^{2n},J_N,\omega_N)$ denote a closed almost Kähler smooth manifold. Assume that $u\in W^{1,2}(B,N)$ is weakly $(J,J_N)$-holomorphic and locally approximable.

Then, there exists $A=A(\Lip(\Omega))\ge 0$ such that
\begin{align}
\label{monotonicity formula below}
\nonumber
    e^{A\rho}(1+A\rho)\frac{\m\big(T_u\res B_{\rho}(x_0)\big)}{\rho^{2m-2}}&-e^{A\sigma}(1+A\sigma)\frac{\m\big(T_u\res B_{\sigma}(x_0)\big)}{\sigma^{2m-2}}\\
    &\ge\int_{B_{\rho}(x_0)\smallsetminus B_{\sigma}(x_0)}\frac{1}{|\cdot-x_0|^{2m-2}}u^*\omega_N\wedge\frac{\Omega_{t,x_0}^{m-1}}{(m-1)!}
\end{align}
and 
\begin{align}
\label{monotonicity formula above}
\nonumber
    e^{-A\rho}(1-A\rho)\frac{\m\big(T_u\res B_{\rho}(x_0)\big)}{\rho^{2m-2}}&-e^{-A\sigma}(1-A\sigma)\frac{\m\big(T_u\res B_{\sigma}(x
    _0)\big)}{\sigma^{2m-2}}\\
    &\hspace{-3mm}\le\int_{B_{\rho}(x_0)\smallsetminus B_{\sigma}(x_0)}\frac{1}{|\cdot-x_0|^{2m-2}}u^*\omega_N\wedge\frac{\Omega_{t,x_0}^{m-1}}{(m-1)!},
\end{align}
for every $x_0\in B$ and $0<\sigma<\rho\le r_{x_0}:=\dist(x_0,\partial B)$.
\begin{proof}
A direct computations leads immediately to
\begin{align*}
    \big(\Omega^{m-1}\big)_{t,x_0}=\Omega_{t,x_0}^{m-1}.
\end{align*}
Hence, the statement follows directly by Lemma \ref{calibratedlem} and \cite[Proposition 1]{pumberger-riviere} by simply noticing that exactly the same proof works when $\Omega$ is just Lipschitz.
\end{proof}
\end{prop}
\begin{rem}
\label{remark density}
Fix any $x_0\in B$. Notice that
\begin{align*}
    *\bigg(u^*\omega_N(x)\wedge\frac{\Omega_{t,x_0}(x)^{m-1}}{(m-1)!}\bigg)&=\left<\frac{\Omega_{t,x_0}(x)^{m-1}}{(m-1)!},\ast u^*\omega_N\right>\\
    &=\frac{1}{2}\left|\frac{\partial u}{\partial\nu_{x_0}}\right|_g^2, \qquad \mbox{ for } \L^{2m}\mbox{-a.e. } x\in B,
\end{align*}
where $\nu_{x_0}:=(dR_{x_0})^{\sharp}$ as above with $R_{x_0}:=\lvert\cdot-x_0\rvert$. 
Hence, by equation \eqref{monotonicity formula below} we conclude that the function 
\begin{align*}
    (0,r_{x_0})\ni\rho\mapsto e^{A\rho}(1+A\rho)\frac{\m(T_u\res B_{\rho})}{\rho^{2m-2}}
\end{align*}
is non-decreasing. As 
\begin{align*}
    \lim_{\rho\rightarrow 0^+}e^{A\rho}(1+A\rho)=1,
\end{align*}
we conclude that the limit
\begin{align}
\label{definition density}
    \theta(x_0,u):=\lim_{\rho\rightarrow 0^+}\frac{\m\big(T_u\res B_{\rho}(x_0)\big)}{\rho^{2m-2}}
\end{align} 
exists and is finite. We say that  $\theta(x_0,u)$ is the \textbf{density of the map $u$ at the point $x_0$}.
\end{rem}

We conclude the section by discussing the existence and the structure of tangent cones for the current $T_u$. Let's pick any sequence $\rho_k\rightarrow 0^+$ as $k\rightarrow +\infty$ and the relative blow-up sequence
$\{T_{\rho_k}:=(\Phi_{\rho_k})_{\ast}T_u\}_{k\in\n}$. Since $T_{\rho_k}$ is a cycle for every $k\in\n$ and 
\begin{align*}
    \m(T_{\rho_k})=\frac{\m(T_u\res B_{\rho_k})}{\rho^{2m-2}}\le e^{A}(1+A)\m(T_u)<+\infty,
\end{align*}
by Federer-Fleming compactness theorem we know that there exists a subsequence $\{\rho_{k_j}\}_{j\in\n}$ of $\{\rho_k\}_{k\in\n}$ such that $T_{\rho_{k_j}}\rightharpoonup C_{\infty}$ as $j\rightarrow +\infty$ in the sense of currents. Moreover, by exploiting the almost monotonicity formula, we get that any tangent cone $C_{\infty}$ is a $(2m-2)$-cycle calibrated by $\Omega_0$ and invariant under dilations, i.e. $(\Phi_{\rho})_{\ast}C_{\infty}=C_{\infty}$, for every $\rho\in(0,1)$ (see e.g. \cite[Section 3]{pumberger-riviere}).


\section{Smoothness at points with small density}
In the present section, we will assume that $\Omega$ is a symplectic form and we prove that weakly holomorphic and locally approximable maps are stationary harmonic in this particular case. Therefore, we conclude that such maps are smooth at points of small density via standard $\varepsilon$-regularity for stationary harmonic maps (Theorem \ref{epsilonreg}). 

The almost symplectic case is completely treated in \cite[Propositions 1, 3, 4]{bellettini-tian}, where it is shown that similar results hold in the almost stationary scenario. Hence, the conclusions of the present section hold even if $d\Omega\neq 0$. We just present here a simplified case in order to deal with less technicalities and draw some light on the key ideas and concepts.  
\allowdisplaybreaks
\begin{lem}[Wirtinger's inequality]
\label{weakharmoniclem}
Let $n\in\n_0$ and let $(N^{2n},J_N,\omega_N)$ be a closed almost Kähler smooth manifold.

Then, for every map $v\in W^{1,2}(B,N)$ it holds that
\begin{align}
\label{weaklyharmoniceq}
    *\bigg(v^*\omega_{N}\wedge\frac{\Omega^{m-1}}{(m-1)!}\bigg)\le\frac{|dv|_g^2}{2}, \qquad \vol_g\mbox{-a.e. on } B.
\end{align}
\begin{proof}
Let $E\subset B$ be the set of the Lebesgue points of $dv$. If $x\in E$ is such that $dv_x=0$, then \eqref{weaklyharmoniceq} holds trivially. Assume then that $x\in E$ is such that $dv_x\neq0$. Fix a $g$-orthonormal basis $\{e_{2k-1},e_{2k}:=Je_{2k-1}\}_{k=1}^m$ of $T_xB$ and an $\omega_{N}$-orthonormal basis $\{\xi_{2i-1},\xi_{2i}:=j\xi_{2i-1}\}_{i=1}^n$ of $T_{v(x)}N$, so that
\begin{align*}
    \Omega_x=\sum_{k=1}^me_{2k-1}^*\wedge e_{2k}^*
\end{align*}
and
\begin{align*}  
    (\omega_{N})_{v(x)}=\sum_{i=1}^n \xi_{2i-1}^*\wedge\xi_{2i}^*.
\end{align*}
Then, we compute
    \begin{align*}
    &*\bigg((v^*\omega_{N})_x\wedge\frac{\Omega_x^{m-1}}{(m-1)!}\bigg)\\
    &=*\sum_{k=1}^m\sum_{i=1}^nv^*\xi_{2i-1}^*\wedge v^*\xi_{2i}^*\wedge e_1^*\wedge e_2^*\wedge...\wedge\widehat{e_{2k-1}^*}\wedge\reallywidehat{e_{2k}^*}\wedge...\wedge e_{2m-1}\wedge e_{2m}\\
    &=\sum_{k=1}^m\sum_{i=1}^n(v^*\xi_{2i-1}^*\wedge v^*\xi_{2i}^*)(e_{2k-1}\wedge e_{2k})\\
    &=\sum_{k=1}^m\sum_{i=1}^n(v^*\xi_{2i-1}^*\wedge v^*\xi_{2i}^*)(e_{2k-1},e_{2k})\\
    &\le\sum_{k=1}^m\sum_{i=1}^n\left|(v^*\xi_{2i-1}^*\wedge v^*\xi_{2i}^*)(e_{2k-1},e_{2k})\right|_g\\
    &=\sum_{k=1}^m\sum_{i=1}^n\left|v^*\xi_{2i-1}^*(e_{2k-1})v^*\xi_{2i}^*(e_{2k})-v^*\xi_{2i-1}^*(e_{2k})v^*\xi_{2i}^*(e_{2k-1})\right|_g\\
    &\le\frac{1}{2}\sum_{k=1}^{2m}\sum_{i=1}^n\big(|v^*\xi_{2i-1}^*(e_k)|_g^2+|v^*\xi_{2i}^*(e_k)|_g^2\big)=\frac{|dv_x|_g^2}{2}.
\end{align*}
Thus, the statement follows.
\end{proof}
\end{lem}

\begin{lem}[Weakly holomorphic maps are weakly harmonic]
\label{weaklyharmonic}
Let $n\in\n_0$ and let $(N^{2n},J_N,\omega_N)$ be a closed almost Kähler smooth manifold.

If $u\in W^{1,2}(B,N)$ is weakly $(J,J_N)$-holomorphic, then $u$ is weakly harmonic.
\begin{proof}
Recall that we always identify $N$ as a smooth submanifold of $\r^k$, for $k$ large enough, through the smooth isometric embedding $\Phi:N\hookrightarrow\r^k$ (see Section 1). Let the map $\pi_{N}:W\subset\r^k\rightarrow N$ be the nearest point projection from a tubular neighbourhood $W$ of $\Phi(N)$ onto $N$. Fix any vector field $X\in C_c^{\infty}(B;\r^k)$. As for every $t\in\r$ the map $\pi_{N}\circ(\Phi\circ u+tX)$ belongs to $W^{1,2}(B,N)$, by exploiting Lemma \ref{weakharmoniclem} we get that
\begin{align}
\label{weakharmoniceq1}
    \int_B\left|d\big(\pi_{N}\circ(\Phi\circ u+tX)\big)\right|_g^2\, d\vol_g\ge2\int_{B}(\Phi\circ u+tX)^*\pi_{N}^*\omega_{N}\wedge\frac{\Omega^{m-1}}{(m-1)!},
\end{align}
for $t\in\r$ such that $\lvert t\rvert<\delta$ with $\delta>0$ sufficiently small. Moreover, the equality holds for $t=0$ by virtue of equation \eqref{weakholeq}. We claim that
\begin{align}
\label{weakharmoniceq2}
    \int_{B}(\Phi\circ u+tX)^*\pi_{N}^*\omega_{N}\wedge\frac{\Omega^{m-1}}{(m-1)!}=\int_{B}u^*\omega_{N}\wedge\frac{\Omega^{m-1}}{(m-1)!},
\end{align}
for every $\lvert t\rvert<\delta$. Indeed, it holds that
\begin{align*}
    &\frac{d}{dt}\Bigg(\int_{B}(\Phi\circ u+tX)^*\pi_{N}^*\omega_{N}\wedge\frac{\Omega^{m-1}}{(m-1)!}\Bigg)\\
    &=\int_{B}\frac{d}{dt}\big((\Phi\circ u+tX)^*\pi_{N}^*\omega_{N}\big)\wedge\frac{\Omega^{m-1}}{(m-1)!}\\
    &=\left<d\big(u^*(\pi_{N}^*\omega_{N})\res X\big),\ast\frac{\Omega^{m-1}}{(m-1)!}\right>\\
    &=-\int_Bu^*(\pi_{N}^*\omega_{N})\res X\wedge d\bigg(\frac{\Omega^{m-1}}{(m-1)!}\bigg)=0,
\end{align*}
where "$\res$" stands for the interior product. Hence, equation \eqref{weakharmoniceq2} follows. By using together equation \eqref{weakharmoniceq1} and  \eqref{weakharmoniceq2}, we get
\begin{align*}
\int_B\left|d\big(\pi_{N}\circ(\Phi\circ u+tX)\big)\right|_g^2\, d\vol_g\ge\int_B|du|_g^2\, d\vol_g \qquad \mbox{ for every } \lvert t\rvert<\delta,
\end{align*}
and the equality holds for $t=0$. Thus, $t=0$ is a global minimum for the differentiable function 
\begin{align*}
    t \mapsto \int_B\left|d\big(\pi_{N}\circ(\Phi\circ u+tX)\big)\right|_g^2\, d\vol_g.
\end{align*} 
Hence, we conclude that
\begin{align*}
    \left.\ \frac{d}{dt}\int_B\left|d\big(\pi_{N}\circ(\Phi\circ u+tX)\big)\right|_g^2\, d\vol_g\right|_{t=0}=0.
\end{align*}
Since our choice of $X\in C_c^{\infty}(B;\r^k)$ was arbitrary, the statement follows. 
\end{proof}
\end{lem}

\begin{lem}[Weakly holomorphic and locally approximable maps are stationary harmonic]
\label{stationaryharmonic}
Let $n\in\n_0$ and let $(N^{2n},J_N,\omega_N)$ be a closed almost Kähler smooth manifold.

If $u\in W^{1,2}(B,N)$ is weakly $(J,J_N)$-holomorphic and locally approximable, then $u$ is stationary harmonic.
\begin{proof}
Fix any vector field $X\in C_c^{\infty}(B;\r^{2m})$. Notice that the map $u\circ (\Id + tX)$ belongs to $W^{1,2}(B;N)$ for $t\in\r$ such that $\lvert t\rvert<\delta$ with $\delta>0$ sufficiently small. Then, by Lemma \ref{weakharmoniclem}, it holds that
\begin{align}
\label{stationaryeq1}
    \int_B\left|d\big(u\circ(\Id+tX)\big)\right|_g^2\, d\vol_g\ge2\int_{B}(\Id+tX)^*u^*\omega_{N}\wedge\frac{\Omega^{m-1}}{(m-1)!}
\end{align}
for every $\lvert t\rvert<\delta$ and the equality holds for $t=0$ by virtue of equation \eqref{weakholeq}. We claim that 
\begin{align}
\label{stationaryeq2}
    \int_{B}(\Id+tX)^*u^*\omega_{N}\wedge\frac{\Omega^{m-1}}{(m-1)!}=\int_{B}u^*\omega_{N}\wedge\frac{\Omega^{m-1}}{(m-1)!}, 
\end{align}
for every $\lvert t\rvert<\delta$. Indeed, as $d(u^*\omega_{N})=0$ distributionally on $B$, it holds that
\begin{align*}
    \frac{d}{dt}\Bigg(\int_{B}(\Id+tX)^*u^*\omega_{N}\wedge\frac{\Omega^{m-1}}{(m-1)!}\Bigg)&=\int_{B}\frac{d}{dt}\big((\Id+tX)^*u^*\omega_{N}\big)\wedge\frac{\Omega^{m-1}}{(m-1)!}\\
    &=\left<\L_X(u^*\omega_{N}),\ast\frac{\Omega^{m-1}}{(m-1)!}\right>\\
    &=\left<d\big(u^*\omega_{N}\res X\big),\ast\frac{\Omega^{m-1}}{(m-1)!}\right>\\
    &=-\int_{B}u^*\omega_{N}\res X\wedge d\bigg(\frac{\Omega^{m-1}}{(m-1)!}\bigg)=0.
\end{align*}
Hence, equation \eqref{stationaryeq2} follows. By using together equation \eqref{stationaryeq1} and \eqref{stationaryeq2} we get that 
\begin{align}
    \int_B\left|d\big(u\circ(\Id+tX)\big)\right|_g^2\, d\vol_g\ge\int_B|du|_g^2\, d\vol_g, \quad \mbox{ for every } \lvert t\rvert<\delta.
\end{align}
and the equality holds for $t=0$. Thus, $t=0$ is a global minimum for the differentiable function
\begin{align*}
    t\mapsto \int_B\left|d\big(u\circ(\Id +tX)\big)\right|_g^2\, d\vol_g.
\end{align*}
Hence, we conclude that
\begin{align*}
    \left.\ \frac{d}{dt}\int_B\left|d\big(u\circ(\Id+tX)\big)\right|_g^2\, d\vol_g\right|_{t=0}=0.
\end{align*}
Since our choice of $X\in C_c^{\infty}(B;\r^k)$ was arbitrary, the statement follows.  
\end{proof}
\end{lem}

The following $\varepsilon$-regularity statement follows immediately by Lemma \ref{stationaryharmonic} and \cite[Theorem 2.1]{riviere-struwe}.

\begin{thm}[$\varepsilon$-regularity for weakly holomorphic and locally approximable maps]
\label{epsilonreg}
Let $n\in\n_0$ and let $(N^{2n},J_N,\omega_N)$ be a closed almost Kähler smooth manifold.

Let $u\in W^{1,2}(B,\cp^{n})$ be weakly $(J,J_N)$-holomorphic and locally approximable. Then, there exists a threshold $\varepsilon_0=\varepsilon_0(m,n)>0$ such that whenever $\theta(x_0,u)<\varepsilon_0$ there exists ball $B_{\rho}(x_0)\subset B$ such that
$u$ is Hölder continuous (and hence smooth) on $B_{\rho}(x_0)$.
\end{thm}
\noindent
We define
\begin{align*}
    \sing(u):=\{x_0\in B \mbox{ s.t. } \theta(x_0,u)\ge\varepsilon_0\}
\end{align*}
and we say that $\sing(u)$ is the \textbf{singular set} of $u$. Moreover, by stationarity of $u$, it follows that 
\begin{align*}
    \H^{2m-2}\big(\sing(u)\big)=0.
\end{align*}


\section{The fundamental Morrey type estimate}

We aim to collect here the proofs of the (mostly technical) tools and estimates that will be used in section 5.
Throughout the present paper, given any $\H^k$-rectifiable subset $\Sigma\subset B$ equipped with an orienting $\H^k$-measurable field of unit and simple $k$-vectors $\vec\Sigma$ we denote by $[\Sigma]$ the current of integration on $\Sigma$, i.e. the $k$-dimensional current given by
\begin{align*}
    \langle[\Sigma],\alpha\rangle:=\int_{\Sigma}\langle\alpha,\vec\Sigma\rangle\,d\H^k, \qquad\forall\,\alpha\in\D^k(B).
\end{align*}

\subsection{Some technical lemmata}
\begin{dfn}[$J$-holomorphic curves]
\label{definition of pseudoholomorphic curve}
A locally $\H^2$-rectifiable subset $\Sigma\subset B$ equipped with an orienting $\H^2$-measurable field of unit and simple $2$-vectors $\vec\Sigma$ is a $J$-\textbf{holomorphic curve} if $\vec\Sigma(x)$ is $J$-invariant for $\H^2$-a.e. $x\in\Sigma$.

Moreover, if $\partial[\Sigma]=0$ we say that $\Sigma$ is \textbf{closed}.
\end{dfn}
\begin{dfn}[Almost $J$-holomorphic curves]
\label{definition of almost pseudoholomorphic curve}
A locally $\H^2$-rectifiable subset $\Sigma\subset B$ equipped with an orienting $\H^2$-measurable field of unit and simple $2$-vectors $\vec\Sigma$ is a \textbf{almost} $J$-\textbf{holomorphic curve} if there exists some $\H^2$-measurable and $J$-invariant field of $2$-vectors $\vec\Sigma_J:\Sigma\rightarrow\bigwedge_2\r^{2m}$ such that for some $\gamma\in(0,1]$, $\ell\ge 0$ it holds that
\begin{align}
    |\vec\Sigma(x)-\vec\Sigma_J(x)|\le\ell|x|^{\gamma}, \qquad\mbox{ for } \H^2\mbox{- a.e. } x\in\Sigma.
\end{align}

Moreover, if $\partial[\Sigma]=0$ we say that $\Sigma$ is \textbf{closed}.
\end{dfn}
\begin{rem}
\label{remark about almost pseudoholomorphic curves}
Given an almost $J$-holomorphic curve in $B$, we can build a $2$-dimensional varifold on $B$ associated to it in the following way. 

Let $\mathscr{G}_2(B):=B\times\text{Gr}(2,\r^{2m})$, where $\text{Gr}(2,\r^{2m})$ is the Grassmannian of the real $2$-planes in $\r^{2m}$. Notice that $\mathscr{G}_2(B)$ can be given the structure of a smooth manifold, since it is the product of two smooth manifolds.
Following the notation by W.K. Allard and L. Simon (see \cite[Chapter 8]{Simon}, \cite{allard}), a general $2$-dimensional varifold on $B$ is simply a Radon measure on $\mathscr{G}_2(B)$. Then, we may associate to an almost $J$-holomorphic curve $\Sigma\subset B$ the Radon measure on $\mathscr{G}_2(B)$ given by 
\begin{align*}
    \H^2\res\Sigma\otimes\delta_{\text{span}\{\vec\Sigma_J\}},
\end{align*}
where by $\otimes$ we denote the usual tensor product of measures and $\text{span}\{\vec\Sigma_J\}$ denotes the field of $2$-planes associated with the field of $2$-vectors $\vec\Sigma_J$.

Such objects are very close to being rectifiable varifolds but the almost tangent space of $\Sigma$ is "tilted", conveniently with respect to the purposes that will be clear in the forthcoming discussion. 

We point out explicitely that the form of these new objects is built (and therefore meaningful) just to work around the origin. We would need to consider a "shifted" version of almost $J$-holomorphic curves in order to work around an arbitrary point $x_0\in B$. 
\end{rem}
\begin{rem}
\label{pseudoholomorphic vs almost pseudoholomorphic case}
All the estimates and the results that will be presented in this section concerning closed almost $J$-holomorphic curves in $B$ are still valid for closed $J$-holomorphic curves. Indeed, any $J$-holomorphic curve is trivially almost $J$-holomorphic (just pick $\vec\Sigma_J=\vec\Sigma$, $\ell=0$ and $\gamma=1$ in Definition \ref{definition of almost pseudoholomorphic curve}).

Hence, in order to get the corresponding estimates for closed $J$-holomorphic curves it will always be sufficient to set $\vec\Sigma_J=\vec\Sigma$, $\ell=0$ and $\gamma=1$ in what follows.
\end{rem}
From now on, we will denote by $\nu_0$ the vector field on $B\smallsetminus\{0\}$ given by $\nu_0(x)=x/|x|$. Moreover, we notice that since $\Omega$ is Lipschitz and $\Omega(0)=\Omega_0$, there exists a constant $\tilde L>0$ depending only on $\Lip(\Omega)$ such that
\begin{align*}
    |\nu-\nu_0|\le\tilde L|\cdot|.
\end{align*}
\begin{prop}[Almost monotonicity formula]
\label{monotonicity formula sigma}
Let $\Sigma$ be a closed almost $J$-holomorphic curve in $B$, according to Definition \ref{definition of almost pseudoholomorphic curve}. Then, there exists a positive constant $A\ge 0$ depending only on the Lipschitz constant of $\Omega$ such that
\begin{align}
\label{monotonicity sigma below}
\nonumber
e^{A\rho+\ell\frac{\rho^{\gamma}}{\gamma}}(1+A\rho)\frac{\H^{2}(\Sigma\cap B_{\rho})}{\rho^2}&-e^{A\sigma+\ell\frac{\sigma^{\gamma}}{\gamma}}(1+A\sigma)\frac{\H^{2}(\Sigma\cap B_{\sigma})}{\sigma^2}\\
&\ge\int_{\Sigma\cap(B_{\rho}\smallsetminus B_{\sigma})}\frac{1}{|\cdot|^2}\big|\vec\Sigma_J\wedge\nu\big|_g^2\, d\H^2
\end{align}
and
\begin{align}
\label{monotonicity sigma above}
\nonumber
e^{-\big(A\rho+\ell\frac{\rho^{\gamma}}{\gamma}\big)}(1-A\rho)\frac{\H^{2}(\Sigma\cap B_{\rho})}{\rho^2}&-e^{-\big(A\sigma+\ell\frac{\sigma^{\gamma}}{\gamma}\big)}(1-A\sigma)\frac{\H^{2}(\Sigma\cap B_{\sigma})}{\sigma^2}\\
&\le\int_{\Sigma\cap\left(B_{\rho}\smallsetminus B_{\sigma}\right)}\frac{1}{|\cdot|^2}\big|\vec\Sigma_J\wedge\nu\big|_g^2\, d\H^2,
\end{align}
\begin{proof}
Throughout this proof, $R$ will denote the smooth radial vector field on $B$ given by $R(x):=x$, for every $x\in B$. We denote by $\Omega_0$ the standard symplectic form on $B^{2m}$ and we define $\Omega_1:=\Omega-\Omega_0$. Moreover, given any arbitrary form $\alpha\in\Omega^2(B)$ we denote by $\alpha_t$ the tangential part of a form with respect to the vector field $\nu$, given by
\begin{align*}
    \alpha_t:=(dr\wedge\alpha)\res\nu,
\end{align*}
according to the notation used in Proposition \ref{monotonicity formula}.

Define the normal $2$-current on $B$ given by 
\begin{align*}
    \left<[\Sigma]_J,\alpha\right>:=\int_{\Sigma}\langle\alpha,\vec\Sigma_J\rangle\, d\H^2, \qquad\forall\, \alpha\in\D^2(B).
\end{align*}
As $[\Sigma]_J$ is semicalibrated by $\Omega$, we will apply the same method that is used in \cite[Proposition 1]{pumberger-riviere}. Nevertheless, we need to take into account the fact that the $2$-current $[\Sigma]_J$ is not a cycle (though not far from being one).

Let $\varphi:[0,+\infty)\rightarrow [0,+\infty)$ be smooth, non-increasing and such that:
\begin{enumerate}
    \item $\varphi\equiv 1$ on $[0,1/2]$;
    \item $|\varphi'|\le 4$ on $[0,+\infty)$.
    \item $\varphi\equiv 0$ on $[1,+\infty)$.
\end{enumerate}
For every $0<\rho<1$, define $\varphi_{\rho}(x):=\varphi(|x|/\rho)$, for every $x\in\r^{2m}$. Notice that $\varphi_{\rho}\equiv 1$ on $B_{\rho/2}$, $\varphi_{\rho}\equiv 0$ on $\r^{2m}\smallsetminus B_{\rho}$ and $|\nabla\varphi_{\rho}|\le 4/\rho$ on $\r^{2m}$. Define
\begin{align*}
    I(\rho)&:=\int_{\Sigma}\varphi_{\rho}\langle\Omega,\vec\Sigma_J\rangle\, d\H^2=\int_{\Sigma}\varphi_{\rho}\, d\H^2,\\
    J(\rho)&:=\int_{\Sigma}\varphi_{\rho}\langle\Omega_t,\vec\Sigma_J\rangle\, d\H^2.
\end{align*}
Recall that $\L_R\Omega_0=2\Omega_0$ and compute 
\begin{align*}
    2I(\rho)&=2\int_{\Sigma}\varphi_{\rho}\langle\Omega_0,\vec\Sigma_J\rangle\, d\H^2+2\int_{\Sigma}\varphi_{\rho}\langle\Omega_1,\vec\Sigma_J\rangle\, d\H^2\\
    &=\int_{\Sigma}\langle\varphi_{\rho}d(\Omega_0\res R),\vec\Sigma_J\rangle\, d\H^2+2\int_{\Sigma}\varphi_{\rho}\langle\Omega_1,\vec\Sigma_J\rangle\, d\H^2\\
    &=\int_{\Sigma}\langle d\big(\varphi_{\rho}(\Omega_0\res R)\big),\vec\Sigma_J\rangle\, d\H^2-\int_{\Sigma}\langle d\varphi_{\rho}\wedge(\Omega_0\res R),\vec\Sigma_J\rangle\, d\H^2\\
    &\quad+2\int_{\Sigma}\varphi_{\rho}\langle\Omega_1,\vec\Sigma_J\rangle\, d\H^2\\
    &=\int_{\Sigma}\langle d\big(\varphi_{\rho}(\Omega_0\res R)\big),\vec\Sigma-\vec\Sigma_J\rangle\, d\H^2\\
    &\quad+\rho\int_{\Sigma} \frac{\partial\varphi_{\rho}}{\partial\rho}\langle dr\wedge(\Omega_0\res(\nu_0-\nu)),\vec\Sigma_J\rangle\, d\H^2\\
    &\quad+\rho\int_{\Sigma} \frac{\partial\varphi_{\rho}}{\partial\rho}\langle dr\wedge(\Omega_0\res\nu),\vec\Sigma_J\rangle\, d\H^2+2\int_{\Sigma}\varphi_{\rho}\langle\Omega_1,\vec\Sigma_J\rangle\, d\H^2\\
    &=\int_{\Sigma}\langle d\big(\varphi_{\rho}(\Omega_0\res R)\big),\vec\Sigma-\vec\Sigma_J\rangle\, d\H^2\\
    &\quad+\rho\int_{\Sigma} \frac{\partial\varphi_{\rho}}{\partial\rho}\langle dr\wedge(\Omega_0\res(\nu_0-\nu)),\vec\Sigma_J\rangle\, d\H^2\\
    &\quad+\rho\int_{\Sigma} \frac{\partial\varphi_{\rho}}{\partial\rho}\langle\Omega_0-(\Omega_0)_t,\vec\Sigma_J\rangle\, d\H^2+2\int_{\Sigma}\varphi_{\rho}\langle\Omega_1,\vec\Sigma_J\rangle\, d\H^2\\
    &=\int_{\Sigma}\langle d\big(\varphi_{\rho}(\Omega_0\res R)\big),\vec\Sigma-\vec\Sigma_J\rangle\, d\H^2+\rho\int_{\Sigma} \frac{\partial\varphi_{\rho}}{\partial\rho}\langle\Omega-\Omega_t,\vec\Sigma_J\rangle\, d\H^2\\
    &\quad+\rho\int_{\Sigma} \frac{\partial\varphi_{\rho}}{\partial\rho}\langle dr\wedge(\Omega_0\res(\nu_0-\nu)),\vec\Sigma_J\rangle\, d\H^2\\
    &\quad+2\int_{\Sigma}\varphi_{\rho}\langle\Omega_1,\vec\Sigma_J\rangle\, d\H^2-\rho\int_{\Sigma} \frac{\partial\varphi_{\rho}}{\partial\rho}\langle\Omega_1-(\Omega_1)_t,\vec\Sigma_J\rangle\, d\H^2\\
    &=\rho I'(\rho)-\rho J'(\rho)+\int_{\Sigma}\langle d\big(\varphi_{\rho}(\Omega_0\res R)\big),\vec\Sigma-\vec\Sigma_J\rangle\, d\H^2\\
    &\quad+\rho\int_{\Sigma} \frac{\partial\varphi_{\rho}}{\partial\rho}\langle dr\wedge(\Omega_0\res(\nu_0-\nu)),\vec\Sigma_J\rangle\, d\H^2\\
    &\quad-\rho\int_{\Sigma} \frac{\partial\varphi_{\rho}}{\partial\rho}\langle\Omega_1-(\Omega_1)_t,\vec\Sigma_J\rangle\, d\H^2+2\int_{\Sigma}\varphi_{\rho}\langle\Omega_1,\vec\Sigma_J\rangle\, d\H^2,
\end{align*}
which leads to
\begin{align*}
    -2\frac{I(\rho)}{\rho^3}+\frac{I'(\rho)}{\rho^2}-\frac{J'(\rho)}{\rho^2}&=-\frac{1}{\rho^3}\int_{\Sigma}\langle d\big(\varphi_{\rho}(\Omega_0\res R)\big),\vec\Sigma-\vec\Sigma_J\rangle\, d\H^2\\
    &\quad+\frac{1}{\rho^2}\int_{\Sigma} \frac{\partial\varphi_{\rho}}{\partial\rho}\langle dr\wedge(\Omega_0\res(\nu_0-\nu)),\vec\Sigma_J\rangle\, d\H^2\\
    &\quad+\frac{1}{\rho^2}\int_{\Sigma} \frac{\partial\varphi_{\rho}}{\partial\rho}\langle\Omega_1-(\Omega_1)_t,\vec\Sigma_J\rangle\, d\H^2\\
    &\quad-\frac{2}{\rho^3}\int_{\Sigma}\varphi_{\rho}\langle\Omega_1,\vec\Sigma_J\rangle\, d\H^2.
\end{align*}
Notice that
\begin{align*}
    \bigg|\frac{1}{\rho^3}\int_{\Sigma}\langle d\big(\varphi_{\rho}(\Omega_0\res R)\big),\vec\Sigma-\vec\Sigma_J\rangle\, d\H^2\bigg|&\le\ell\rho^{\gamma-1}\bigg(\frac{\H^2(\Sigma\cap B_{\rho})}{\rho^2}\bigg),
\end{align*}
\begin{align*}
    \bigg|\frac{1}{\rho^2}\int_{\Sigma} \frac{\partial\varphi_{\rho}}{\partial\rho}\langle\Omega_1-(\Omega_1)_t,\vec\Sigma_J\rangle\, d\H^2\bigg|\le\frac{2\Lip(\Omega)}{\rho}\int_{\Sigma}\frac{\partial\varphi_{\rho}}{\partial\rho}\, d\H^2=2\Lip(\Omega)\frac{I'(\rho)}{\rho},
\end{align*}
\begin{align*}
    \bigg|\frac{1}{\rho^2}\int_{\Sigma} \frac{\partial\varphi_{\rho}}{\partial\rho}\langle dr\wedge(\Omega_0\res(\nu_0-\nu)),\vec\Sigma_J\rangle\, d\H^2\bigg|\le\frac{\tilde L}{\rho}\int_{\Sigma}\frac{\partial\varphi_{\rho}}{\partial\rho}\, d\H^2=\tilde L\frac{I'(\rho)}{\rho}
\end{align*}
and
\begin{align*}
    \bigg|\frac{2}{\rho^3}\int_{\Sigma}\varphi_{\rho}\langle\Omega_1,\vec\Sigma_J\rangle\, d\H^2\bigg|\le\frac{2\Lip(\Omega)}{\rho^2}\int_{\Sigma}\varphi_{\rho}\langle\Omega_1,\vec\Sigma_J\rangle\, d\H^2=2\Lip(\Omega)\frac{I(\rho)}{\rho^2}.
\end{align*}
Hence, we conclude
\begin{align*}
    \bigg|\frac{d}{d\rho}\bigg(\frac{I(\rho)}{\rho^2}\bigg)-\frac{J'(\rho)}{\rho^2}\bigg|&=\bigg|-2\frac{I(\rho)}{\rho^3}+\frac{I'(\rho)}{\rho^2}-\frac{J'(\rho)}{\rho^2}\bigg|\\
    &\le \big(2\Lip(\Omega)+\tilde L)\frac{I'(\rho)}{\rho}+2\Lip(\Omega)\frac{I(\rho)}{\rho^2}\\
    &\quad+\ell\rho^{\gamma-1}\bigg(\frac{\H^2(\Sigma\cap B_{\rho})}{\rho^2}\bigg)\\
    &\le A\frac{I(\rho)}{\rho^2}+A\frac{d}{d\rho}\bigg(\frac{I(\rho)}{\rho}\bigg)+\ell\rho^{\gamma-1}\bigg(\frac{\H^2(\Sigma\cap B_{\rho})}{\rho^2}\bigg),
\end{align*}
where $A:=2\Lip(\Omega)+\tilde L$. From the previous estimate, we immediately conclude that
\begin{align}
\label{main estimate below}
 \nonumber
 \frac{d}{d\rho}\bigg(\frac{I(\rho)}{\rho^2}\bigg)+(A+\ell\rho^{\gamma-1})\frac{I(\rho)}{\rho^2}&\ge\frac{J'(\rho)}{\rho^2}-\frac{d}{d\rho}\bigg(A\rho\frac{I(\rho)}{\rho^2}\bigg)\\
 &\quad+\ell\rho^{\gamma-1}\bigg(\frac{I(\rho)}{\rho^2}-\frac{\H^2(\Sigma\cap B_{\rho})}{\rho^2}\bigg)
\end{align}
and 
\begin{align}
\label{main estimate above}
 \nonumber
 \frac{d}{d\rho}\bigg(\frac{I(\rho)}{\rho^2}\bigg)-(A+\ell\rho^{\gamma-1})\frac{I(\rho)}{\rho^2}&\le\frac{J'(\rho)}{\rho^2}+\frac{d}{d\rho}\bigg(A\rho\frac{I(\rho)}{\rho^2}\bigg)\\
 &\quad+\ell\rho^{\gamma-1}\bigg(\frac{\H^2(\Sigma\cap B_{\rho})}{\rho^2}-\frac{I(\rho)}{\rho^2}\bigg).
\end{align}
By letting $\varphi$ increase to the characteristic function of the interval $[0,1]$ in \eqref{main estimate below}, the above estimate passes to the limit in the sense of distributions and we obtain
\begin{align*}
     \frac{d}{d\rho}\bigg(\frac{\H^2(\Sigma\cap B_{\rho})}{\rho^2}\bigg)&+(A+\ell\rho^{\gamma-1})\frac{\H^2(\Sigma\cap B_{\rho})}{\rho^2}\\
     &\ge\frac{d}{d\rho}\bigg(\int_{\Sigma\cap B_{\rho}}\frac{\langle\Omega_t,\vec\Sigma_J\rangle}{|\cdot|^2}\, d\H^2\bigg)-\frac{d}{d\rho}\bigg(A\rho\frac{\H^2(\Sigma\cap B_{\rho})}{\rho^2}\bigg).
\end{align*}
Multiplying both of the last inequality sides by the factor $\displaystyle{e^{A\rho+\ell\frac{\rho^{\gamma}}{\gamma}}}$ and taking into account the fact that the first term on the right-hand-side is non-negative, we get
\begin{align*}
     \frac{d}{d\rho}\bigg(e^{A\rho+\ell\frac{\rho^{\gamma}}{\gamma}}\frac{\H^2(\Sigma\cap B_{\rho})}{\rho^2}\bigg)&\ge\frac{d}{d\rho}\bigg(\int_{\Sigma\cap B_{\rho}}\frac{\langle\Omega_t,\vec\Sigma_J\rangle}{|\cdot|^2}\, d\H^2\bigg)\\
     &\quad-\frac{d}{d\rho}\bigg(e^{A\rho+\ell\frac{\rho^{\gamma}}{\gamma}}A\rho\frac{\H^2(\Sigma\cap B_{\rho})}{\rho^2}\bigg),
\end{align*}
which turns into
\begin{align*}
    \frac{d}{d\rho}\bigg(e^{A\rho+\ell\frac{\rho^{\gamma}}{\gamma}}(1+A\rho)\frac{\H^2(\Sigma\cap B_{\rho})}{\rho^2}\bigg)&\ge\frac{d}{d\rho}\bigg(\int_{\Sigma\cap B_{\rho}}\frac{1}{|\cdot|^2}\langle\Omega_t,\vec\Sigma_J\rangle\, d\H^2\bigg).
\end{align*}
By integration of the previous inequality, we get
\begin{align*}
e^{A\rho+\ell\frac{\rho^{\gamma}}{\gamma}}(1+A\rho)\frac{\H^{2}(\Sigma\cap B_{\rho})}{\rho^2}&-e^{A\sigma+\ell\frac{\sigma^{\gamma}}{\gamma}}(1+A\sigma)\frac{\H^{2}(\Sigma\cap B_{\sigma})}{\sigma^2}\\
&\ge\int_{\Sigma\cap(B_{\rho}\smallsetminus B_{\sigma})}\frac{1}{|\cdot|^2}\langle\Omega_t,\vec\Sigma_J\rangle\, d\H^2,
\end{align*}
for every $0<\sigma<\rho<1$. Since
\begin{align*}
    \langle\Omega_t,\vec\Sigma_J\rangle=|\vec\Sigma_J\wedge\nu|_{g}^2,
\end{align*}
the estimate \eqref{monotonicity sigma below} follows.

By applying the same techniques to \eqref{main estimate above}, we get \eqref{monotonicity sigma above} and the statement follows. 
\end{proof}
\end{prop}
\begin{rem}
Proposition \ref{monotonicity formula sigma} immediately implies that the function
\begin{align*}
    (0,1)\ni\rho\mapsto e^{A\rho+\ell\frac{\rho^{\gamma}}{\gamma}}(1+A\rho)\frac{\H^2(\Sigma\cap B_{\rho})}{\rho^2}
\end{align*}
is non-decreasing. In particular, the limit
\begin{align*}
    \theta(0,\Sigma):=\lim_{\rho\rightarrow 0^+}\frac{\H^2(\Sigma\cap B_{\rho})}{\rho^2}=\lim_{\rho\rightarrow 0^+}e^{A\rho+\ell\frac{\rho^{\gamma}}{\gamma}}(1+A\rho)\frac{\H^2(\Sigma\cap B_{\rho})}{\rho^2}
\end{align*}
exists and it is finite.
\end{rem}
\begin{lem}
\label{lemma sigma zero}
Let $\Sigma$ be a closed almost $J$-holomorphic curve in $B$, according to Definition \ref{definition of almost pseudoholomorphic curve}. Then, there exist an $\H^2$-measurable field of unit and simple $2$-vectors $\vec\Sigma_0\in\bigwedge_2\r^{2m}$ on $\Sigma$ and a constant $L>0$ depending only on $\ell$ and on $\Lip(\Omega)$ such that for $\H^2$-a.e. $x\in\Sigma$ the following facts hold:
\begin{enumerate}
    \item $\vec\Sigma_0(x)$ is a unit simple $2$-vector calibrated by $\Omega_0$;
    \item $|\vec\Sigma(x)-\vec\Sigma_0(x)|\le L|x|^{\gamma/2}$;
    \item if $\vec\Sigma(x)$ is $J_0$-invariant, then $\vec\Sigma_0(x)=\vec\Sigma(x)$;
    \item if $\vec\Sigma(x)$ is not $J_0$-invariant, then
        \begin{align*}
            |\vec\Sigma_0(x)\wedge\nu_0(x)\wedge J_0\nu(x)|=\max_{v\in S_x}|v\wedge J_0v\wedge\nu_0(x)\wedge J_0\nu_0(x)|,
        \end{align*}
        where $S_x$ denotes the unit sphere in the approximate tangent space $T_x\Sigma$.
\end{enumerate}
\begin{proof}
Recall the definition of the vector field $\nu_0$, given at the beginning of the present section. If $x\in\Sigma$ is such that $\vec\Sigma(x)$ is $J_0$-invariant, we set $\vec\Sigma_0(x):=\vec\Sigma(x)$ and all the required properties are satisfied. 
Otherwise, if $x\in\Sigma$ is such that $\vec\Sigma(x)$ is not $J_0$-invariant, we first claim that there exists some $\Omega_0$-orthonormal basis  
\begin{align*}
    \big\{e_1(x),J_0e_1(x),...,e_m(x),J_0e_m(x)\big\}
\end{align*}
of $\r^{2m}$ such that 
\begin{align*}
    |e_1(x)\wedge J_0e_1(x)\wedge\nu_0(x)\wedge J_0\nu_0(x)|=\max_{v\in S_x}|v\wedge J_0v\wedge\nu_0(x)\wedge J_0\nu_0(x)|
\end{align*}
and we can write $\vec\Sigma(x)$ as 
\begin{align*}
    \vec\Sigma(x)=\cos\phi(x)e_1(x)\wedge J_0e_1(x)+\sin\phi(x)e_1(x)\wedge e_2(x).
\end{align*}
for some angle $\phi(x)\in[0,2\pi]$. Indeed, since $S_x$ is compact, there exists $e_1(x)\in S_x$ maximizing the continuous function $v\mapsto|v\wedge J_0v\wedge\nu_0(x)\wedge J_0\nu_0(x)|$. We complete $\{e_1(x)\}$ to an $\Omega_0$-orthonormal basis $\{e_1(x),\xi(x)\}$ of $T_x\Sigma$ and we write $\vec\Sigma(x)=e_1(x)\wedge\xi(x)$.

Notice that the set $\{e_1(x), J_0e_1(x),\xi(x)-J_0e_1(x)\}$ is linearly independent, otherwise $\vec\Sigma(x)$ would be $J_0$-invariant. We define $e_2(x)$ as the unique vector such that $\{e_1(x), J_0e_1(x),e_2(x)\}$ is an orthonormal set and
\begin{align*}
    \text{span}\{e_1(x), J_0e_1(x),e_2(x)\}=\text{span}\{e_1(x), J_0e_1(x), \xi(x)-J_0e_1(x)\}.
\end{align*} 
Eventually, notice that 
\begin{align*}
    \vec\Sigma(x)=e_1(x)\wedge \xi(x)\in\text{span}\big\{e_1(x)\wedge J_0e_1(x), e_1(x)\wedge e_2(x)\big\}
\end{align*} 
and our initial claim follows since 
\begin{align*}
    |e_1(x)\wedge J_0e_1(x)|=|e_1(x)\wedge e_2(x)|=1.
\end{align*}
Then, we define $\vec\Sigma_0(x):=e_1(x)\wedge J_0e_1(x)$. Clearly, $\vec\Sigma_0(x)$ satisfies (1) and (4). For what concerns (2), notice that $|\Omega_x-\Omega_0|\le\Lip(\Omega)|x|$. In particular, since $\langle\vec\Sigma_J,\Omega_x\rangle=1$ and $\lvert\vec\Sigma_J-\vec\Sigma\rvert\le\ell\lvert\,\cdot\,\rvert^{\gamma}$, it follows that $|\big<\Omega_0,\vec\Sigma(x)\big>-1|\le C|x|^{\gamma}$, for some constant $C>0$ depending only on $\ell$ and $\Lip(\Omega)$. Then,
\begin{align*}
    1+C|x|^{\gamma}&\ge\big<\Omega_0,\vec\Sigma(x)\big>\\
    &=\big<\Omega_0,\cos\phi(x)\,e_1(x)\wedge J_0e_1(x)+\sin\phi(x)\, e_1(x)\wedge e_2(x)\big>\\
    &=\cos\phi(x)\ge 1-C|x|^{\gamma}.
\end{align*}
Hence, we eventually obtain
\begin{align*}
    \big|\vec\Sigma(x)-\vec\Sigma_0(x)\big|^2=(1-\cos\phi(x))^2+\sin^2\phi(x)=2(1-\cos\phi)\le 2C|x|^{\gamma}
\end{align*}
and the statement follows with $L:=\sqrt{2C}$.
\end{proof}
\end{lem}
For the purposed of the following lemma, we recall that with $\pi:\c^m\smallsetminus\{0\}\rightarrow\cp^{m-1}$ we denote the standard projection to the quotient (see section 1.3). 
\begin{lem}
\label{lemma projected sigma}
Under the same hypothesis and notation of Lemma \ref{lemma sigma zero}, let $\tilde r\in(0,1)$. Then, there exists some constant $\Xi=\Xi\big(m,\tilde r,\Lip(\Omega),\ell,\gamma\big)>0$ such that
\begin{enumerate}
    \item for every $0<\rho\le\tilde r$ and every open set $U\subset B_{\rho}$ it holds that 
    \begin{align*}
        \bigg|\m\big(\pi_{*}([\Sigma]\res U)\big)-\int_{\Sigma\cap U}|\opwedge_2d\pi(\vec\Sigma_0)|\, d\H^2\bigg|\le\Xi\H^2(\Sigma\cap B_{\tilde r})\rho^{\gamma/2}.
    \end{align*}
    \item there exists constants $G>0$ (depending only on the metric) and $K_m>0$ (obtained as $C_{2m}$ in \cite[Lemma 1]{pumberger-riviere}) such that in for every $0<\rho\le\tilde r$, we have
    \begin{align*}
        \int_{\Sigma\cap B_{\rho}}|\opwedge_2d\pi(\vec\Sigma_0)|\, d\H^2&\le 2K_mG\bigg(e^{A\rho+\ell\frac{\rho^{\gamma}}{\gamma}}(1+A\rho)\frac{\H^2(\Sigma\cap B_{\rho})}{\rho^2}-\theta(0,\Sigma)\bigg)\\
        &\quad+\Xi\H^2(\Sigma\cap B_{\tilde r})\rho^{\gamma}
    \end{align*}
\end{enumerate}
\begin{proof}
First, we aim to prove (1). By Lemma \ref{lemma sigma zero} and by the definition of almost $J$-holomorphic curve, it follows that
\begin{align*}
    \big||\opwedge_2d\pi_x(\vec\Sigma(x))|-|\opwedge_2d\pi_x(\vec\Sigma_0(x))|\big|&\le|\opwedge_2d\pi_x(\vec\Sigma(x))-\opwedge_2d\pi_x(\vec\Sigma_0(x))|\\
    &\le\frac{1}{|x|^2}|\vec\Sigma(x)-\vec\Sigma_0(x)|\le\frac{L}{|x|^{2-\gamma/2}},
\end{align*}
for $\H^2$-a.e. every $x\in\Sigma\smallsetminus\{0\}$. By integrating on $U$ both sides in the previous inequality, we get 
\begin{align}
\label{estimate projected sigma}
    \nonumber
     \bigg|\m\big(\pi_{*}([\Sigma]\res U)\big)-\int_{\Sigma\cap U}|\opwedge_2d\pi(\vec\Sigma_0)|\, d\H^2\bigg|&\le L\int_{\Sigma\cap U}\frac{1}{|x|^{2-\gamma/2}}\, d\H^2\\
     &\le L\int_{\Sigma\cap B_{\rho}}\frac{1}{|x|^{2-\gamma/2}}\, d\H^2.
\end{align}
Notice that, by exploiting \eqref{monotonicity sigma below}, we get
\begin{align*}
    \int_{\Sigma\cap(B_{\rho}\smallsetminus B_{\rho/2})}\frac{1}{|x|^{2-\gamma/2}}\, d\H^2(x)&\le 2^{2-\gamma/2}\,\frac{\H^2(\Sigma\cap B_{\rho})}{\rho^{2-\gamma/2}}\\
    &\le 4e^{A\rho+\ell\frac{\rho^{\gamma}}{\gamma}}(1+A\rho)\frac{\H^2(\Sigma\cap B_{\rho})}{\rho^2}\frac{\rho^{\gamma/2}}{2^{\gamma/2}}\\
    &\le\frac{4e^{A+\frac{\ell}{\gamma}}(1+A)}{\tilde r^2}\H^2(\Sigma\cap B_{\tilde r})\frac{\rho^{\gamma/2}}{2^{\gamma/2}},
\end{align*}
for every $0<\rho\le\tilde r$. By iteration, we obtain
\begin{align*}
    &\int_{\Sigma\cap(B_{\rho}\smallsetminus B_{\rho/2^n})}\frac{1}{|x|^{2-\gamma/2}}\, d\H^2(x)\\
    &\le\frac{4e^{A+\frac{\ell}{\gamma}}(1+A)}{\tilde r^2}\H^2(\Sigma\cap B_{\tilde r})\Bigg(\sum_{j=0}^{n-1}\frac{1}{(2^{\gamma/2})^j}\Bigg)\rho^{\gamma/2}
\end{align*}
and, by passing to the limit as $n\rightarrow+\infty$, we have
\begin{align}
\label{useful estimate projected sigma 1}
    \int_{\Sigma\cap B_{\rho}}\frac{1}{|x|^{2-\gamma/2}}\, d\H^2(x)\le\frac{4e^{A+\frac{\ell}{\gamma}}(1+A)}{\tilde r^2(1-2^{-\gamma/2})}\H^2(\Sigma\cap B_{\tilde r})\rho^{\gamma/2},
\end{align}
for every $0<\rho\le\tilde r$. By combining the previous estimate with \eqref{estimate projected sigma}, estimate (1) follows with
\begin{align*}
    \Xi_1:=\frac{4e^{A+\frac{\ell}{\gamma}}(1+A)L}{\tilde r^2(1-2^{-\gamma/2})}.
\end{align*}
For what concerns (2), we simply notice that by \eqref{monotonicity sigma below}, by point (2) in Lemma \ref{lemma sigma zero} and by \cite[Lemma 1]{pumberger-riviere}, we get
\begin{align*}
    \int_{\Sigma\cap B_{\rho}}|\opwedge_2d\pi(\vec\Sigma_0)|\, d\H^2&\le K_m\int_{\Sigma\cap B_{\rho}}\frac{1}{|\cdot|^2}|\vec\Sigma_0\wedge\nu_0|^2\, d\H^2\\
    &\le 4K_m\int_{\Sigma\cap B_{\rho}}\frac{1}{|\cdot|^2}|\vec\Sigma_J\wedge\nu|^2\, d\H^2\\
    &\quad+4K_m\int_{\Sigma\cap B_{\rho}}\frac{1}{|\cdot|^2}|\nu-\nu_0|^2\, d\H^2\\
    &\quad+4K_m\int_{\Sigma\cap B_{\rho}}\frac{1}{|\cdot|^2}|(\vec\Sigma_J-\vec\Sigma_0)\wedge\nu_0|^2\, d\H^2\\
    &\le 4K_mG\bigg(e^{A\rho+\ell\rho^{\gamma}}(1+A\rho)\frac{\H^{2}(\Sigma\cap B_{\rho})}{\rho^2}-\theta(0,\Sigma)\bigg)\\
    &\quad+4K_m\tilde L\H^2(\Sigma\cap B_{\rho})\\
    &\quad+4K_m(\ell^2+L^2)\int_{\Sigma\cap B_{\rho}}\frac{1}{|\cdot|^{2-\gamma}}\, d\H^2.
\end{align*}
By using the same method that we have used in order to prove the decay in \eqref{useful estimate projected sigma 1}, we can show that
\begin{align}
\label{useful estimate projected sigma 2}
    \int_{\Sigma\cap B_{\rho}}\frac{1}{|x|^{2-\gamma}}\, d\H^2(x)\le\frac{4e^{A+\frac{\ell}{\gamma}}(1+A)}{\tilde r^2}\H^2(\Sigma\cap B_{\tilde r})\rho^{\gamma},
\end{align}
for very $0<\rho\le\tilde r$. Moreover, we clearly have that 
\begin{align*}
    \frac{\H^2(\Sigma\cap B_{\rho})}{\rho^{\gamma}}&\le\frac{\H^2(\Sigma\cap B_{\rho})}{\rho^2}\le\frac{4e^{A\rho+\ell\frac{\rho^{\gamma}}{\gamma}}(1+A\rho)}{\rho^2}\H^2(\Sigma\cap B_{\rho})\\
    &\le \frac{4e^{A+\frac{\ell}{\gamma}}(1+A)}{\tilde r^2}\H^2(\Sigma\cap B_{\tilde r}),
\end{align*}
for very $0<\rho\le\tilde r$. Thus, we get that (2) holds with
\begin{align*}
    \Xi_2:=\frac{4K_me^{A+\frac{\ell}{\gamma}}(1+A)(\tilde L+\ell^2+L^2)}{\tilde r^2}.
\end{align*}
Hence, the statement follows with $\Xi:=\max\{\Xi_1,\Xi_2\}$.
\end{proof}
\end{lem}
\begin{rem}
A first remarkable consequence of Lemma \ref{lemma projected sigma} and Proposition \ref{monotonicity formula sigma} is that
\begin{align}
    \m\big(\pi_{*}([\Sigma]\res B_{\rho})\big)\rightarrow 0 \qquad \mbox{ as } \rho\rightarrow 0^+.
\end{align}
\end{rem}
\begin{lem}[Good slicing]
\label{lemma good slices}
Under the same hypotheses and notation of Lemma \ref{lemma sigma zero}, let $\tilde r\in(0,1)$. Then, for every $r\in(0,\tilde r]$ there exist $\tilde\rho\in[r/2,r]$ and a constant $\Theta=\Theta\big(m,\tilde r,\Lip(\Omega),\ell,\gamma\big)>0$ such that:
\begin{enumerate}
\item $\H^1(\Sigma\cap\partial B_{\tilde\rho})\le\Theta\H^2(\Sigma\cap B_{\tilde r})\tilde\rho$;
\item 
    $\begin{aligned}
        \int_{\Sigma\cap\partial B_{\tilde\rho}}|\opwedge_2d\pi(\vec\Sigma)|\, d\H^1\le\frac{\Theta}{\tilde\rho}\int_{\Sigma\cap(B_r\smallsetminus B_{r/2})}|\opwedge_2d\pi(\vec\Sigma)|\, d\H^2;
    \end{aligned}$
\item 
    $\begin{aligned}
        \m\big(\pi_{\ast}\partial([\Sigma]\res B_{\tilde\rho})\big)\le\Theta\sqrt{K_m}\fint_{r/2}^r\frac{1}{\rho}\int_{\Sigma\cap\partial B_{\rho}}|\vec\Sigma_0\wedge\nu_0|\, d\H^1\, d\L^1(\rho).
    \end{aligned}$
\end{enumerate}
\begin{proof}
First, we notice that, by the coarea formula and the monotonicity formula \eqref{monotonicity sigma below}, it holds that
\begin{align*}
    \int_{r/2}^r\frac{\H^1(\Sigma\cap\partial B_{\rho})}{\rho}\, d\L^1(\rho)&\le\frac{2}{r}\H^2(\Sigma\cap B_r)\\
    &\le 2e^{Ar+\ell\frac{r^{\gamma}}{\gamma}}(1+Ar)\frac{\H^2(\Sigma\cap B_r)}{r^2}r\\
    &\le\frac{4e^{A+\frac{\ell}{\gamma}}(1+A)}{\tilde r^2}\H^2(\Sigma\cap B_{\tilde r})\frac{r}{2}.
\end{align*}
Hence,
\begin{align}
\label{good slice estimate 1}
\fint_{r/2}^r\frac{1}{\Theta_1\H^2(\Sigma\cap B_{\tilde r})}\frac{\H^1(\Sigma\cap\partial B_{\rho})}{\rho}\, d\L^1(\rho)\le 1,
\end{align}
with 
\begin{align*}
    \Theta_1:=\frac{4e^{A+\frac{\ell}{\gamma}}(1+A)}{\tilde r^2}.
\end{align*}
Moreover, again by the coarea formula, we get
\begin{align*}
    &\fint_{r/2}^r\rho\int_{\Sigma\cap\partial B_{\rho}}|\opwedge_2d\pi(\vec\Sigma)|\, d\H^1\, d\L^1(\rho)\\
    &\le 2\int_{r/2}^r\int_{\Sigma\cap\partial B_{\rho}}|\opwedge_2d\pi(\vec\Sigma)|\, d\H^1\, d\L^1(\rho)\\
    &=2\int_{\Sigma\cap(B_r\smallsetminus B_{r/2})}|\opwedge_2d\pi(\vec\Sigma)|\, d\H^2=:a,
\end{align*}
which leads to 
\begin{align}
\label{good slice estimate 2}
\fint_{r/2}^r\frac{1}{a}\rho\int_{\Sigma\cap\partial B_{\rho}}|\opwedge_2d\pi(\vec\Sigma)|\, d\H^1\, d\L^1(\rho)&\le 1.
\end{align}
Eventually, by \cite[Lemma 7.6.1]{krantzparks}, we know that for a.e. $\rho\in(0,1)$ the slice $\Sigma\cap \partial B_{\rho}$ is a $1$-rectifiable subset of $B$ and the vector field $\vec\Sigma_{\rho}$ orienting its approximate tangent space at $x$ belongs to $S_x$ (see notation of Lemma \ref{lemma sigma zero}). Then, by \cite[Lemma 1]{pumberger-riviere} and by points (3) and (4) of Lemma \ref{lemma sigma zero}, it follows that
\begin{align*}
    \big|\vec\Sigma_{\rho}\wedge\nu_0\wedge J_0\nu_0\big|^2&=\big|\vec\Sigma_{\rho}\wedge J_0\vec\Sigma_{\rho}\wedge\nu_0\wedge J_0\nu_0\big|\\
    &\le|\vec\Sigma_0\wedge\nu_0\wedge J_0\nu_0\big|\\
    &\le K_m\big|\vec\Sigma_0\wedge\nu_0\big|^2, \qquad\H^1\mbox{-a.e. on } \Sigma\cap\partial B_{\rho}.
\end{align*}
Hence, we get
\begin{align*}
    \m\big(\pi_{\ast}\partial([\Sigma]\res B_{\rho})\big)&\le\int_{\Sigma\cap\partial B_{\rho}}\big|d\pi(\vec\Sigma_{\rho})\big|\, d\H^1\\
    &=\frac{1}{\rho}\int_{\Sigma\cap\partial B_{\rho}}|\vec\Sigma_{\rho}\wedge\nu_0\wedge J\nu_0\big|\, d\H^1\\
    &\le\frac{\sqrt{K_m}}{\rho}\int_{\Sigma\cap\partial B_{\rho}}\big|\vec\Sigma_0\wedge\nu_0\big|\, d\H^1.
\end{align*}
Thus, by averaging the previous inequality on $[r/2,r]$, we obtain
\begin{align*}
    \fint_{r/2}^r\m\big(\pi_{\ast}\partial([\Sigma]\res B_{\rho})\, d\L^1(\rho)\le\sqrt{K_m}\fint_{r/2}^r\frac{1}{\rho}\int_{\Sigma\cap\partial B_{\rho}}\big|\vec\Sigma_0\wedge\nu_0\big|\, d\H^1=:b,
\end{align*}
which leads to 
\begin{align}
\label{good slice estimate 3}
    \fint_{r/2}^r\frac{1}{b}\m\big(\pi_{\ast}\partial([\Sigma]\res B_{\rho})\, d\L^1(\rho)\le 1.
\end{align}
By summing up the three inequalities \eqref{good slice estimate 1}, \eqref{good slice estimate 2} and \eqref{good slice estimate 3} we obtain
\begin{align*}
    &\fint_{r/2}^r\bigg(\frac{1}{\Theta_1\H^2(\Sigma\cap B_{\tilde r})}\frac{\H^1(\Sigma\cap\partial B_{\rho})}{\rho}+\frac{1}{a}\rho\int_{\Sigma\cap\partial B_{\rho}}|\opwedge_2d\pi(\vec\Sigma)|\, d\H^1\\
    &+\frac{1}{b}\m\big(\pi_{\ast}\partial([\Sigma]\res B_{\rho})\bigg)\, d\L^1(\rho)\le 3.
\end{align*}
Then, we conclude that there exists $\rho\in[r/2,r]$ such that
\begin{align*}
    &\frac{1}{\Theta_1\H^2(\Sigma\cap B_{\tilde r})}\frac{\H^1(\Sigma\cap\partial B_{\rho})}{\rho}+\frac{1}{a}\rho\int_{\Sigma\cap\partial B_{\rho}}|\opwedge_2d\pi(\vec\Sigma)|\, d\H^1\\
    &+\frac{1}{b}\m\big(\pi_{\ast}\partial([\Sigma]\res B_{\rho})\le 3
\end{align*}
and the statement follows with $\Theta:=\max\{\Theta_1,6\}$. 
\end{proof}
\end{lem}
\begin{lem}[Controlling the mass of the projected boundaries]
\label{lemma projected boundaries}
Under the same hypotheses and notation of Lemma \ref{lemma sigma zero}, let $\tilde r\in(0,1)$. Let $r\in(0,\tilde r]$ be such that
\begin{align}
\label{estimate projected boundaries 1}
    \m\big(\pi_{\ast}([\Sigma]\res B_r)\big)<2K_m^2\m\big(\pi_{\ast}([\Sigma]\res B_{r/2})\big)
\end{align}
and 
\begin{align}
\label{estimate projected boundaries 2}
    \int_{\Sigma\cap B_r}|\opwedge_2d\pi(\vec\Sigma_0)|\, d\H^2>\zeta^{-1}\Xi\H^2(\Sigma\cap B_{\tilde r})r^{\gamma/2},
\end{align}
for some $\zeta\in(0,1)$. If $\rho_j\in[r/2,r]$ is such that $\Sigma\cap\partial B_{\rho_j}$ is a good slice of $\Sigma$ in the sense of Lemma \ref{lemma good slices}, then
\begin{align}
    \nonumber
    \m\big(\pi_{\ast}\partial([\Sigma]\res B_{\rho_j})\big)&\le\Lambda\sqrt{\H^2(\Sigma\cap B_{\tilde r})}\sqrt{\m\big(\pi_{\ast}([\Sigma]\res B_{\rho_j})\big)}\\
    &\quad+\Lambda\H^2(\Sigma\cap B_{\tilde r})\rho_j^{\gamma/4},
\end{align}
for a constant $\Lambda=\Lambda\big(m,\tilde r,\Lip(\Omega),\ell,\gamma\big)>0$.
\begin{proof}
Let $\rho_j\in[r/2,r]$ be such that $\Sigma\cap\partial B_{\rho_j}$ is a good slice of $\Sigma$. We apply twice the Cauchy-Schwarz inequality in the right-hand side of (3) in Lemma \ref{lemma good slices} and the coarea formula to get
\allowdisplaybreaks
\begin{align*}
    \m\big(\pi_{\ast}\partial([\Sigma]\res B_{\rho_j})\big)&\le\Theta\sqrt{K_m}\fint_{r/2}^r\frac{1}{\rho}\int_{\Sigma\cap\partial B_{\rho}}|\vec\Sigma_0\wedge\nu_0|\, d\H^1\, d\L^1(\rho)\\
    &=\Theta\sqrt{K_m}\fint_{r/2}^r\int_{\Sigma\cap\partial B_{\rho}}\frac{1}{\rho}|\vec\Sigma_0\wedge\nu_0|\, d\H^1\, d\L^1(\rho)\\
    &\le\Theta\sqrt{K_m}\fint_{r/2}^r\sqrt{\H^1(\Sigma\cap\partial B_{\rho})}\\
    &\quad\cdot\sqrt{\int_{\Sigma\cap\partial B_{\rho}}\frac{1}{\rho^2}|\vec\Sigma_0\wedge\nu_0|^2\, d\H^1}\, d\L^1(\rho)\\
    &=\Theta\sqrt{K_m}\frac{2}{r}\int_{r/2}^r\sqrt{\H^1(\Sigma\cap\partial B_{\rho})}\\
    &\quad\cdot\sqrt{\int_{\Sigma\cap\partial B_{\rho}}\frac{1}{\rho^2}|\vec\Sigma_0\wedge\nu_0|^2\, d\H^1}\, d\L^1(\rho)\\
    &\le \Theta\sqrt{K_m}\frac{2}{r}\sqrt{\int_{r/2}^r\H^1(\Sigma\cap\partial B_{\rho})\, d\L^1(\rho)}\\
    &\quad\cdot\sqrt{\int_{r/2}^r\int_{\Sigma\cap\partial B_{\rho}}\frac{1}{\rho^2}|\vec\Sigma_0\wedge\nu_0|^2\, d\H^1\, d\L^1(\rho)}\\
    &=\Theta\sqrt{K_m}\sqrt{\frac{2}{r}}\sqrt{\fint_{r/2}^r\H^1(\Sigma\cap\partial B_{\rho})\, d\L^1(\rho)}\\
    &\quad\cdot\sqrt{\int_{\Sigma\cap(B_r\smallsetminus B_{r/2})}\frac{1}{|\cdot|^2}|\vec\Sigma_0\wedge\nu_0|^2\, d\H^2}.\\
\end{align*}
We notice that, by point (1) in Lemma \ref{lemma projected sigma} and by our assumption \eqref{estimate projected boundaries 2}, it holds that
\begin{align*}
    \bigg|\m\big(\pi_{*}([\Sigma]\res B_r)\big)-\int_{\Sigma\cap B_r}|\opwedge_2d\pi(\vec\Sigma_0)|\, d\H^2\bigg|&\le\Xi\H^2(\Sigma\cap B_{\tilde r})r^{\gamma/2}\\
    &<\zeta\int_{\Sigma\cap B_r}|\opwedge_2d\pi(\vec\Sigma_0)|\, d\H^2
\end{align*}
which implies
\begin{align*}
    \int_{\Sigma\cap B_r}|\opwedge_2d\pi(\vec\Sigma_0)|\, d\H^2\le\frac{1}{1-\zeta}\m\big(\pi_{*}([\Sigma]\res B_r)\big)
\end{align*}
Hence, \cite[Lemma 1]{pumberger-riviere} we have 
\begin{align*}
    \int_{\Sigma\cap(B_r\smallsetminus B_{r/2})}\frac{1}{|\cdot|^2}|\vec\Sigma_0\wedge\nu_0|^2\, d\H^2&\le K_m\int_{\Sigma\cap(B_r\smallsetminus B_{r/2})}|\opwedge_2d\pi(\vec\Sigma_0)|\, d\H^2\\
    &\le\frac{K_m}{1-\zeta}\m\big(\pi_{*}([\Sigma]\res B_{r})\big).
\end{align*}
Moreover, by \eqref{good slice estimate 1}, it follows that
\begin{align*}
    \sqrt{\fint_{r/2}^r\H^1(\Sigma\cap\partial B_{\rho})\, d\L^1(\rho)}&\le\sqrt{r}\sqrt{\fint_{r/2}^r\frac{\H^1(\Sigma\cap\partial B_{\rho})}{\rho}\, d\L^1(\rho)}\\
    &\le\sqrt{2\Theta}\sqrt{\H^2(\Sigma\cap B_{\tilde r})}\sqrt{\frac{r}{2}}.
\end{align*}
Thus, 
\begin{align*}
   \m\big(\pi_{\ast}\partial([\Sigma]\res B_{\rho_j})\big)\le\frac{\sqrt{2}\Theta^{3/2}K_m}{\sqrt{1-\zeta}}\sqrt{\H^2(\Sigma\cap B_{\tilde r})}\sqrt{\m\big(\pi_{\ast}([\Sigma]\res B_r)\big)}. 
\end{align*}
By our hypothesis \eqref{estimate projected boundaries 1} and since $\rho_j>r/2$, we obtain that 
\begin{align*}
    \m\big(\pi_{\ast}([\Sigma]\res B_r)\big)&<2K_m^2\,\m\big(\pi_{\ast}([\Sigma]\res B_{r/2})\big)\\
    &\le 2K_m^2\,\m\big(\pi_{\ast}([\Sigma]\res B_{\rho_j})\big)+4K_m^2\Xi\H^2(\Sigma\cap B_{\tilde r})\rho_j^{\gamma/2}.
\end{align*}
We point out that the last inequality follows a direct application of point (1) in Lemma \ref{lemma projected sigma}. Then, the statement follows with
\begin{align*}
    \Lambda:=\max\bigg\{\frac{2\Theta^{3/2}K_m^2}{\sqrt{1-\zeta}},\frac{2\sqrt{2}\Theta^2K_m^2}{\sqrt{1-\zeta}}\bigg\}.
\end{align*}
\end{proof}
\end{lem}
We recall the following general fact about integral $2$-currents on $\cp^{m-1}$ with small mass which are $\zeta$-almost semicalibrated by $\omega_{\cp^{m-1}}$, whose proof can be found in \cite[Lemma 11]{pumberger-riviere}. Recall that a current $T\in\D^2(\cp^{m-1})$ is said to be $\zeta$-almost semicalibrated by $\omega_{\cp^{m-1}}$ for some constant $\zeta\in(0,1)$ if 
\begin{align*}
    (1-\zeta)|\left<T\res U,\omega_{\cp^{m-1}}\right>|\le\m(T\res U)\le(1+\zeta)|\left<T\res U,\omega_{\cp^{m-1}}\right>|,
\end{align*}
for every open set $U\subset\cp^{m-1}$.
\begin{lem}
\label{lemma tubular neighbourhood}
Let $\zeta\in(0,1)$. Given any couple of constants $\tilde\Lambda>0$ and $\lambda>0$, there exist $\delta>0$ and $\varepsilon>0$ satisfying what follows. For every integral $2$-current $T\in\D^2(\cp^{m-1})$ such that
\begin{enumerate}
\item $T$ is $\zeta$-almost semicalibrated by $\omega_{\cp^{m-1}}$,
\item $\m(T)+\m(\partial T)<\delta$,
\item $\m(\partial T)\le\tilde\Lambda\sqrt{\m(T)}$,
\end{enumerate}
there is a complex projective $(m-2)$-hyperplane $H\subset\cp^{m-1}$ and a tubular neighbourhood $H_{\varepsilon}\subset\cp^{m-1}$ of $H$ with width $\varepsilon$ such that
\begin{align}
\frac{\m(T\res H_{\varepsilon})}{\varepsilon^2}\le\lambda\m(T).
\end{align}
\end{lem}
As a last tool, we need to establish that given any complex projective $(m-2)$-hyperplane $H\subset\cp^{m-1}$ and a tubular neighbourhood $H_{\varepsilon}\subset\cp^{m-1}$ of $H$ with width $\varepsilon$, we can approximate the symplectic form $\omega_{\cp^{m-1}}$ on $\cp^{m-1}$ with an exact form $d\alpha$ that coincides with $\omega_{\cp^{m-1}}$ on the complement of $H_{\varepsilon}$ and vanishes on $H_{\varepsilon/2}$. We achieve this approximation through the following lemma, whose proof is again in \cite[Lemma 6]{pumberger-riviere}.
\begin{lem}
\label{lemma exact approximation}
Let $H\subset\cp^{m-1}$ be any complex projective $(m-2)$-hyperplane and let $H_{\varepsilon}\subset\cp^{m-1}$ be a tubular neighbourhood of $H$ with width $\varepsilon$. Then there exists a $1$-form $\alpha\in\Omega^1(\cp^{m-1})$ and a universal constant $\kappa>0$ such that:
\begin{enumerate}
\item $\omega_{\cp^{m-1}}=d\alpha$ on $\cp^{m-1}\smallsetminus H_{\varepsilon}$;
\item $\alpha=0$ on $H_{\varepsilon/2}$;
\item $||\alpha||_{\ast}\le\kappa$;
\item $\displaystyle{||\omega_{\cp^{m-1}}-d\alpha||_{\ast}\le\frac{\kappa}{\varepsilon^2}}$.
\end{enumerate}
\end{lem} 

\subsection{Proof of the fundamental Morrey type estimate}
Recall that $\kappa,\Xi>0$ are positive constants introduced in Lemma \ref{lemma exact approximation} and Lemma \ref{lemma projected sigma} respectively.

Fix any $j_0\in\n\smallsetminus\{0\}$ and let $\ell\ge0$ be a constant depending only on $\Lip(\Omega)$. Let $\delta>0$ and $\varepsilon>0$ be the constants given by applying Lemma \ref{lemma tubular neighbourhood} with $\tilde\Lambda=\Lambda\big(m,2^{-j_0},\Lip(\Omega),\ell,1/2\big)$ from Lemma \ref{lemma projected boundaries} and $\lambda:=5(24\kappa)^{-1}$. Let $\delta'>0$ be such that
\begin{align*}
    \Lambda\sqrt{\delta'}+\frac{\delta'}{2}+\delta'<\delta
\end{align*}
and choose $\tilde r\in\big(0,\min\{2^{-j_0},\delta'(2\Xi)^{-1}\}\big)$ such that 
\begin{align*}
   \max\{\Lambda,\Xi\}\big(e^{A+2\ell}(1+A)\big)^{-1}\tilde r^{1/4}<\frac{\delta'}{2}.
\end{align*}

Assume that $\F$ is a family of closed almost $J$-holomorphic curves in $B$ such that every element $\Sigma\in\F$ satisfies the following properties:
\begin{enumerate}
    \item $|\vec\Sigma(x)-\vec\Sigma_J(x)|\le\ell|x|^{1/2}$, for $\H^2$-a.e. $x\in\Sigma$,
    \item $\H^2(\Sigma\cap B_{2^{-j_0}})<\big(e^{A+2\ell}(1+A)\big)^{-1}$,
    \item   $\begin{aligned}
                \int_{\Sigma\cap B_{2^{-j_0}}}|\opwedge_2d\pi(\vec\Sigma_0)|\, d\H^2<\frac{\delta'}{2}.
            \end{aligned}$
\end{enumerate}
\begin{rem}
\label{remark good slices satisfy the hypotheses}
For every $\Sigma\in\F$, the hypotheses (2) and (3) combined with point (1) in Lemma \ref{lemma projected sigma} imply that:
\begin{enumerate}
    \item   $\begin{aligned}
                \m\big(\pi_{*}([\Sigma]\res B_{\rho})\big)+\m\big(\pi_{*}\partial([\Sigma]\res B_{\rho})\big)<\delta,
            \end{aligned}$
    \item   $\begin{aligned}
                \m\big(\pi_{*}\partial([\Sigma]\res B_{\rho})\big)\le\Lambda\sqrt{\m\big(\pi_{*}([\Sigma]\res B_{\rho})\big)},
            \end{aligned}$
\end{enumerate}
for every good slice $\Sigma\cap\partial B_{\rho}$ of $\Sigma$, where $\rho\in[r/2,r]$ with $r\in(0,\tilde r]$ satisfying the hypotheses \eqref{estimate projected boundaries 1} and \eqref{estimate projected boundaries 2} of Lemma \ref{lemma projected boundaries}.
\end{rem}
We want to show that for every $\Sigma\in\F$ there exist constants $C>0$ and $0<\alpha<1$ depending on $m$, $j_0$, $\Lip(\Omega)$ such that
\begin{align}
\label{morrey decrease}
    \bigg|\int_{\Sigma\cap B_{\rho}}\pi^*\omega_{\cp^{m-1}}|_{\Sigma}\bigg|\le C\rho^{\alpha}, \qquad\forall\,\rho\in(0,2^{-j_0}).
\end{align}
By definition of mass it holds that
\begin{align}
\label{estimate with the mass from above}
\nonumber
    \bigg|\int_{\Sigma\cap B_{\rho}}\pi^*\omega_{\cp^{m-1}}|_{\Sigma}\bigg|&=|\left<\pi_{*}([\Sigma]\res B_{\rho}),\omega_{\cp^{m-1}}\right>|\\
    &\le\m\big(\pi_{*}([\Sigma]\res B_{\rho})\big),
\end{align}
for every $\rho\in (0,1)$. Hence, in order to prove \eqref{morrey decrease} it is enough to show that 
\begin{align}
\label{morrey decrease true}
    \m\big(\pi_{*}([\Sigma]\res B_{\rho})\big)\le C\rho^{\alpha}, \qquad \mbox{ for every } \rho\in(0,\tilde r).
\end{align}
Moreover, by exploiting point (1) in Lemma \ref{lemma projected sigma}, we realize that we if we show
\begin{align}
\label{morrey decrease very true}
     \int_{\Sigma\cap B_{\rho}}|\opwedge_2d\pi(\vec\Sigma_0)|\, d\H^2\le\tilde C\rho^{\tilde\alpha}, \qquad \mbox{ for every } \rho\in(0,\tilde r),
\end{align}
then \eqref{morrey decrease true} will follow with $C:=\tilde C+\Xi$ and $\alpha:=\min\{\tilde\alpha,1/4\}$. Thus, we just need to show \eqref{morrey decrease very true}. 

Fix any $\Sigma\in\F$. Let 
\begin{align*}
    E(\rho):=\int_{\Sigma\cap B_{\rho}}|\opwedge_2d\pi(\vec\Sigma_0)|\, d\H^2, \qquad\forall\,\rho\in(0,1)
\end{align*}
and define
\begin{align*}
    I&:=\bigg\{j\in\n \mbox{ s.t. } 2^{-j}\le\tilde r \mbox{ and } \, E\big(2^{-(j+1)}\big)\le\frac{1}{2}E\big(2^{-j}\big)\bigg\},\\
    J&:=\big\{j\in\n \mbox{ s.t. } 2^{-j}\le\tilde r \mbox{ and } \, E\big(2^{-(j+1)}\big)\le 5\Xi2^{-(j+1)/4}\big\}.
\end{align*}
First, we claim that there exists $\theta=\theta\big(m,j_0,\Lip(\Omega)\big)\in (0,1)$ such that 
\begin{align}
\label{step theta}
    E\big(2^{-(j+1)}\big)\le\theta\Big(E\big(2^{-j}\big)+2^{-j/4}\Big),
\end{align}
for every $j\in (I\cup J)^c$. Fix $j\in (I\cup J)^c$ and set $r:=2^{-j}$. 
Pick a radius $\rho_j\in[r/2,r]$ such that $\Sigma\cap\partial B_{\rho_j}$ is a good slice of $\Sigma$ (see Lemma \ref{lemma good slices}). By Lemma \ref{lemma projected sigma} and Lemma \ref{lemma projected boundaries}, it follows that we can choose a sequence of radii $\{s_k\}_{k\in\n}\in (0,r/2)$ such that $s_k\rightarrow 0^+$ as $k\rightarrow +\infty$ and
\begin{align}
\label{good choice of s_k}
    \m\big(\pi_{\ast}\partial([\Sigma]\res B_{s_k})\big)&\le\m\big(\pi_{\ast}([\Sigma]\res B_r\smallsetminus B_{r/2})\big), \qquad\forall\, k\in\n.
\end{align}
Since $\Sigma\in\F$, by Remark \ref{remark good slices satisfy the hypotheses} and since $j\in (I\cup J)^c$ it follows that the current $T=\pi_{\ast}([\Sigma]\res B_{\rho_j})$ satisfies the hypotheses of Lemma \ref{lemma tubular neighbourhood} with $\zeta=1/2$. Thus, there exists a complex projective $(m-2)$-hyperplane $H\subset\cp^{m-1}$ and a tubular neighbourhood $H_{\varepsilon}\subset\cp^{m-1}$ of $H$ with width $\varepsilon>0$ such that 
\begin{align*}
    \frac{\m\big(\pi_{\ast}([\Sigma]\res B_{\rho_j})\res H_{\varepsilon}\big)}{\varepsilon^2}\le\lambda\m\big(\pi_{\ast}([\Sigma]\res B_{\rho_j})\big).
\end{align*}
We let $\alpha\in\Omega^1(\cp^{m-1})$ be a smooth $1$-form given by Lemma \ref{lemma exact approximation} relatively to $H,H_{\varepsilon}$. Following the proof of point (1) in Lemma \ref{lemma projected sigma}, we notice that
\begin{align}
\label{estimate morrey decrease 1}
    \nonumber
    &\int_{\Sigma\cap (B_{r/2}\smallsetminus B_{s_k})}|\opwedge_2d\pi(\vec\Sigma_0)|\, d\H^2\\
    \nonumber
    &\le\int_{\Sigma\cap (B_{\rho_j}\smallsetminus B_{s_k})}|\opwedge_2d\pi(\vec\Sigma_0)|\, d\H^2\\
    \nonumber
    &=2\bigg|\int_{\Sigma\cap (B_{\rho_j}\smallsetminus B_{s_k})}\left<\pi^*\omega_{\cp^{m-1}},\vec\Sigma\right>\, d\H^2\bigg|\\
    \nonumber
    &=2\bigg|\int_{\Sigma\cap (B_{\rho_j}\smallsetminus B_{s_k})}\left<\pi^*d\alpha,\vec\Sigma\right>\, d\H^2\bigg|\\
    &\quad+2\bigg|\int_{\Sigma\cap (B_{\rho_j}\smallsetminus B_{s_k})}\left<\pi^*(\omega_{\cp^{m-1}}-d\alpha),\vec\Sigma\right>\, d\H^2\bigg|.
\end{align}
For what concerns the second term in the last sum, by Lemmas \ref{lemma tubular neighbourhood}, \ref{lemma exact approximation} and (1) in Lemma \ref{lemma projected sigma} we see that
\begin{align}
\label{estimate morrey decrease 2}
    \nonumber
    &\bigg|\int_{\Sigma\cap (B_{\rho_j}\smallsetminus B_{s_k})}\left<\pi^*(\omega_{\cp^{m-1}}-d\alpha),\vec\Sigma\right>\, d\H^2\bigg|\\
    \nonumber
    &\le ||\omega_{\cp^{m-1}}-d\alpha||_{\ast}\m\big(\pi_{\ast}([\Sigma]\res B_{\rho_j})\res H_{\varepsilon}\big)\\
    \nonumber
    &\le\kappa\,\frac{\m\big(\pi_{\ast}([\Sigma]\res B_{\rho_j})\res H_{\varepsilon}\big)}{\varepsilon^2}\\
    \nonumber
    &\le\kappa\lambda\,\m\big(\pi_{\ast}([\Sigma]\res B_{\rho_j})\big)\\
    \nonumber
    &\le\frac{1}{4}\int_{\Sigma\cap B_{\rho_j}}|\opwedge_2d\pi(\vec\Sigma_0)|\, d\H^2+\frac{\Xi}{4}\H^2(\Sigma\cap B_{\rho_j})\rho_j^{1/4}\\
    &\le\frac{1}{4}\int_{\Sigma\cap B_r}|\opwedge_2d\pi(\vec\Sigma_0)|\, d\H^2+\tilde\Xi r^{1/4},
\end{align}
with
\begin{align*}
    \tilde\Xi:=\frac{\Xi}{4}\big(e^{A+2\ell}(1+A)\big)^{-1}.
\end{align*}
Now we want to estimate 
\begin{align*}
    \int_{\Sigma\cap (B_{\rho_j}\smallsetminus B_{s_k})}\left<\pi^*d\alpha,\vec\Sigma\right>\, d\H^2.
\end{align*}
Since $\pi$ is a smooth map on $B_{\rho_j}\smallsetminus B_{s_k}$, by Stokes theorem we get that
\begin{align*}
    &\bigg|\int_{\Sigma\cap (B_{\rho_j}\smallsetminus B_{s_k})}\left<\pi^*d\alpha,\vec\Sigma\right>\, d\H^2\bigg|\\
    &=\bigg|\int_{\Sigma\cap (B_{\rho_j}\smallsetminus B_{s_k})}\pi^*d\alpha|_\Sigma\bigg|=\bigg|\int_{\Sigma\cap (B_{\rho_j}\smallsetminus B_{s_k})}d(\pi^*\alpha)|_\Sigma\bigg|\\
    &=\bigg|\int_{\Sigma\cap\partial B_{\rho_j}}\pi^*\alpha|_{\Sigma\cap\partial B_{\rho_j}}-\int_{\Sigma\cap\partial B_{s_k}}\pi^*\alpha|_{\Sigma\cap\partial B_{s_k}}\bigg|\\
    &\le\left|\int_{\Sigma\cap\partial B_{\rho_j}}\pi^*\alpha|_{\Sigma\cap\partial B_{\rho_j}}\right|+\left|\int_{\Sigma\cap\partial B_{s_k}}\pi^*\alpha|_{\Sigma\cap\partial B_{s_k}}\right|.
\end{align*}
Since $j\in(I\cup J)^c$, by \eqref{good choice of s_k} we get that
\begin{align*}
    \left|\int_{\Sigma\cap\partial B_{s_k}}\pi^*\alpha|_{\Sigma\cap\partial B_{s_k}}\right|&=\left|\left<\pi_{\ast}\partial([\Sigma]\res B_{s_k}),\alpha\right>\right|\le ||\alpha||_{\ast}\m\big(\pi_{\ast}\partial([\Sigma]\res B_{s_k})\big)\\
    &\le\kappa\m\big(\pi_{\ast}([\Sigma]\res (B_{r}\smallsetminus B_{r/2}))\big)\\
    &\le\frac{3\kappa}{2}\int_{\Sigma\cap(B_r\smallsetminus B_{r/2})}|\opwedge_2d\pi(\vec\Sigma_0)|\, d\H^2. 
\end{align*}
Thus, we have obtained 
\begin{align}
\label{estimate morrey decrease 3}
    \nonumber
    \bigg|\int_{\Sigma\cap (B_{\rho_j}\smallsetminus B_{s_k})}\left<\pi^*d\alpha,\vec\Sigma\right>\, d\H^2\bigg|&\le\left|\int_{\Sigma\cap\partial B_{\rho_j}}\pi^*\alpha|_{\Sigma\cap\partial B_{\rho_j}}\right|\\
    &\quad+\frac{3\kappa}{2}\int_{\Sigma\cap(B_r\smallsetminus B_{r/2})}|\opwedge_2d\pi(\vec\Sigma_0)|\, d\H^2
\end{align}
and we just need to bound
\begin{align*}
    \left|\int_{\Sigma\cap\partial B_{\rho_j}}\pi^*\alpha|_{\Sigma\cap\partial B_{\rho_j}}\right|.
\end{align*}
To do this, we write the $1$-rectifiable closed curve $\Sigma\cap\partial B_{\rho_j}$ as
\begin{align*}
\Sigma\cap\partial B_{\rho_j}=\bigcup_{i=0}^{\infty}\Gamma_i,
\end{align*}
where $\Gamma_i$ is a Lipschitz connected closed curve in $B$. We let $\gamma_i:[0,\H^1(\Gamma_i)]\rightarrow B$ be the parametrization of $\Gamma_i$ through its arc-length, so that $|\gamma_i'|\equiv 1$ a.e. on $[0,\H^1(\Gamma_i)]$. First, fix $i\in\n$ and notice that for every smooth function $f:B\smallsetminus\{0\}\rightarrow\r$ such that $\bar f^i=0$, where
\begin{align*}
    \bar f^i:=\int_{\Gamma_i} f\, d\H^1,
\end{align*}
the following Poincaré type inequality holds for $\Gamma_i$:
\allowdisplaybreaks
\begin{align}
\label{poincare inequality}
    \nonumber
    \bigg(\int_{\Gamma_i} |f|^2\, d\H^1\bigg)^{1/2}&=\bigg(\int_0^{\H^1(\Gamma_i)} |f\circ\gamma_i|^2|\gamma_i'|\, d\L^1\bigg)^{1/2}\\
    \nonumber
    &=\bigg(\int_0^{\H^1(\Gamma_i)} |f\circ\gamma_i|^2\, d\L^1\bigg)^{1/2}\\
    \nonumber
    &\le\H^1(\Gamma_i)\bigg(\int_0^{\H^1(\Gamma_i)} |(f\circ\gamma_i)'|^2\, d\L^1\bigg)^{1/2}\\
    \nonumber
    &=\H^1(\Gamma_i)\bigg(\int_0^{\H^1(\Gamma_i)} |df_{\gamma_i}(\gamma_i')|^2\, d\L^1\bigg)^{1/2}\\
    \nonumber
    &=\H^1(\Gamma_i)\bigg(\int_0^{\H^1(\Gamma_i)} |df_{\gamma_i}(\vec\Sigma_{\rho_j})|^2|\gamma_i'|\, d\L^1\bigg)^{1/2}\\
    &=\H^1(\Gamma_i)\bigg(\int_{\Gamma_i}|df|_{\Sigma_{\rho_j}}|^2\, d\H^1\bigg)^{1/2}.
\end{align}
Secondly, since $\spt(\alpha)\subset\cp^{m-1}\smallsetminus H_{\varepsilon/2}$ and $\cp^{m-1}\smallsetminus H_{\varepsilon/2}$ is diffeomorphic to $\r^{2m-2}$ we can write $\alpha$ in coordinates $\{y_1,...y_{2m-2}\}$ on $\cp^{m-1}$ as 
\begin{align*}
    \alpha=\sum_{a=1}^{2m-2}a\alpha_{a}dy_{a}
\end{align*}
in order to get the expansion
\begin{align*}
    \pi^*\alpha|_{\Gamma_i}=\left<\pi^*\alpha,\gamma_i'\right>=\sum_{a=1}^{2m-2}(\alpha_a\circ\pi)\left<dy_a,d\pi|_{\Gamma_i}\right>.
\end{align*}
Moreover, we notice that 
\begin{align*}
    \big|d(\alpha_a\circ\pi)|_{\Sigma_{\rho_j}}\big|&\le\big|d\alpha_a\circ\pi\big|^2\big|d\pi|_{\Sigma_{\rho_j}}\big|^2\le\big|d\alpha\circ\pi\big|^2\big|d\pi|_{\Sigma_{\rho_j}}\big|^2\\
    &\le\max\bigg\{||\omega_{\cp^{m-1}}||_{\infty},\frac{\kappa}{\varepsilon^2}\bigg\}\big|d\pi|_{\Sigma_{\rho_j}}\big|^2\\
    &\le\max\bigg\{||\omega_{\cp^{m-1}}||_{\infty},\frac{\kappa}{\varepsilon^2}\bigg\}\big|d\pi|_{\Sigma_{\rho_j}}\big|^2\\
    &=M_m\big|d\pi|_{\Sigma_{\rho_j}}\big|^2,
\end{align*}
where
\begin{align*}
    M_m:=\max\bigg\{||\omega_{\cp^{m-1}}||_{\infty},\frac{\kappa}{\varepsilon^2}\bigg\}
\end{align*}
depends only on $m$.
Then, by \eqref{poincare inequality}, Hölder's inequality and point (3) in Lemma \ref{lemma sigma zero}, we estimate
\begin{align*}
    &\left|\int_{\Gamma_i}\pi^*\alpha|_{\Gamma_i}\right|\\
    &=\left|\sum_{a=1}^{2m-2}\int_{\Gamma_i}(\alpha_a\circ\pi)\left<dy_a,d\pi(\gamma_i')\right>\right|\\
    &=\left|\sum_{a=1}^{2m-2}\int_{\Gamma_i}(\alpha_a\circ\pi-\bar\alpha_k^i\circ\pi)\left<dy_a,d\pi(\gamma_i')\right>\right|\\
    &\le\sum_{a=1}^{2m-2}\int_{\Gamma_i}\big|\alpha_a\circ\pi-\bar\alpha_k^i\circ\pi\big|\big|d\pi|_{\Gamma_i}\big|^2\\
    &\le\sum_{a=1}^{2m-2}\bigg(\int_{\Gamma_i}\big|\alpha_a\circ\pi-\bar\alpha_k^i\circ\pi\big|^2\, d\H^1\bigg)^{1/2}\bigg(\int_{\Gamma_i}\big|d\pi|_{\Sigma_{\rho_j}}\big|\, d\H^1\bigg)^{1/2}\\
    &\le\H^1(\Gamma_i)\sum_{a=1}^{2m-2}\bigg(\int_{\Gamma_i}\big|d(\alpha_a\circ\pi)|_{\Sigma_{\rho_j}}\big|^2\, d\H^1\bigg)^{1/2}\bigg(\int_{\Gamma_i}\big|d\pi|_{\Sigma_{\rho_j}}\big|^2\, d\H^1\bigg)^{1/2}\\
    &\le\tilde M_m\H^1(\Gamma_i)\int_{\Gamma_i}\big|d\pi|_{\Sigma_{\rho_j}}\big|^2\, d\H^1\\
    &\le\tilde M_m\H^1(\Sigma\cap\partial B_{\rho_j})\int_{\Gamma_i}|\opwedge_2d\pi(\vec\Sigma_0)|\, d\H^1,
\end{align*}
where $\tilde M_m:=(2m-2)M_m$ and the last inequality follows by our choice of $\vec\Sigma_0$ (see point (4) in Lemma \ref{lemma sigma zero}).
Summing up over $i\in\n$ in the previous inequality, using the properties of good slices established in Lemma \ref{lemma good slices} and since $\Sigma\in\F$, we eventually get 
\begin{align}
\label{estimate morrey decrease 4}
    \nonumber
    &\left|\int_{\Sigma\cap\partial B_{\rho_j}}\pi^*\alpha|_{\Sigma\cap\partial B_{\rho_j}}\right|\\
    \nonumber
    &\le\tilde M_m\H^1(\Sigma\cap\partial B_{\rho_j})\int_{\Sigma\cap\partial B_{\rho_j}}|\opwedge_2d\pi(\vec\Sigma_0)|\, d\H^1\\
    \nonumber
    &\le\tilde M_m\Theta^2\H^2(\Sigma\cap B_{2^{-j_0}})\int_{\Sigma\cap(B_r\smallsetminus B_{r/2})}|\opwedge_2d\pi(\vec\Sigma_0)|\, d\H^2\\
    &\le\tilde M_m\Theta^2\big(e^{A+j_0}(1+A)\big)^{-1}\int_{\Sigma\cap(B_r\smallsetminus B_{r/2})}|\opwedge_2d\pi(\vec\Sigma_0)|\, d\H^2
\end{align}
Plugging \eqref{estimate morrey decrease 4} in \eqref{estimate morrey decrease 3} we get
\begin{align}
\label{estimate morrey decrease 5}
    \nonumber
    &\bigg|\int_{\Sigma\cap (B_{\rho_j}\smallsetminus B_{s_k})}\left<\pi^*d\alpha,\vec\Sigma\right>\, d\H^2\bigg|\\
    &\le\bigg(\tilde M_m\Theta^2\big(e^{A+2\ell}(1+A)\big)^{-1}+\frac{3\kappa}{2}\bigg)\int_{\Sigma\cap(B_r\smallsetminus B_{r/2})}|\opwedge_2d\pi(\vec\Sigma_0)|\, d\H^2.
\end{align}
Combining \eqref{estimate morrey decrease 5}, \eqref{estimate morrey decrease 2} and \eqref{estimate morrey decrease 1} and setting 
\begin{align*}
    \hat C:=\max\bigg\{\frac{2\tilde M_m\Theta^2}{e^{A+2\ell}(1+A)}+3\kappa,1,2\tilde\Xi\bigg\},
\end{align*}
we obtain
\begin{align*}
    \int_{\Sigma\cap (B_{r/2}\smallsetminus B_{s_k})}|\opwedge_2d\pi(\vec\Sigma_0)|\, d\H^2&\le\hat C\bigg(\int_{\Sigma\cap(B_r\smallsetminus B_{r/2})}|\opwedge_2d\pi(\vec\Sigma_0)|\, d\H^2+r^{1/4}\bigg)\\
    &\quad+\frac{1}{2}\int_{\Sigma\cap B_r}|\opwedge_2d\pi(\vec\Sigma_0)|\, d\H^2.
\end{align*}
By letting $k\rightarrow+\infty$ in the previous inequality, we obtain
\begin{align*}
    \int_{\Sigma\cap B_{r/2}}|\opwedge_2d\pi(\vec\Sigma_0)|\, d\H^2&\le\hat C\bigg(\int_{\Sigma\cap(B_r\smallsetminus B_{r/2})}|\opwedge_2d\pi(\vec\Sigma_0)|\, d\H^2+r^{1/4}\bigg)\\
    &\quad+\frac{1}{2}\int_{\Sigma\cap B_r}|\opwedge_2d\pi(\vec\Sigma_0)|\, d\H^2.
\end{align*}
and by subtracting from both sides the quantity 
\begin{align*}
    \frac{1}{2}\int_{\Sigma\cap B_{r/2}}|\opwedge_2d\pi(\vec\Sigma_0)|\, d\H^2,
\end{align*}
we get
\begin{align*}
    \int_{\Sigma\cap B_{r/2}}|\opwedge_2d\pi(\vec\Sigma_0)|\, d\H^2\le \bar C\bigg(\int_{\Sigma\cap(B_r\smallsetminus B_{r/2})}|\opwedge_2d\pi(\vec\Sigma_0)|\, d\H^2+r^{1/4}\bigg),
\end{align*}
where $\bar C>0$ is chosen big enough so that $\bar C\ge 2\hat C$ and
\begin{align*}
    \frac{\bar C}{\bar C+1}\ge 2^{-1/4}.
\end{align*}
By the hole filling technique and recalling that $r=2^{-j}$, we obtain
\begin{align*}
    \int_{\Sigma\cap B_{2^{-(j+1)}}}|\opwedge_2d\pi(\vec\Sigma_0)|\, d\H^2\le\theta\int_{\Sigma\cap B_{2^{-j}}}|\opwedge_2d\pi(\vec\Sigma_0)|\, d\H^2+\theta^{j+1},
\end{align*}
with $\theta\in(0,1)$ given by $\theta:=\bar C/(\bar C+1)\ge2^{-1/4}$ and our claim \eqref{step theta} follows. 

By \eqref{step theta} we obtain that
\begin{align*}
    E\big(2^{-j}\big)&\le\theta E\big(2^{-(j-1)}\big)+\theta^{j}\le\theta^2 E\big(2^{-(j-2)}\big)+2\theta^j\\
    &\le...\le\theta^{j_0}E\big(2^{-j_0}\big)\theta^{-j}+\big((j-j_0)\theta^{-j/2}\big)\theta^{-j/2}\\
    &\le\hat\Xi\tilde\theta^{-j}
\end{align*}
with $\tilde\theta:=\theta^{1/2}$ and 
\begin{align*}
    \hat\Xi:=\frac{\delta'\theta^{j_0}}{2}+\sup_{j\ge j_0}\big\{(j-j_0)\tilde\theta^{-j}\big\}<+\infty.
\end{align*}

In order to get \eqref{morrey decrease very true}, we notice that if $j\in I\cup J$ then either $j\in I$ or $j\in J\smallsetminus I$. In the first case, we have
\begin{align*}
    \int_{\Sigma\cap B_{2^{-(j+1)}}}|\opwedge_2d\pi(\vec\Sigma_0)|\, d\H^2&\le\frac{1}{2}\int_{\Sigma\cap B_{2^{-j}}}|\opwedge_2d\pi(\vec\Sigma_0)|\, d\H^2\\
    &\le\theta\int_{\Sigma\cap B_{2^{-j}}}|\opwedge_2d\pi(\vec\Sigma_0)|\, d\H^2.
\end{align*}
In the second case, by definition of $J$, it holds that
\begin{align*}
    \int_{\Sigma\cap B_{2^{-(j+1)}}}|\opwedge_2d\pi(\vec\Sigma_0)|\, d\H^2\le 5\Xi\,2^{-(j+1)/4}.
\end{align*}
By setting $\tilde\alpha=\min\{-\log_2\tilde\theta,1/4\}\in (0,1)$, we get 
\begin{align*}
    \int_{\Sigma\cap B_{2^{-j}}}|\opwedge_2d\pi(\vec\Sigma_0)|\, d\H^2\le 5\Xi\bigg(\frac{1}{2^j}\bigg)^{\alpha},
\end{align*}
which leads to \eqref{morrey decrease very true} with $\tilde C:=\max\{5\theta^{1/4}\Xi,\hat\Xi\}$.


\section{Almost pseudo-holomorphic foliations}

\begin{lem}
\label{level sets of locally approximable maps are closed}
Let $m\ge 2$ and let $(X,J_X,\omega_X)$ be a closed, almost Kähler, smooth $(2m-2)$-dimensional manifold. Assume that $v\in W^{1,2}(B,X)$ satisfies $v^*\vol_X\in L^1(B)$ and $d(v^*\vol_X)=0$ in $\D'(B)$. Then, there exists a representative of $v$ such that the co-area formula holds. Moreover, given such a representative, for $\vol_X$-a.e. $z\in X$ the following facts hold:
\begin{enumerate}
    \item $v^{-1}(z)$ is a countably $\H^{2}$-rectifiable subset of $B$;
    \item $(v^*\vol_X)_x\neq 0$, for $\H^2$-a.e. $x\in v^{-1}(z)$;
    \item $[v^{-1}(z)]$ is a cycle of finite mass.
\end{enumerate}
\begin{proof}
By \cite[Theorem 11, Theorem 12]{coareaformula} there exists a representative of $v$ such that both (1) and the co-area formula hold. Moreover, if we denote by $E\subset B$ the set of all the $x\in B$ such that $(v^*\vol_X)_x=0$, by the coarea formula we get
\begin{align*}
    0=\int_E|v^*\vol_X|_g\, d\vol_g=\int_X\H^2\big(v^{-1}(z)\cap E\big)\, d\vol_X(z),
\end{align*}
which implies that for $\vol_X$-a.e. $z\in X$ the set $v^{-1}(z)\cap E$ has vanishing $\H^2$-measure. Thus, (2) immediately follows. 
    
We are just left to prove (3). By the coarea formula, it follows that 
\begin{align*}
    \int_X\H^{2}\big(v^{-1}(z)\big)\, d\vol_X(z)=\int_B|v^*\vol_X|\, d\L^{2m} <+\infty.
\end{align*}
Hence, the function $X\ni z\xmapsto{f}\H^{2}\big(v^{-1}(z)\big)$ belongs to $L^1(X)$ and we know that a.e. $z\in X$ is a Lebesgue point for $f$ such that $f(z)<+\infty$. Fix any such point $z\in X$. By our choice of $z$, it holds that $\m([v^{-1}(z)])=\H^2\big(v^{-1}(z)\big)=f(z)<+\infty$. Hence, just need to show that $[v^{-1}(z)]$ is a cycle. Let $\exp_z:\r^{2m-2}\rightarrow X$ be the exponential map of $X$ at the point $z$. Denote by $\rho_0\in(0,+\infty)$ the injectivity radius of $X$ at $z$ and we define
\begin{align*}
    B_{\varepsilon}(z):=\exp_z\big(B_{\varepsilon}(0)\big), \qquad \mbox{ for every } \varepsilon\in(0,\rho_0).
\end{align*}
For every $\varepsilon\in(0,\rho_0)$, we let $\{\varphi_{\varepsilon,k}\}_{k\in\n}\subset C^{\infty}(X)$ be a sequence of smooth functions on $X$ such that:
\begin{enumerate}
    \item $\varphi_{\varepsilon,k}\equiv0$ on $X\smallsetminus B_{\varepsilon}(z)$;
    \item $0<\varphi_{\varepsilon,k}\le\big(\hspace{-0.5mm}\vol_X(B_{\varepsilon}(z))\big)^{-1}$ on $B_{\varepsilon}(z)$; 
    \item it holds that 
        \begin{align*}
            \varphi_{\varepsilon,k}\xrightarrow{k\rightarrow\infty}\frac{1}{\vol_{X}\big(B_{\varepsilon}(z)\big)}\rchi_{B_{\varepsilon}(z)}, \qquad \vol_{X}\mbox{-a.e. on } X.
        \end{align*}
\end{enumerate}  
Fix any $\alpha\in\D^{1}(B)$. By the coarea formula, it follows that
\begin{align*}
    \int_{B}d\alpha\wedge v^*(\varphi_{\varepsilon,k}\vol_{X})&=\int_{X}\varphi_{\varepsilon,k}(z)\bigg(\int_{v^{-1}(z)}d\alpha|_{v^{-1}(z)}\bigg)d\vol_{X}(z).
\end{align*}
Hence, by dominated convergence, we get
\begin{align*}
    \lim_{k\rightarrow+\infty}\int_{B}d\alpha\wedge v^*(\varphi_{\varepsilon,k}\vol_{X})&=\fint_{B_{\varepsilon}(z)}\bigg(\int_{v^{-1}(z)}d\alpha|_{v^{-1}(z)}\bigg)d\vol_{X}(z).
\end{align*}
Since $z$ is a Lebesgue point for $f$, we obtain
\begin{align}
\label{boundaryeq1}
    \lim_{\varepsilon\rightarrow 0^+}\lim_{k\rightarrow+\infty}\int_{B}d\alpha\wedge v^*(\varphi_{\varepsilon,k}\vol_{X})=\int_{v^{-1}(z)}d\alpha|_{v^{-1}(z)}=:\langle[v^{-1}(z)],d\alpha\rangle
\end{align}
Moreover, since $v$ is such that $d\big(v^*\vol_X)=0$ distributionally on $B$ and by the upper bound on $\varphi_{\varepsilon,k}$, it holds that
\begin{align}
\label{boundaryeq2}
\nonumber
    \left|\int_{B}d\alpha\wedge v^*(\varphi_{\varepsilon,k}\vol_{X})\right|&=\left|\int_{B}v^*\varphi_{\varepsilon,k}\big(d\alpha\wedge v^*\vol_{X}\big)\right|\\
    &\le\frac{1}{\vol_X(B_{\varepsilon}(z))}\left|\int_{B}d\alpha\wedge v^*\vol_{X}\right|=0,
\end{align}
for every $\varepsilon\in(0,\rho_0)$ and $k\in\n$. 

By \eqref{boundaryeq1} and \eqref{boundaryeq2} we get that $\left<[v^{-1}(x)],d\alpha\right>=0$ and, by arbitrariness of $\alpha\in\D^{1}(B)$, it follows that $\partial[v^{-1}(x)]=0$ in the sense of currents. The statement follows.
\end{proof}
\end{lem}

\begin{lem}
\label{J-holomorphic foliations}
Let $v\in W^{1,2}(B,X)$ be a weakly $(J,J_X)$-holomorphic map such that $v^*\vol_X\in L^1(B)$ and $d(v^*\vol_X)=0$ in $\D'(B)$. Then, there exist a representative of $u$ and a full measure set $\reg(v)\subset X$ such that:
\begin{enumerate}
    \item the co-area formula holds for $v$;
    \item for every $z\in\reg(u)$, the level set $v^{-1}(z)$ is a closed $J$-holomorphic curve in $B$.
\end{enumerate}
\begin{proof}
By Lemma \ref{level sets of locally approximable maps are closed}, it follows immediately that there exists a representative of $v$ such that the co-area formula holds and, for such a representative, $v^{-1}(z)$ is an $\H^2$-rectifiable subset of $B$ with $\partial[v^{-1}(z)]=0$, for $\vol_{X}$-a.e. $z\in X$. Thus, we are just left to show that $v^{-1}(z)$ is $J$-holomorphic, for a.e. $z\in X$. By the co-area formula and since $v$ is weakly $(J,J_X)$-holomorphic, for $\vol_X$-a.e. $z\in X$ the form $v^*\vol_X$ is non-vanishing on $v^{-1}(z)$ and $dv(Jw)=J_Xdv(w)$ for every $w\in\r^{2m}$, up to some $\H^2$-negligible set. For such $z\in X$, the orienting vector field to $v^{-1}(z)$ is given by
\begin{align*}
    \vec\Sigma:=\frac{\ast v^*\vol_X}{|v^*\vol_X|_g}.
\end{align*}
We claim that $\vec\Sigma$ is $J$-invariant for $\H^2$-a.e. $x\in v^{-1}(z)$. Indeed, given any $x\in v^{-1}(z)$ such that $(v^*\vol_X)_x\neq 0$, we pick an orthonormal basis of $T_{v(x)}^*X$ of the form $\{\xi_1,J_X\xi_1,...,\xi_{m-1}, J_X\xi_{m-1}\}$ and we notice that
\begin{align*}
    (v^*\vol_X)_x&=\frac{1}{(m-1)!}v^*(\xi_1\wedge J_X\xi_1\wedge...\wedge\xi_{m-1}\wedge J_X\xi_{m-1})\\
    &=\frac{1}{(m-1)!}v^*\xi_1\wedge v^*J_X\xi_1\wedge...\wedge v^*\xi_{m-1}\wedge v^*J_X\xi_{m-1}\\
    &=\frac{1}{(m-1)!}v^*\xi_1\wedge J(v^*\xi_1)\wedge...\wedge v^*\xi_{m-1}\wedge J(v^*\xi_{m-1}).
\end{align*}
This clearly implies that $\vec\Sigma$ is $J$-invariant and the statement follows.
\end{proof}
\end{lem}

In the following lemma, which generalises the model situation presented in Lemma \ref{J-holomorphic foliations}, we will adopt the notation developed in Appendix A. Moreover, we will denote by $X$ the product space $X:=\cp^1\times\cp^{m-2}$ and by $p_1:X\rightarrow\cp^1$ and $p_2:X\rightarrow\cp^{m-2}$ the canonical projections on the first and on the second factor respectively. We will endow $X$ with the complex structure $J_X:=p_1^*j_1+p_2^*j_{m-2}$ and with symplectic form $\omega_X:=p_1^*\omega_{\cp^1}+p_2^*\omega_{\cp^{m-2}}$ in order to obtain the Kähler manifold $(X,J_X,\omega_X)$.

\begin{lem}
\label{almost J-holomorphic foliations}
Let $m,n\in\n_0$ be such that $m\ge 3$. Let $u\in W^{1,2}(B,\cp^n)$ be weakly $(J,j_n)$-holomorphic and locally approximable. If $n\ge2$, then for a.e. $(q_1,...,q_{n-1},p)\in \cp^{n}\times\cp^{n-1}\times...\times\cp^{2}\times\cp^{m-1}$ the map $v_{q_1,...,q_{n-1},p}:=\big(F_{q_{n-1}}\circ...\circ F_{q_1}\circ u, F_p\circ\pi\big):B\rightarrow X$ has the following properties:
\begin{enumerate}
    \item $v_{q_1,...,q_{n-1},p}\in W^{1,2}(B,X)$;
    \item $v_{q_1,...,q_{n-1},p}^*\vol_X\in L^1(B)$;
    \item there exists a set $\reg(v_{q_1,...,q_{n-1},p})\subset X$ such that
    \begin{align*}
        \vol_X\big((X\smallsetminus\reg(v_{q_1,...,q_{n-1},p})\big)=0
    \end{align*}
    and for every $(y,z)\in\reg(v_{q_1,...,q_{n-1},p})$ the $\H^2$-rectifiable set $v_{q_1,...,q_{n-1},p}^{-1}(y,z)$ is a closed almost $J$-holomorphic curve in $B$, in the sense of Definition \ref{definition of almost pseudoholomorphic curve}. Moreover, the constants $\ell>0$ and $\gamma\in(0,1]$ can be chosen as $\ell=2\sqrt{2\Lip(\Omega)}$ and $\gamma=1/2$.
\end{enumerate}
If $n=1$, analogous properties hold for the map $v_p:=(u,F_p\circ\pi):B\rightarrow X$ and for a.e. $p\in\cp^{m-1}$.
\begin{proof}
Since the techniques are identical both in the case $n=1$ and $n\ge2$, we just focus on the second one. 

Let $Y:=\cp^{n}\times...\times\cp^{2}$. First, we want to prove (1). By Lemma \ref{properties of singular projections}, we know that  $p_2\circ v_{q_1,...,q_{n-1},p}$ belongs to $W^{1,2}(B,\cp^{m-2})$, for every $p\in\cp^{m-1}$. We claim that $p_1\circ v_{1_1,...,q_{n-1},p}=F_{q_{n-1}}\circ...\circ F_{q_1}\circ u$ belongs to $W^{1,2}(B,\cp^1)$ for a.e. $(q_1,...,q_{n-1})\in Y$. Indeed, notice that the map $F_{q_{n-1}}\circ...\circ F_{q_1}\circ u$ is weakly $(J,j_1)$-holomorphic, for every $(q_1,...,q_{n-1})\in Y$. Thus, by Corollary \ref{weakholcor} we have that
\begin{align*}
    &\int_B|d(F_{q_{n-1}}\circ...\circ F_{q_1}\circ u)|_g^2\, d\vol_g\\
    &=2\int_B(F_{q_{n-1}}\circ...\circ F_{q_1}\circ u)^*\omega_{\cp^1}\wedge\frac{\Omega^{m-1}}{(m-1)!}.
\end{align*}
Hence, by Lemma \ref{lemma averaging property} we obtain
\begin{align*}
    &\int_{B}\varphi\bigg(\int_Y|d(F_{q_{n-1}}\circ...\circ F_{q_1}\circ u)|_g^2\, d\vol_{Y}(q_1,...,q_{n-1})\bigg)\, d\vol_g\\
    &=2\int_{B}\varphi\bigg(\int_Y(F_{q_{n-1}}\circ...\circ F_{q_1}\circ u)^*\omega_{\cp^1}\, d\vol_{Y}(q_1,...,q_{n-1})\bigg)\wedge\frac{\Omega^{m-1}}{(m-1)!}\\
    &=2D\int_B\varphi\,u^*\omega_{\cp^n}\wedge\frac{\Omega^{m-1}}{(m-1)!}=D\int_B\varphi|du|_g^2\, d\vol_g<+\infty,
\end{align*}
where $D:=B_{n+1}\cdot...\cdot B_3$, for every $\varphi\in C_c^{\infty}(B)$. Thus, we get that 
\begin{align}
\label{estimate first component 1}
    \int_Y|d(F_{q_{n-1}}\circ...\circ F_{q_1}\circ u)|_g^2\, d\vol_{Y}(q_1,...,q_{n-1})=D|du|_g^2, 
\end{align}
for a.e. $x\in B$. By integrating both sides of \eqref{estimate first component 1} on $B$ and by Fubini's theorem, we get
\begin{align}
\label{estimate first component 2}
    \nonumber
    &\int_B\bigg(\int_Y|d(F_{q_{n-1}}\circ...\circ F_{q_1}\circ u)|_g^2\, d\vol_{Y}(q_1,...,q_{n-1})\bigg)\, d\L^{2m}\\
    \nonumber
    &=\int_Y\bigg(\int_B|d(F_{q_{n-1}}\circ...\circ F_{q_1}\circ u)|_g^2\, d\vol_g\bigg)\, d\vol_{Y}(q_1,...,q_{n-1})\\
    &=D\int_B|du|_g^2\, d\vol_g<+\infty.
\end{align}
Since \eqref{estimate first component 2} directly implies that 
\begin{align*}
    \int_B|d(F_{q_{n-1}}\circ...\circ F_{q_1}\circ u)|_g^2\, d\vol_g<+\infty
\end{align*}
for $\vol_Y$-a.e. $(q_1,...,q_{n-1})\in Y$, point (1) follows.

Next, we turn to show (2). By \eqref{estimate first component 1}, Lemma \ref{lemma averaging property}, \eqref{equivalence of norms}, \eqref{estimate second slicing} and by Fubini's theorem, we have
\begin{align}
\label{estimate 2m-4 number 1}
    \nonumber
    &\int_{Y}\bigg(\int_B\big|v_{q_1,...,q_{n-1},p}^*\vol_X\big|_g\, d\vol_g\bigg)\, d\vol_{Y\times\cp^{m-1}}(q_1,...,q_{n-1},p)\\
    \nonumber
    &\le\int_B\bigg(\int_Y\big|d(F_{q_{n-1}}\circ...\circ F_{q_1}\circ u)\big|_g^2\, d\vol_{Y\times\cp^{m-1}}(q_1,...,q_{n-1},p)\bigg)\\
    \nonumber
    &\quad\cdot|d(F_p\circ\pi)|_g^{2m-4}\, d\vol_g\\
    \nonumber
    &=DG^{m-2}\int_B|du|_g^2\cdot|d(F_p\circ\pi)|^{2m-4}\, d\vol_g\\
    &\le DG^{m-2}\int_B\frac{|du|_g^2}{\dist(\,\cdot\,,L_p)^{2m-4}}\, d\vol_g,
\end{align}
where $L_p$ is defined as in Appendix A. For any $\rho\in(0,1)$, define $L_p^{\rho}:=(L_p+B_{\rho})\cap B$. By the almost monotonicity formula \eqref{monotonicity formula below}, we get that
\begin{align*}
    \int_{L_p^{\rho}\smallsetminus L_p^{\rho/2}}\frac{|du|_g^2}{\dist(\,\cdot\, 
    ,L_p)^{2m-4}}\, d\vol_g&\le \frac{2^{2m-4}}{\rho^{2m-4}}\int_{B_{\rho}}|du|_g^2\, d\vol_g\\
    &\le 2^{2m-2}\frac{e^{A\rho}(1+A\rho)}{\rho^{2m-2}}\int_{B_{\rho}}|du|_g^2\, d\vol_g\frac{\rho^2}{4}\\
    &\le\bigg(2^{2m-2}e^{A}(1+A)\int_{B}|du|_g^2\, d\vol_g\bigg)\frac{\rho^2}{4},
\end{align*}
for every $\rho\in(0,1)$. By iteration (see also the proof of Lemma \ref{lemma projected sigma}) we get
\begin{align}
\label{estimate 2m-4 number 2}
    \int_{L_p^{\rho}}\frac{|du|_g^2}{\dist(\,\cdot\, 
    ,L_p)^{2m-4}}\, d\vol_g\le \frac{2^{2m}e^{A}(1+A)}{3}\bigg(\int_{B}|du|_g^2\, d\vol_g\bigg)\rho^2.
\end{align}
Combining \eqref{estimate 2m-4 number 1} and \eqref{estimate 2m-4 number 2}, we obtain
\begin{align}
\label{integrated estimate}
    \nonumber
    &\int_{Y\times\cp^{m-1}}\bigg(\int_B\big|v_{q_1,...,q_{n-1},p}^*\vol_X\big|_g\, d\vol_g\bigg)\, d\vol_{Y\times\cp^{m-1}}(q_1,...,q_{n-1},p)\\
    &\le C\int_B|du|_g^2\, d\vol_g<\infty,
\end{align}
where $C>0$ is a constant depending on $m$, $n$ and $\Lip(\Omega)$. Again, point (2) follows by Fubini's theorem. 

We are left to prove (3). First, we claim that for a.e. $(q_1,...,q_{n-1},p)\in Y\times\cp^{m-1}$ it holds that 
$d\big(v_{q_1,...,q_{n-1},p}^*\vol_X)=0$ in the sense of distributions. Indeed, we already know that for $\vol_Y$-a.e. $(q_1,...,q_{n-1})\in Y$ the estimate \eqref{estimate first component 2} holds. Fix any $\alpha\in\D^1(B)$. Notice that by estimates \eqref{estimate first component 2} and \eqref{estimate 2m-4 number 2} we obtain
\begin{align*}
    \bigg|\int_Bv_{q_1,...,q_{n-1},p}^*\vol_X\wedge d\alpha\bigg|&\le C|d\alpha|_{\ast}\int_{L_p^{\rho}}\frac{|du|_g^2}{\dist(\,\cdot\,L_p)^{2m-4}}\, d\vol_g\\
    &\le C|d\alpha|_{\ast}\bigg(\int_B|du|_g^2\, d\vol_g\bigg)\rho^2, \qquad\forall\, \rho\in(0,1),
\end{align*}
where $C>0$ is a constant depending only on $m$, $n$ and $\Lip(\Omega)$. By letting $\rho\rightarrow 0^+$, we get
\begin{align*}
    \int_{Y\times\cp^{m-1}}&\bigg|\int_Bv_{q_1,...,q_{n-1},p}^*\vol_X\wedge d\alpha\bigg|\, d\vol_{Y\times\cp^{m-1}}(q_1,...,q_{n-1},p)=0.
\end{align*}
By arbitrariness of $\alpha\in D^1(B)$, our claim follows.

Let $E\subset Y\times\cp^{m-1}$ be the set of all the $n$-tuples $(q_1,...,q_{n-1},p)\in Y\times\cp^{m-1}$ such that
\begin{enumerate}
    \item (1) and (2) hold;
    \item $d\big(v_{q_1,...,q_{n-1},p}^*\vol_X)=0$ in the sense of distributions.
\end{enumerate} 
By what we have shown so far, we have $\,\vol_{Y\times\cp^{m-1}}(E^c)=0$. Fix any $(q_1,...,q_{n-1},p)\in E$. We fix the representative of the map $v_{q_1,...,q_{n-1},p}$ given by Lemma \ref{level sets of locally approximable maps are closed}. Thus, we know that the following facts hold for $\vol_X$-a.e. $(y,z)\in X$:
\begin{enumerate}
    \item the set $v_{q_1,...,q_{n-1},p}^{-1}(y,z)$ is $\H^2$-rectifiable;
    \item $(v_{q_1,...,q_{n-1},p}^*\vol_{X})_x\neq 0$, for $\H^2$-a.e. $x\in v_{q_1,...,q_{n-1},p}^{-1}(y,z)$ and the rectifiable set $v_{q_1,...,q_{n-1},p}^{-1}(y,z)$ is oriented by the $\H^2$-measurable and unitary field $2$-vectors given by:
        \begin{align*}
            \vec\Sigma:=\frac{\big(\hspace{-1mm}*(v_{q_1,...,q_{n-1},p}^*\vol_X)\big)^{\sharp}}{\big|v_{q_1,...,q_{n-1},p}^*\vol_X\big|_g}
        \end{align*}
    \item $\partial[v_{q_1,...,q_{n-1},p}^{-1}(y,z)]=0$.
\end{enumerate}

Hence, we just need to show that $v_{q_1,...,q_{n-1},p}^{-1}(y,z)$ is almost $J$-holomorphic according to Definition \ref{definition of almost pseudoholomorphic curve}, i.e. we claim that there exists some $J$-invariant and $\H^2$-measurable field of $g$-unitary $2$-vectors $\vec\Sigma_J:v_{q_1,...,q_{n-1},p}^{-1}(y,z)\rightarrow\bigwedge_2\r^{2m}$ such that
\begin{align}
\label{estimate sigma_J}
    \big|\vec\Sigma-\vec\Sigma_J\big|\le\ell|\cdot|^{\gamma},
\end{align}
for some $\ell>0$ and $\gamma\in(0,1]$. In order to prove our claim, consider the following $\H^2$-measurable and $g$-unitary fields respectively of $2$-vectors and $4$-vectors on $v_{q_1,...,q_{n-1},p}^{-1}(y,z)$:
\begin{align*}
    \vec\Sigma^1&:=\frac{(p_1\circ v_{q_1,...,q_{n-1},p})^*\omega_{\cp^1}\big)^{\sharp}}{\big|(p_1\circ v_{q_1,...,q_{n-1},p})^*\omega_{\cp^1}\big|_g}\\
    \vec\Sigma^2&:=\frac{\big(\hspace{-1mm}*(p_2\circ v_{q_1,...,q_{n-1},p})^*\vol_{\cp^{m-2}}\big)^{\sharp}}{\big|(p_2\circ v_{q_1,...,q_{n-1},p})^*\vol_{\cp^{m-2}}\big|_g}=\frac{\big(\hspace{-1mm}*(F_p\circ\pi)^*\vol_{\cp^{m-2}}\big)^{\sharp}}{\big|(F_p\circ\pi)^*\vol_{\cp^{m-2}}\big|_g}.
\end{align*}
Notice that they are both well defined $\H^2$-a.e. on $v_{q_1,...,q_{n-1},p}^{-1}(y,z)$, since $(v_{q_1,...,q_{n-1},p}^*\vol_X)_x\neq 0$ for $\H^2$-a.e. $x\in v_{q_1,...,q_{n-1},p}^{-1}(y,z)$. Fix any $x\in v_{q_1,...,q_{n-1},p}^{-1}(y,z)$ such that $(v_{q_1,...,q_{n-1},p}^*\vol_X)_x\neq 0$, so that the subspace $W_1:=\text{span}\{\vec\Sigma^1(x)\}$ is a $J$-holomorphic $2$-plane and $W_2:=\text{span}\{\vec\Sigma^2(x)\}$ is a $J_0$-holomorphic $4$-plane. Let $W:=\text{span}\{\vec\Sigma(x)\}$ and notice that $W=W_1\cap W_2$, $W_2=(W_1^{\perp}\cap W_2)\oplus W$ and $\dim(W)=2$. Moreover, we have that 
\begin{align*}
    4=\dim(W_2)=\dim(W_1^{\perp}\cap W_2)+\dim(W)=\dim(W_1^{\perp}\cap W_2)+2,
\end{align*}
which implies $\dim(W_1^{\perp}\cap W_2)=2$. We let $\{e_1,e_3\}$ be an $g$-orthonormal basis of $W$ and let $\{v,w\}$ be an $g$-orthonormal basis of $W_1^{\perp}\cap W_2$. By construction, $\{e_1,e_3,v,w\}$ is an $\Omega_0$-orthonormal basis of $W_2$ and we can write
\begin{align*}
    \vec\Sigma^2(x):=e_1\wedge e_3\wedge v\wedge w.
\end{align*}
If $\vec\Sigma^2(x)$ is $J$-invariant, we set $\vec\Sigma_J^2(x):=\vec\Sigma^2(x)$. 

If not, notice that $\{e_1,Je_1,e_3,Je_3,v-Je_1,w-Je_3\}$ is a linearly independent set. Let $e_2,e_4$ be the unique unitary vectors such that $\{e_1,Je_1,e_3,Je_3,e_2,e_4\}$ is an $g$-orthonormal set such that 
\begin{align*}
    \text{span}\{e_1,Je_1,e_3,Je_3,v-Je_1,w-Je_3\}=\text{span}\{e_1,Je_1,e_3,Je_3,e_2,e_4\}.
\end{align*}
Exactly as in the proof of Lemma \ref{lemma sigma zero}, it follows that there exist two angles $\phi_1,\phi_2\in[0,2\pi]$ such that
\begin{align*}
    \vec\Sigma^2(x):=e_1\wedge(\cos\phi_1 Je_1+\sin\phi_1 e_2)\wedge e_3\wedge (\cos\phi_2 Je_3+\sin\phi_2 e_4).
\end{align*}
The same computation as in Lemma \ref{lemma sigma zero} leads to 
\begin{align*}
    1-\Lip(\Omega)|x|\le\cos\phi_1\cos\phi_2\le 1+\Lip(\Omega)|x|.
\end{align*}
We set 
\begin{align*}
    \Sigma_J^2(x)=e_1\wedge\cos\phi_1Je_1\wedge e_3\wedge\cos\phi_2Je_3 
\end{align*}
and we compute
\begin{align*}
    \big|\vec\Sigma^2(x)-\vec\Sigma_J^2(x)\big|_g^2&=(1-\cos\phi_1\cos\phi_2)^2+\sin^2\phi_1\\
    &=1+\cos^2\phi_2-2\cos\phi_1\cos\phi_2\\
    &\le 2(1-\cos\phi_1\cos\phi_2)\le 2\Lip(\Omega)|x|,
\end{align*}
which leads to
\begin{align}
\label{estimate sigma_J^2}
    \big|\vec\Sigma^2(x)-\vec\Sigma_J^2(x)\big|_g\le\sqrt{2\Lip(\Omega)}|x|^{1/2}.
\end{align}
Eventually, we define 
\begin{align*}
    \vec\Sigma_J:=\ast\big(\vec\Sigma^1\wedge *\vec\Sigma_J^2\big).
\end{align*}
By construction, $\vec\Sigma_J$ is an $\H^2$-measurable and unitary field of $2$-vectors on $v_{q_1,...,q_{n-1},p}^{-1}(y,z)$. Moreover, by \eqref{estimate sigma_J^2}, we have  
\begin{align*}
    \big|\vec\Sigma^1\wedge\ast\vec\Sigma^2\big|_g&=\big|\vec\Sigma^1\wedge\ast(\vec\Sigma^2-\vec\Sigma_J^2)+\vec\Sigma^1\wedge\ast\vec\Sigma_J^2\big|_g\\
    &\ge\big|\vec\Sigma^1\wedge\ast\vec\Sigma_J^2|_g-\big|\vec\Sigma^1\wedge\ast(\vec\Sigma^2-\vec\Sigma_J^2)|_g\\
    &\ge 1-\big|\vec\Sigma^2-\vec\Sigma_J^2\big|_g\ge 1-\sqrt{2\Lip(\Omega)}|\cdot|^{1/2},
\end{align*}
which leads to 
\begin{align*}
    \big|\vec\Sigma-\vec\Sigma_J\big|&\le G\big|\vec\Sigma-\vec\Sigma_J\big|_g=G\bigg|\frac{\vec\Sigma^1\wedge\ast\vec\Sigma^2}{\big|\vec\Sigma^1\wedge\ast\vec\Sigma^2\big|}-\vec\Sigma^1\wedge\ast\vec\Sigma_J^2\bigg|_g\\
    &\le2\sqrt{2\Lip(\Omega)}G|\cdot|^{1/2}.
\end{align*}
Hence, \eqref{estimate sigma_J} holds with $\ell=2\sqrt{2\Lip(\Omega)}G$ and $\gamma=1/2$. The statement follows.
\end{proof}
\end{lem}


\section{Proof of the main theorem}

This section is entirely devoted to proof Theorems \ref{main theorem introduction} and \ref{main theorem}, whose local versions will be recalled now for the reader's convenience.

\begin{thm*}
Let $m,n\in\n_0$ be such that $m\ge 2$ and let $B$ denote the open unit ball in $\r^{2m}$. Assume that $u\in W^{1,2}(B,\cp^n)$ is weakly $(J,j_n)$-holomorphic and locally approximable. 

Then, $u$ has a unique tangent map at the origin.
\end{thm*}

\begin{thm*}
Let $m,n\in\n_0$ be such that $m\ge 2$ and let $B$ denote the open unit ball in $\r^{2m}$. Assume that $u\in W^{1,2}(B,\cp^n)$ is weakly $(J,j_n)$-holomorphic and locally approximable. 

Then, the $(2m-2)$-cycle $T_u\in\D_{2m-2}(B)$ has a unique tangent cone at the origin.
\end{thm*}

In the first two subsections, we treat the proof of Theorem \ref{main theorem}. We will first address the easy case $m=2$, $n=1$ in subsection 6.1, in order to clarify which will be the main ideas in order to proceed towards higher dimensions and codimensions. The general case will be discussed in subsection 6.2. Finally, in subsection 6.3 we will show how Theorem \ref{main theorem introduction} can be obtained as a consequence of Theorem \ref{main theorem}.

\subsection{A model problem}
Throughout all this subsection, $B\subset\r^4$ will denote the open unit ball in $\r^4$. 

Let $u\in W^{1,2}(B,\cp^1)$ be weakly $(J,j_1)$-holomorphic and locally approximable. As usual, $\pi:B\rightarrow\cp^1$ denotes the Hopf map.

If $\theta(0,u)<\varepsilon_0$, then $u$ is smooth in a neeighbourhood of $0$ by Theorem \ref{epsilonreg} and the statement follows. Assume then that $\theta(0,u)\ge\varepsilon_0$. 

We use the same notations and labeling for the constants as in sections 3 and 4. Since $u^*\omega_{\cp^1}\in L^1(B)$, by using Lemma \ref{J-holomorphic foliations} with $X=\cp^1$, we get that there exists a representative of $u$ and a full measure set $\reg(u)\subset\cp^1$ such that:
\begin{enumerate}
    \item the co-area formula holds for $u$;
    \item for every $y\in\reg(u)$, the level set $u^{-1}(y)$ is a closed $J$-holomorphic curve.
\end{enumerate}
Hence, all the estimates in section 3 will be used assuming $\vec\Sigma_J=\vec\Sigma$, $\ell=0$ and $\gamma=1$, as we stressed out in Remark \ref{pseudoholomorphic vs almost pseudoholomorphic case}.

For every $k\in\n_0$, we consider the set $E_k\subset\cp^1$ given by all the points $y$ in $\reg(u)$ such that 
\begin{enumerate}
    \item $\H^2(u^{-1}(y)\cap B_{2^{-k}})<\big((e^{A}(1+A)\big)^{-1}$;
    \item $\begin{aligned}\int_{u^{-1}(y)\cap B_{2^{-k}}}|\opwedge_2d\pi(\vec\Sigma_0^y)|\, d\H^2<\frac{\delta'}{2},\end{aligned}$
\end{enumerate}
where $\delta'>0$ is the constant introduced in section 3.2 and $\vec\Sigma_0^y$ is built as shown in Lemma \ref{lemma sigma zero} starting from the $J$-holomorphic field of $2$-vectors given by
\begin{align*}
    \vec\Sigma^y:=\frac{\ast(u^*\omega_{\cp^1})^{\sharp}}{|u^*\omega_{\cp^1}|_g},
\end{align*}
which orients the closed $J$-holomorphic curve $u^{-1}(y)$ for every $y\in\reg(u)$.
We notice that $E_{k-1}\subset E_k$ for every $k\in\n_0$. Moreover, since $\reg(u)\subset\cp^1$ has full measure in $\cp^1$, we get
\begin{align*}
    \vol_{\cp^1}\bigg(\cp^1\smallsetminus\bigcup_{k=1}^{+\infty} E_k\bigg)=0
\end{align*}
For every $k\in\n_0$, we define the localized current $T_k:=T_u\res u^{-1}(E_k)$, i.e. 
\begin{align*}
   \left<T_k,\alpha\right>:=\int_{u^{-1}(E_k)}u^*\omega_{\cp^1}\wedge \alpha \qquad \forall\, \alpha\in\D^2(B).
\end{align*}
\textbf{Claim}. We claim that every $T_k$ has a unique tangent cone at the origin. First, notice that $T_k$ is a normal $2$-cycle on $B$ semicalibrated by $\Omega$. By definition of $E_k$ and by Proposition \ref{monotonicity formula sigma}, for every $y\in E_k$ we get that 
\begin{align*}
e^{A\rho}(1+A\rho)\frac{\H^{2}\big(u^{-1}(y)\cap B_{\rho}\big)}{\rho^2}&-e^{A\sigma}(1+A\sigma)\frac{\H^{2}\big(u^{-1}(y)\cap B_{\sigma}\big)}{\sigma^2}\\
&\ge\int_{u^{-1}(y)\cap(B_{\rho}\smallsetminus B_{\sigma})}\frac{1}{|\cdot|^2}\big<\Omega_t,\vec\Sigma^y\big>\, d\H^2
\end{align*}
and
\begin{align*}
e^{-A\rho}(1-A\rho)\frac{\H^{2}\big(u^{-1}(y)\cap B_{\rho}\big)}{\rho^2}&-e^{-A\sigma}(1-A\sigma)\frac{\H^{2}\big(u^{-1}(y)\cap B_{\sigma}\big)}{\sigma^2}\\
&\le\int_{u^{-1}(y)\cap\left(B_{\rho}\smallsetminus B_{\sigma}\right)}\frac{1}{|\cdot|^2}\big<\Omega_t,\vec\Sigma_J^y\big>\, d\H^2,
\end{align*}
for every $0<\sigma<\rho<1$. Since a direct computation leads to
\begin{align*}
    \frac{\m(T_k\res B_{\rho})}{\rho^2}=\frac{1}{\rho^2}\int_{E_k}\H^2\big(u^{-1}(y)\cap B_{\rho}\big)\,d\vol_{\cp^1}(y),
\end{align*}
by integrating on $E_k$ the two previous inequalities we get the following almost monotonicity formulas for the current $T_k$:
\begin{align}
\label{monotonicity formula below T_k}
\nonumber
&e^{A\rho}(1+A\rho)\frac{\m(T_k\res B_{\rho})}{\rho^2}-e^{A\sigma}(1+A\sigma)\frac{\m(T_k\res B_{\sigma})}{\sigma^2}\\
&\ge\int_{E_k}\bigg(\int_{u^{-1}(y)\cap(B_{\rho}\smallsetminus B_{\sigma})}\frac{1}{|\cdot|^2}\big<\Omega_t,\vec\Sigma^y\big>\, d\H^2\bigg)\,d\vol_{\cp^1}(y),
\end{align}
\begin{align}
\label{monotonicity formula above T_k}
\nonumber
&e^{-A\rho}(1-A\rho)\frac{\m(T_k\res B_{\rho})}{\rho^2}-e^{-A\sigma}(1-A\sigma)\frac{\m(T_k\res B_{\sigma})}{\sigma^2}\\
&\le\int_{E_k}\bigg(\int_{u^{-1}(y)\cap(B_{\rho}\smallsetminus B_{\sigma})}\frac{1}{|\cdot|^2}\big<\Omega_t,\vec\Sigma^y\big>\, d\H^2\bigg)\,d\vol_{\cp^1}(y),
\end{align}
for every $0<\sigma<\rho<1$. Equation \eqref{monotonicity formula below T_k} immediately implies that function
\begin{align*}
    (0,1)\ni\rho\mapsto e^{A\rho}(1+A\rho)\frac{\m(T_k\res B_{\rho})}{\rho^2}
\end{align*}
is monotonically non-decreasing. Thus, the density of the current $T_k$ at zero, which is given by 
\begin{align*}
    \theta(T_k,0):=\lim_{\rho\rightarrow 0^+}\frac{\m(T_k\res B_{\rho})}{\rho^2}=\lim_{\rho\rightarrow 0^+}e^{A\rho}(1+A\rho)\frac{\m(T_k\res B_{\rho})}{\rho^2}
\end{align*}
exists and is finite. Moreover, by \eqref{monotonicity formula above T_k}, the coarea formula, \eqref{morrey decrease} and the estimate \eqref{useful estimate projected sigma 2}, it follows that
\begin{align}
\label{morrey decrease T_k}
    \nonumber
    &\bigg|\frac{\m(T_k\res B_{\rho})}{\rho^2}-\theta(0,T_k)\bigg|\\
    \nonumber
    &\le C\bigg|\int_{u^{-1}(E_k)\cap B_{\rho}}u^*\omega_{\cp^1}\wedge\frac{\Omega_t}{|\cdot|^2}\bigg|\\
    \nonumber
    &\le C\bigg|\int_{u^{-1}(E_k)\cap B_{\rho}}u^*\omega_{\cp^1}\wedge\pi^*\omega_{\cp^1}\bigg|+C\bigg|\int_{u^{-1}(E_j)\cap B_{\rho}}u^*\omega_{\cp^1}\wedge\frac{(\Omega-\Omega_0)_t}{|\cdot|^2}\bigg|\\
    \nonumber
    &\le C\int_{E_k}\bigg|\int_{u^{-1}(y)\cap B_{\rho}}\pi^*\omega_{\cp^1}\bigg|\, d\vol_{\cp^1}(y)\\
    \nonumber
    &\quad+C\int_{E_k}\int_{u^{-1}(y)\cap B_{\rho}}\frac{|\Omega-\Omega_0|}{|\cdot|^2}\, d\H^2\, d\vol_{\cp^1}(y)\\
    \nonumber
    &\le C\vol_{\cp^1}(\cp^1)\rho^{\alpha}+C\int_{E_k}\int_{u^{-1}(y)\cap B_{\rho}}\frac{1}{|\cdot|}\, d\H^2\, d\vol_{\cp^1}(y)\\
    &\le C\vol_{\cp^1}(\cp^1)\rho^{\alpha}+C\vol_{\cp^1}(\cp^1)\rho\le C\vol_{\cp^1}(\cp^1)\rho^{\alpha},
\end{align}
for every $\rho\in(0,\tilde r)$, where the constant $C>0$ and $\alpha,\tilde r\in (0,1)$ all depend just on $k$ and on $\Lip(\Omega)$. From the Morrey decay \eqref{morrey decrease T_k}, uniqueness of tangent cone for $T_k$ follows by standard arguments. 

\textbf{Conclusion}. For every $j\in\n_0$, consider the residual current $R_j:=T_u-T_j$. By the same arguments that we have used in the proof of the previous claim, we conclude that $R_j$ is a normal $2$-cycle in $B$ which is semicalibrated by $\Omega$. In particular, the quantity
\begin{align}
\label{residual mass}
    e^{A\rho}(1+A\rho)\frac{\m(R_j\res B_{\rho})}{\rho^2}
\end{align}
is non-decreasing in $\rho\in (0,1)$. Therefore, the limit as $\rho\rightarrow 0^+$ of the quantity \eqref{residual mass} exists and it is finite. Then, since the quantity \eqref{residual mass} is also non-increasing in $j\in\n$ and going to $0$ as $j\rightarrow +\infty$, we are allowed to exchange the limits in the following chain of equalities and we get
\begin{align}
\label{no neck property}
    \nonumber
    \lim_{j\rightarrow+\infty}\lim_{\rho\rightarrow 0^+}\frac{\m(R_j\res B_{\rho})}{\rho^2}&=\lim_{j\rightarrow+\infty}\lim_{\rho\rightarrow 0^+}e^{A\rho}(1+A\rho)\frac{\m(R_j\res B_{\rho})}{\rho^2}\\
    &=\lim_{\rho\rightarrow 0^+}e^{A\rho}(1+A\rho)\lim_{j\rightarrow+\infty}\frac{\m(R_j\res B_{\rho})}{\rho^2}=0.
\end{align}
Fix any $\varepsilon>0$. By \eqref{no neck property}, we can pick $j\in\n_0$ sufficiently large so that
\begin{align}
\label{eqest}
\lim_{\rho\rightarrow 0^+}\frac{\m(R_j\res B_{\rho})}{\rho^2}<\frac{\varepsilon}{2}.
\end{align}
Now assume that $\{\rho_k\}_{k\in\n}\subset(0,1)$ and $\{\rho_k'\}_{k\in\n}\subset(0,1)$ are two sequences converging to $0$ as $k\rightarrow +\infty$ and both
\begin{align*}
    (\Phi_{\rho_k})_{\ast}T_u&\rightharpoonup C_{\infty},\\
    (\Phi_{\rho_k'})_{\ast}T_u&\rightharpoonup C_{\infty}',
\end{align*}
where for every $\rho\in(0,1)$ the map $\Phi_{\rho}$ is defined as in subsection 1.2.
By further extracting subsequences if needed, we assume also that the sequences $\{(\Phi_{\rho_k})_{\ast}T_j\}_{k\in\n}$ and $\{(\Phi_{\rho_k'})_{\ast}T_j\}_{k\in\n}$ converge weakly in the sense of currents. By our previous claim, they converge to the same limit and then we have
\begin{align*}
    C_{\infty}'-C_{\infty}&=\lim_{k\rightarrow+\infty}\big((\Phi_{\rho_k'})_{\ast}R_j-(\Phi_{\rho_k})_{\ast}R_j\big)+\lim_{k\rightarrow+\infty}(\Phi_{\rho_k'})_{\ast}T_j\\
    &\quad-\lim_{k\rightarrow+\infty}(\Phi_{\rho_k})_{\ast}T_j\\
    &=\lim_{k\rightarrow+\infty}\big((\Phi_{\rho_k'})_{\ast}R_j-(\Phi_{\rho_k})_{\ast}R_j\big),
\end{align*}
in the sense of currents. By sequential lower semicontinuity of mass with the respect to weak convergence of currents, and by \eqref{eqest}, we eventually get
\begin{align*}
    \m(C_{\infty}'-C_{\infty})&\le\liminf_{k\rightarrow+\infty}\m\big((\Phi_{\rho_k'})_{\ast}R_j-(\Phi_{\rho_k})_{\ast}R_j\big)\\
    &\le\liminf_{k\rightarrow+\infty}\m\big((\Phi_{\rho_k'})_{\ast}R_j\big)+\liminf_{k\rightarrow+\infty}\m\big((\Phi_{\rho_k})_{\ast}R_j\big)\\
    &=\lim_{k\rightarrow+\infty}\frac{\m\big(R_j\res B_{\rho_k'}\big)}{(\rho_k')^2}+\lim_{k\rightarrow+\infty}\frac{\m\big(R_j\res B_{\rho_k}\big)}{\rho_k^2}<\varepsilon.
\end{align*}
By arbitrariness of $\varepsilon>0$, we obtain that $\m(C_{\infty}'-C_{\infty})=0$ and the conclusion follows.

\subsection{The general case}
Let $m,n\in\n_0$ be such that $m\ge 3$. Assume that $u\in W^{1,2}(B,\cp^n)$ is weakly $(J,j_n)$-holomorphic and locally approximable, where $B\subset\r^{2m}$ is the open unit ball in $\r^{2m}$. As usual, $\pi:B\rightarrow\cp^{m-1}$ denotes the Hopf map.

If $\theta(0,u)<\varepsilon_0$, then $u$ is smooth in a neeighbourhood of $0$ by Theorem \ref{epsilonreg} and the statement follows. Assume then that $\theta(0,u)\ge\varepsilon_0$. Moreover, since the case $n=1$ can be done exactly using the same method, we just focus on the case $n\ge2$. 

Let $T\in\D_2(B)$ be the $2$-current given by 
\begin{align*}
    \left<T,\alpha\right>:=\frac{1}{(m-2)!}\int_Bu^*\omega_{\cp^n}\wedge\pi^*\omega_{\cp^{m-1}}^{m-2}\wedge\alpha \qquad\forall\,\alpha\in\D^2(B).
\end{align*}
Notice that $T$ is well-defined and normal, since
\begin{align*}
    \left|\left<T,\alpha\right>\right|&\le\frac{1}{(m-2)!}|\alpha|_{\ast}\int_B|du|_g^2|\opwedge_2d\pi|_{g}^{2m-4}\, d\vol_g\\
    &\le C\frac{1}{(m-2)!}|\alpha|_{\ast}\int_B\frac{|du|_g^2}{|\cdot|^{2m-4}}\, d\vol_g\\
    &\le C\frac{1}{(m-2)!}|\alpha|_{\ast}\int_B|du|_g^2\, d\vol_g<+\infty, \qquad\forall\,\alpha\in\D^2(B),
\end{align*}
where $C=C\big(\Lip(\Omega)\big)>0$ is a constant and the last inequality follows from \eqref{monotonicity formula below} exactly in the same way as estimate \eqref{estimate 2m-4 number 2}.
Let $Y:=\cp^{n}\times...\times\cp^2\times\cp^{m-1}$ and $X:=\cp^1\times\cp^{m-2}$. Notice that by Lemma \ref{lemma averaging property}, Fubini's theorem, Lemma \ref{almost J-holomorphic foliations} and the co-area formula, we can write the action of $T$ as 
\begin{align*}
    \left<T,\alpha\right>&=\frac{1}{(m-2)!}\int_{Y}\bigg(\int_Bv_{y}^*\vol_X\wedge\alpha\bigg)\, d\vol_Y(y)\\
    &=\frac{1}{(m-2)!}\int_{Y}\int_X\bigg(\int_{v_{y}^{-1}(z)}\big<\alpha,\vec\Sigma^y\big>\, d\H^2\bigg)\, d\vol_X(z)\,d\vol_Y(y),
\end{align*}
for every $\alpha\in\D^2(B)$, where $y:=(q_1,...,q_{n-1},p)\in Y$ is any point in $Y$ such that (1), (2) and (3) of Lemma \ref{almost J-holomorphic foliations} hold, $v_y:=v_{q_1,...,q_{n-1},p}$ and 
\begin{align*}
    \vec\Sigma^y:=\frac{\big(\hspace{-1mm}\ast(v_y^*\vol_X)\big)^{\sharp}}{|v_y^*\vol_X|_g},
\end{align*}
following the notation that is used in Lemma \ref{almost J-holomorphic foliations}, is the $g$-unitary field of $2$ vectors orienting $v_y^{-1}(z)$, for every $z\in\reg(v_y)\subset X$. We define the "tilted current" $T_J\in \D_2(B)$ by
\begin{align*}
    \left<T_J,\alpha\right>=\frac{1}{(m-2)!}\int_{Y}\int_X\bigg(\int_{v_{y}^{-1}(z)}\big<\alpha,\vec\Sigma_J^y\big>\, d\H^2\bigg)\, d\vol_X(z)\,d\vol_Y(y),
\end{align*}
for every $\alpha\in\D^2(B)$, where $\vec\Sigma_J^y$ is the $J$-holomorphic field of $2$-vectors that we have built in the proof of Lemma \ref{almost J-holomorphic foliations}.

\textbf{Step 1}. We want to show that that $T_J$ has a unique tangent cone at the origin. First, for every $k\in\n$ we define the set $E_k\subset Y\times X$ given by
\begin{enumerate}
    \item points (1), (2) and (3) in Lemma \ref{almost J-holomorphic foliations} hold for the map $v_y$ and the level set $v_y^{-1}(z)$;
    \item $\H^2\big(v_y^{-1}(z)\cap B_{2^{-k}}\big)<\big(e^{A+2\ell}(1+A)\big)^{-1}$;
    \item $\begin{aligned}
                \int_{v_y^{-1}(z)\cap B_{2^{-k}}}|\opwedge_2d\pi(\vec\Sigma_0^y)|\, d\H^2<\frac{\delta'}{2},
            \end{aligned}$
\end{enumerate}
where we are using the notation of subsection 3.2 and $\ell=\ell\big(\Lip(\Omega)\big)>0$ is the constant provided in Lemma \ref{almost J-holomorphic foliations}. Notice that, by Lemma \ref{almost J-holomorphic foliations}, Lemma \ref{lemma projected sigma} and Fubini's theorem, it holds that
\begin{align*}
    \vol_{Y\times X}\bigg(Y\times X\smallsetminus\bigcup_{k\in\n}E_k\bigg).
\end{align*}
Fix any $k\in\n$. Define the truncated current $T_J^k\in\D_2(B)$ by 
\begin{align*}
    \left<T_J^k,\alpha\right>=\int_{E_k}\bigg(\int_{v_{y}^{-1}(z)}\big<\alpha,\vec\Sigma_J^y\big>\, d\H^2\bigg)\, d\vol_{Y\times X}(y,z), \qquad\forall\,\alpha\in\D^2(B).
\end{align*}
Notice that, by Proposition \ref{monotonicity formula sigma} and by definition of $E_k$, for every $(y,z)\in E_k$ it holds that
\begin{align*}
&e^{A\rho+\ell\rho^{1/2}}(1+A\rho)\frac{\H^{2}\big(v_y^{-1}(z)\cap B_{\rho}\big)}{\rho^2}\\
&\quad-e^{A\sigma+\ell\sigma^{1/2}}(1+A\sigma)\frac{\H^{2}\big(v_y^{-1}(z)\cap B_{\sigma}\big)}{\sigma^2}\\
&\ge\int_{v_y^{-1}(z)\cap(B_{\rho}\smallsetminus B_{\sigma})}\frac{1}{|\cdot|^2}\big<\Omega_t,\vec\Sigma_J^y\big>\, d\H^2
\end{align*}
and
\begin{align*}
&e^{-(A\rho+\ell\rho^{1/2})}(1-A\rho)\frac{\H^{2}\big(v_y^{-1}(z)\cap B_{\rho}\big)}{\rho^2}\\
&\quad-e^{-(A\sigma+\ell\sigma^{1/2})}(1-A\sigma)\frac{\H^{2}\big(v_y^{-1}(z)\cap B_{\sigma}\big)}{\sigma^2}\\
&\le\int_{v_y^{-1}(z)\cap\left(B_{\rho}\smallsetminus B_{\sigma}\right)}\frac{1}{|\cdot|^2}\big<\Omega_t,\vec\Sigma_J^y\big>\, d\H^2,
\end{align*}
for every $0<\sigma<\rho<1$. Since a direct computation leads to
\begin{align*}
    \frac{\m(T_J\res B_{\rho})}{\rho^2}=\frac{1}{\rho^2}\int_{E_k}\H^2\big(v_y^{-1}(z)\cap B_{\rho}\big)\,d\vol_{Y\times X}(y,z),
\end{align*}
by integrating on $E_k$ the two previous inequalities we get the following almost monotonicity formulas for the current $T_J^k$:
\begin{align}
\label{monotonicity formula below T_J^k}
\nonumber
&e^{A\rho+\ell\rho^{1/2}}(1+A\rho)\frac{\m(T_J^k\res B_{\rho})}{\rho^2}-e^{A\sigma+\ell\sigma^{1/2}}(1+A\sigma)\frac{\m(T_J^k\res B_{\sigma})}{\sigma^2}\\
&\ge\int_{E_k}\bigg(\int_{v_y^{-1}(z)\cap(B_{\rho}\smallsetminus B_{\sigma})}\frac{1}{|\cdot|^2}\big<\Omega_t,\vec\Sigma_J^y\big>\, d\H^2\bigg)\,d\vol_{Y\times X}(y,z),
\end{align}
\begin{align}
\label{monotonicity formula above T_J^k}
\nonumber
&e^{-(A\rho+\ell\rho^{1/2})}(1-A\rho)\frac{\m(T_J^k\res B_{\rho})}{\rho^2}-e^{-(A\sigma+\ell\sigma^{1/2})}(1-A\sigma)\frac{\m(T_J^k\res B_{\sigma})}{\sigma^2}\\
&\le\int_{E_k}\bigg(\int_{v_y^{-1}(z)\cap(B_{\rho}\smallsetminus B_{\sigma})}\frac{1}{|\cdot|^2}\big<\Omega_t,\vec\Sigma_J^y\big>\, d\H^2\bigg)\,d\vol_{Y\times X}(y,z),
\end{align}
for every $0<\sigma<\rho<1$. The inequality \eqref{monotonicity formula below T_J^k} immediately implies that the function
\begin{align*}
    (0,1)\ni\rho\mapsto e^{A\rho+\ell\rho^{1/2}}(1+A\rho)\frac{\m(T_J\res B_{\rho})}{\rho^2}
\end{align*}
is monotonically non-decreasing. Thus, the density of $T_J^k$ at $0$, given by
\begin{align*}
    \theta(T_J^k,0):=\lim_{\rho\rightarrow 0^+}\frac{\m(T_J^k\res B_{\rho})}{\rho^2}=\lim_{\rho\rightarrow 0^+}e^{A\rho+\ell\rho^{1/2}}(1+A\rho)\frac{\m(T_J^k\res B_{\rho})}{\rho^2}
\end{align*}
exists and is finite. 

We claim that $T_J^k$ has a unique tangent cone at the origin, for every given $k\in\n$. The fact that $T_J$ itself has a unique tangent cone at the origin will follow directly by the same method that is used in the conclusion of the previous subsection. By using \eqref{monotonicity formula above T_J^k}, the fact that $\Omega$ is Lipschitz, point (3) in Lemma \ref{almost J-holomorphic foliations}, the estimates \eqref{morrey decrease} and \eqref{useful estimate projected sigma 2}, we get
\begin{align*}
    &\bigg|\frac{\m(T_J^k\res B_{\rho})}{\rho^2}-\theta(T_J^k,0)\bigg|\\
    &\le C\int_{E_k}\bigg(\int_{v_y^{-1}(z)\cap B_{\rho}}\frac{1}{|\cdot|^2}\big<\Omega_t,\vec\Sigma_J^y\big>\, d\H^2\bigg)\,d\vol_{Y\times X}(y,z)\\
    &\le C\int_{E_j}\int_{v_y^{-1}(z)\cap B_{\rho}}\frac{|\Omega-\Omega_0|}{|\cdot|^2}\, d\H^2\,d\vol_{Y\times X}(y,z)\\
    &\quad+C\int_{E_j}\bigg|\int_{v_y^{-1}(z)\cap B_{\rho}}\frac{|\vec\Sigma^y-\vec\Sigma_J^y|}{|\cdot|^2}\, d\H^2\bigg|\,d\vol_{Y\times X}(y,z)\\
    &\quad+C\int_{E_j}\bigg|\int_{v_y^{-1}(z)\cap B_{\rho}}\pi^*\omega_{\cp^{m-1}}|_{\Sigma^y}\, d\H^2\bigg|\,d\vol_{Y\times X}(y,z)\\
    &\le C\vol_{Y\times X}(Y\times X)\rho^{\alpha},
\end{align*}
for every $\rho\in(0,\tilde r)$, where the constant $C>0$ and $\alpha,\tilde r\in (0,1)$ all depend just on $k$ and on $\Lip(\Omega)$. The fact that $T_J^k$ has a unique tangent cone at the origin than follows by standard arguments and step 1 is proved, due to the arbitrariness of $k\in\n$.

\textbf{Step 2}. We claim $T$ has a unique tangent cone at the origin. A direct computation using the estimate in point (3) of Lemma \ref{almost J-holomorphic foliations} leads to
\begin{align*}
    \m\big((T-T_J)\res B_{\rho}\big)\le\int_Y\int_X\int_{v_y^{-1}(z)\cap B_{\rho}}|\cdot|^{1/2}\,d\H^2\,d\vol_X(z)\,\vol_Y(y).
\end{align*}
Hence,
\begin{align}
\label{estimate T-T_J}
    \frac{\m\big((T-T_J)\res B_{\rho}\big)}{\rho^2}\le\int_Y\int_X\int_{v_y^{-1}(z)\cap B_{\rho}}\frac{1}{|\cdot|^{3/2}}\,d\H^2\,d\vol_X(z)\,\vol_Y(y).
\end{align}
By \eqref{useful estimate projected sigma 2} (recall that $\gamma=1/2)$), we get
\begin{align*}
    \int_{v_y^{-1}(z)\cap B_{\rho}}\frac{1}{|\cdot|^{3/2}}\,d\H^2\le C\H^2(v_y^{-1}(z))\rho^{1/2},
\end{align*}
for every $\rho\in(0,1)$ and for $\vol_{Y\times X}$-a.e. $(y,z)\in Y\times X$. By integrating the previous equality on $Y\times X$, \eqref{estimate T-T_J} and \eqref{integrated estimate}, we get that
\begin{align*}
   \frac{\m\big((T-T_J)\res B_{\rho}\big)}{\rho^2}&\le C\bigg(\int_Y\int_B|v_y^*\vol_X|_g\, d\vol_g\bigg)\rho^{1/2}\\
   &\le C\bigg(\int_B|du|_g^2\, d\vol_g\bigg)\rho^{1/2},
\end{align*}
where $C>0$ is a constant depending only on $m$, $n$ and $\Lip(\Omega)$. This implies that the density of $T-T_J$ at $0$, given by
\begin{align*}
    \theta(T-T_J,0):=\lim_{\rho\rightarrow 0^+}\frac{\m\big((T-T_J)\res B_{\rho}\big)}{\rho^2}=0
\end{align*}
and there is a Morrey decrease of the mass ratio to the limiting density zero. Thus, $T-T_J$ has a unique tangent cone at the origin. Since by step 1 we know that $T_J$ has a unique tangent cone at the origin and $T=T_J+(T-T_J)$, our claim follows.

\textbf{Conclusion}. Notice that $T=(m-2)!\,T_u\res\pi^*\omega_{\cp^{m-1}}^{m-2}$ and recall that $\pi^*\omega_{\mathbb{CP}^{m-1}}^{m-2}$ is invariant under $\Phi_{\rho}^*$. We address the reader to \cite[Section 7.2]{krantzparks} for the definition of the standard operations "$\res$" and "$\wedge$" when the arguments are a current and a form. Pick any two sequences of radii $\{\rho_k\}_{k\in\n}\subset(0,1)$ and $\{\rho_k'\}_{k\in\n}\subset(0,1)$ such that $\rho_k,\rho_k'\rightarrow 0^+$ as $k\rightarrow+\infty$ and
\begin{align*}
    (\Phi_{\rho_k})_{*}T_u&\rightharpoonup C_{\infty},\\
    (\Phi_{\rho_k'})_{*}T_u&\rightharpoonup C_{\infty}'.
\end{align*}
Since
\begin{align*}
    \left<(\Phi_{\rho})_{*}T,\alpha\right>&=\left<T,(\Phi_{\rho})^*\alpha\right>=(m-2)!\left<T_u\res\pi^*\omega_{\cp^{m-1}}^{m-2},(\Phi_{\rho})^*\alpha\right>\\
    &=(m-2)!\left<T_u,\pi^*\omega_{\cp^{m-1}}^{m-2}\wedge(\Phi_{\rho})^*\alpha\right>\\
    &=(m-2)!\left<T_u,(\Phi_{\rho})^*(\pi^*\omega_{\cp^{m-1}}^{m-2}\wedge\alpha)\right>\\
    &=(m-2)!\left<(\Phi_{\rho})_{*}T_u,\pi^*\omega_{\cp^{m-1}}^{m-2}\wedge\alpha\right>,
\end{align*}
for every $\alpha\in\D^2(B)$ and for every $\rho\in(0,1)$, we get that 
\begin{align*}
     (\Phi_{\rho_k})_{*}T&\rightharpoonup (m-2)!\, C_{\infty}\res\pi^*\omega_{\cp^{m-1}}^{m-2}\\
    (\Phi_{\rho_k'})_{*}T&\rightharpoonup (m-2)!\, C_{\infty}'\res\pi^*\omega_{\cp^{m-1}}^{m-2}.
\end{align*}
Since the tangent cone to $T$ at the origin is unique, we conclude that 
\begin{align*}
    C_{\infty}\res\pi^*\omega_{\cp^{m-1}}^{m-2}=C_{\infty}'\res\pi^*\omega_{\cp^{m-1}}^{m-2},
\end{align*} 
which implies
\begin{align}
\label{uniqueness tangential part}
    \big(C_{\infty}\res\pi^*\omega_{\cp^{m-1}}^{m-2}\big)\wedge\pi^*\omega_{\cp^{m-1}}^{m-2}=\big(C_{\infty}'\res\pi^*\omega_{\cp^{m-1}}^{m-2}\big)\wedge\pi^*\omega_{\cp^{m-1}}^{m-2}.
\end{align}
Notice that
\begin{align*}
    C_{\infty}&=\big(C_{\infty}\res\pi^*\omega_{\cp^{m-1}}^{m-2}\big)\wedge\pi^*\omega_{\cp^{m-1}}^{m-2}+\big(C_{\infty}\wedge\pi^*\omega_{\cp^{m-1}}^{m-2}\big)\res\pi^*\omega_{\cp^{m-1}}^{m-2},\\
    C_{\infty}'&=\big(C_{\infty}'\res\pi^*\omega_{\cp^{m-1}}^{m-2}\big)\wedge\pi^*\omega_{\cp^{m-1}}^{m-2}+\big(C_{\infty}'\wedge\pi^*\omega_{\cp^{m-1}}^{m-2}\big)\res\pi^*\omega_{\cp^{m-1}}^{m-2}.
\end{align*}
Since $m\ge 3$ we have, by dimensional considerations, that
\begin{align}
\label{vanishing projection radial part1}
    \pi_{*}\big(\big(C_{\infty}\wedge\pi^*\omega_{\cp^{m-1}}^{m-2}\big)\res\pi^*\omega_{\cp^{m-1}}^{m-2}\big)&=0,\\
\label{vanishing projection radial part2}
    \pi_{*}\big(\big(C_{\infty}'\wedge\pi^*\omega_{\cp^{m-1}}^{m-2}\big)\res\pi^*\omega_{\cp^{m-1}}^{m-2}\big)&=0.
\end{align}
Thus, by \eqref{uniqueness tangential part}, \eqref{vanishing projection radial part1} and \eqref{vanishing projection radial part2}, we get $\pi_{*}C_{\infty}=\pi_{*}C_{\infty}'$. Since $C_{\infty}$ and $C_{\infty}'$ are $J_0$-holomorphic cones, we get $C_{\infty}=C_{\infty}'$ and the statement of Theorem \ref{main theorem} follows. 
\subsection{Recovering uniqueness of tangent maps for $u$}
\,%
\newline
\textbf{The case $n=1$}. Let $m\ge 3$ and let $u\in W^{1,2}(B^{2m},\cp^1)$ be weakly $(J,j_1)$-holomorphic and locally approximable. By the methods that we have introduced in the previous subsection, it follows that uniqueness of tangent cone holds for every $(2m-2)$-dimensional current $T_{u,\psi}$ of the form
\begin{align*}
\langle T_{u,\psi},\alpha\rangle:=\int_{B^{2m}}u^*(\psi\,\omega_{\cp^1})\wedge\alpha\qquad\forall\,\alpha\in\D^{2m-2}(B),
\end{align*}
with $\psi\in C^{\infty}(\cp^1)$. 

Pick any two sequences of radii $\{\rho_k\}_{k\in\n}\subset(0,1)$ and $\{\rho_k'\}_{k\in\n}\subset(0,1)$ such that $\rho_k,\rho_k'\rightarrow 0^+$ as $k\rightarrow+\infty$ and
\begin{align*}
    u_{\rho_k}&\rightharpoonup u_{\infty},\\
    u_{\rho_k'}&\rightharpoonup u_{\infty}',
\end{align*}
weakly in $W^{1,2}(B,\cp^1)$. By uniqueness of tangent cone for $T_{u,\psi}$ we get immediately that
\begin{align*}
    \int_{B^{2m}}u_{\infty}^*(\psi\,\omega_{\cp^1})\wedge(\varphi\Omega_0)=\int_{B^{2m}}(u_{\infty}')^*(\psi\,\omega_{\cp^1})\wedge (\varphi\Omega_0),
\end{align*}
for every $\psi\in C^{\infty}(\cp^1),\,\varphi\in C_c^{\infty}(B)$. Since both $u_{\infty}$ and $u_{\infty}'$ are weakly $(J_0,j_1)$-holomorphic, by the coarea formula and by Corollary 2.1 we get 
\begin{align*}
    \int_{B}u_{\infty}^*(\varphi\,\omega_{\cp^1})\wedge(\varphi\Omega_0)&=\int_{\cp^1}\psi(y)\bigg(\int_{B}\varphi\rchi_{u_{\infty}^{-1}(y)}\, d\H^{2m-2}\bigg)\,\vol_{\cp^1}(y)\\
    \int_{B}(u_{\infty}')^*(\varphi\,\omega_{\cp^1})\wedge(\varphi\Omega_0)&=\int_{\cp^1}\psi(y)\bigg(\int_{B}\varphi\rchi_{(u_{\infty}')^{-1}(y)}\, d\H^{2m-2}\bigg)\,\vol_{\cp^1}(y)
\end{align*}
for every $\psi\in C^{\infty}(\cp^1),\,\varphi\in C_c^{\infty}(B)$. Hence, 
\begin{align*}
    \int_{\cp^1}\psi(y)\bigg(\int_{B}\varphi\big(\rchi_{u_{\infty}^{-1}(y)}-\rchi_{(u_{\infty}')^{-1}(y)}\big)\, d\H^{2m-2}\bigg)\,\vol_{\cp^1}(y)=0
\end{align*}
for every $\psi\in C^{\infty}(\cp^1),\,\varphi\in C_c^{\infty}(B)$. This implies that for $\vol_{\cp^1}$-a.e. $y\in\cp^1$ the sets $u_{\infty}^{-1}(y)$ and $(u_{\infty}')^{-1}(y)$ coincide up to $\H^{2m-2}$-negligible sets. We conclude that $u_{\infty}=u_{\infty}'$ $\L^{2m}$-a.e. on $B$ and the statement of Theorem \ref{main theorem introduction} follows.
\newline
\textbf{The case $n>1$}. Let $m\ge 3,\,n\ge 2$ and let $u\in W^{1,2}(B^{2m},\cp^n)$ be weakly $(J,j_n)$-holomorphic and locally approximable. By the methods that we have introduced in the previous subsection, it follows that uniqueness of tangent cone holds for every $(2m-2)$-dimensional current $T_{u,\psi}^{q_1,...,q_{n-1}}$ of the form
\begin{align*}
\langle T_{u,\psi}^{q_1,...,q_{n-1}},\alpha\rangle:=\int_{B^{2m}}(F_{q_1}\circ...\circ F_{q_{n-1}}\circ u)^*(\psi\,\omega_{\cp^1})\wedge\alpha\quad\forall\,\alpha\in\D^{2m-2}(B),
\end{align*}
with $\psi\in C^{\infty}(\cp^n)$ and for every choice of $(q_1,...,q_{n-1})\in\cp^2\times...\times\cp^{n}$. 

Pick any two sequences of radii $\{\rho_k\}_{k\in\n}\subset(0,1)$ and $\{\rho_k'\}_{k\in\n}\subset(0,1)$ such that $\rho_k,\rho_k'\rightarrow 0^+$ as $k\rightarrow+\infty$ and
\begin{align*}
    u_{\rho_k}&\rightharpoonup u_{\infty},\\
    u_{\rho_k'}&\rightharpoonup u_{\infty}',
\end{align*}
weakly in $W^{1,2}(B,\cp^1)$. By using the technique that we have shown for the case $n=1$, we get that 
\begin{align*}
    F_{q_1}\circ...\circ F_{q_{n-1}}\circ u_{\infty}=F_{q_1}\circ...\circ F_{q_{n-1}}\circ u_{\infty}', \qquad\L^{2m}\mbox{- a.e. on } B,
\end{align*}
for every choice of $(q_1,...,q_{n-1})\in\cp^2\times...\times\cp^{n}$. By using iteratively Lemma \ref{appendix projection lemma}, we obtain $u_{\infty}=u_{\infty}'$ $\L^{2m}$-a.e. on $B$ and the statement of Theorem \ref{main theorem introduction} follows.

\begin{rem}
\label{remark bellettini}
The advantage of the previous approach relies in the fact the we don't get only uniqueness of tangent cone for the current $T_u$ but also for its "localizations" $T_{u,\psi}$ (see the beginning of subsection 6.3) through smooth functions $\psi\in C^{\infty}(\cp^n)$. This allows more flexibility and we would like to drag the attention of the reader on the fact we could exploit such flexibility in order to get a new proof the result in \cite{bellettini-p}. Given an integer $(p,p)$-cycle $\Sigma\subset B^{2m}$, we could consider a weakly holomorphic and locally approximable map $u\in W^{1,2}(B^{2m},\cp^{m-p})$ such that $u(\Sigma)=\{y\}\in\cp^{m-p}$. By localizing the associated cycle $T_u$ through a sequence $\{\psi_{\varepsilon}\}\subset C^{\infty}(\cp^{m-p})$ such that $\psi_{\varepsilon}\rightarrow\delta_y$ in $\D'(\cp^{m-p})$, we could exploit our techniques to get uniqueness of tangent cone for $\Sigma$ ultimately.
\end{rem}
\section*{Acknowledgments}
The authors would like to thank the referees for the careful reading and for the very helpful suggestions which have contributed a lot to improve the original version of the present paper. 


\appendix

\section{Slicing through singular meromorphic maps}

Let $m\in\n$ be such that $m\ge 3$ and fix any point $p\in\cp^{m-1}$. As usual, let $\pi:\c^m\smallsetminus\{0\}\rightarrow\cp^{m-1}$ be the quotient map given by
\begin{align*}
    \pi(z_1,...,z_m):=[z_1;...;z_m], \qquad \forall z\in\c^m\smallsetminus\{0\}.
\end{align*}
Denote by $L_p$ the complex line generated by $p$ in $\cp^{m-1}$ and consider the map $T_p:\c^m\smallsetminus L_p\rightarrow L_p^{\perp}\smallsetminus\{0\}$ given by the restriction to $\c^m\smallsetminus L_p$ of the standard orthogonal projection from $\c^m$ into $L_p^{\perp}$. Fix a complex orthonormal basis $\{e_1^p,...,e_{m-1}^p\}$ of $L_p^{\perp}$ and let $\varphi_p:L_p^{\perp}\rightarrow\c^{m-1}$ be the following linear isomorphism:
\begin{align*}
    \varphi_p\Bigg(\sum_{j=1}^{m-1}\alpha_je_j^{p}\Bigg):=(\alpha_1,...,\alpha_m), \qquad \forall\, (\alpha_1,...,\alpha_{m-1})\in\c^{m-1}.
\end{align*}
Let $\pi_p:L_p^{\perp}\smallsetminus\{0\}\rightarrow\cp^{m-2}$ be the smooth submersion given by $\pi_p:=\tilde\pi\circ\varphi_p$, where
\begin{align*}
    \tilde\pi(\alpha_1,...,\alpha_{m-1}):=[\alpha_1;...;\alpha_{m-1}], \qquad \forall\, (\alpha_1,...\alpha_{m-1})\in\c^{m-1}\smallsetminus\{0\}.
\end{align*}
Eventually, notice that the map $F_p:\cp^{m-1}\smallsetminus\{p\}\rightarrow\cp^{m-2}$ given by
\begin{align*}
    F_p([z_1;...;z_m])=(\pi_p\circ T_p)(z_1,...,z_m), \qquad \forall\, [z_1,...,z_m]\in\cp^{m-1}\smallsetminus\{p\},
\end{align*}
is well-defined and smooth, since the map $\pi_p\circ T_p$ is constant on the fibres of $\pi$.

\begin{lem}
\label{properties of singular projections}
Let $m\in\n$ be such that $m\ge 3$. Then, for every $p\in\cp^{m-1}$ the following facts hold:
\begin{enumerate}
    \item the map $F_p\circ\pi$ belongs to $W^{1,2m-4}(B,\cp^{m-2})$;
    \item $F_p\circ\pi$ is weakly $(J_0,j_{m-2})$-holomorphic, where $j_{m-2}$ is the standard complex structure on $\cp^{m-2}$;
    \item $F_p\circ\pi$ is such that $d\big((F_p\circ\pi)^*\vol_{\cp^{m-2}}\big)=0$, distributionally on $B$.
\end{enumerate}
\begin{proof}
Fix any $p\in\cp^{m-1}$ and notice that the complex line $L_p$ is indeed a real $2$-plane in $\r^{2m}$. Thus, $\H^{2m-\alpha}(L_p\cap B)=0$, for every $\alpha\in[1,2m-2)$. Hence, $L_p\cap B$ has vanishing $(2m-4)$-capacity. Since $F_p\circ\pi\in L^{\infty}(B)\cap C^{\infty}(B\smallsetminus L_p)$ and the classical differential of $F_p\circ\pi$ on $B\smallsetminus L_p$ can be estimated by
\begin{align}
\label{estimate second slicing}
    |d(F_p\circ\pi)|=|d(\pi_p\circ T_p)|\le |\opwedge_2d\pi_p\circ T_p|\le\frac{C}{\dist(\,\cdot\,,L_p)},
\end{align} 
we obtain that $d(F_p\circ\pi)\in L^{2m-4}(B\smallsetminus L_p)$. Point (1) immediately follows.

For what concerns (2), we know that the weak differential of $F_p\circ\pi$ coincides $\L^{2m}$-a.e. with its classical differential on $B\smallsetminus L_p$, where $F_p\circ\pi$ is smooth. Moreover  $F_p\circ\pi=\pi_p\circ\ T_p$ on $B\smallsetminus L_p$. Since both $\pi_p$ and $T_p$ are holomorphic maps, then $F_p\circ\pi$ is holomorphic on $B\smallsetminus L_p$. Then, since $B\smallsetminus L_p$ has full $\L^{2m}$-measure in $B$, the fact that the weak differential of $F_p\circ\pi$ commutes with the complex structures $J_0$ and $j_{m-2}$ for $\L^{2m}$-a.e. $x\in B$ follows and we have proved (2).

We are just left to prove (3). Fix any $\alpha\in\D^3(B^{2m})$. For every any $\varepsilon>0$ we define
\begin{align}
    L_p^{\varepsilon}:=\big(L_p+B_{\varepsilon}(0)\big)\cap B^{2m}
\end{align}
and we notice that
\begin{align*}
    \left|\int_{B^{2m}}(F_p\circ\pi)^*\vol_{\cp^{m-2}}\wedge d\alpha\right|&=\left|\int_{L_p^{\varepsilon}}(F_p\circ\pi)^*\vol_{\cp^{m-2}}\wedge d\alpha\right|\\
    &\le\int_{L_p^{\varepsilon}}\big|d\alpha\big|\big|(F_p\circ\pi)^*\vol_{\cp^{m-2}}\big|\, d\L^{2m}\\
    &\le||d\alpha||_{L^{\infty}}\int_{L_p^{\varepsilon}}|d(F_p\circ\pi)|^{2m-4}d\L^{2m}\\
    &\le||d\alpha||_{L^{\infty}}||d(F_p\circ\pi)||_{L^{2m-3}}^{2m-4}\L^{2m}(L_p^{\varepsilon})^{1/q'}\rightarrow 0
\end{align*}
 as $\varepsilon\rightarrow 0^+$, where $q:=(2m-3)/(2m-4)$ and $q'=2m-3$ is the conjugate exponent of $q$. By arbitrariness of $\alpha\in\D^3(B^{2m})$, point (3) follows. 
\end{proof}
\end{lem} 
\begin{lem}
\label{lemma averaging property}
For every $m\ge 3$, there exists a constant $B_m>0$ such that
\begin{align*}
    \omega_{\cp^{m-1}}=B_m\int_{\cp^{m-1}}F_p^*\omega_{\cp^{m-2}}\, dp.
\end{align*}
\begin{proof}
Throughout this proof, given any $m\ge 1$ and a unitary matrix $A\in U(m)$, we will denote by $\bar A:\cp^{m-1}\rightarrow\cp^{m-1}$ the map $[z]\mapsto [Az]$.\\
It is well known that, up to rescalings by constant factors, the Fubini-Study metric is the only $U(m)$-invariant symplectic form on $\cp^{m-1}$, for every $m\ge 2$. Thus, it is enough to show that
\begin{align*}
    \bar A^*\bigg(\int_{\cp^{m-1}}F_p^*\omega_{\cp^{m-2}}\, dp\bigg)=\int_{\cp^{m-1}}F_p^*\omega_{\cp^{m-2}}\, dp, \qquad \forall\, A\in U(m).
\end{align*}
Fix any $A\in U(m)$.Given $p\in\cp^{m-1}$, define
\begin{align*}
    B_p:=\varphi_p\circ T_p\circ A\circ S_p\circ\varphi_p^{-1}:\c^{m-1}\rightarrow\c^{m-1},
\end{align*}
where $S_p:L_p^{\perp}\rightarrow\c^{m}$ is the left inverse of the orthogonal projection map $T_p:\c^{m}\rightarrow L_p^{\perp}$. As composition of linear and unitary maps, $B_p\in U(m-1)$. Moreover, by construction it holds that $F_p\circ\bar A=\bar B_p\circ F_p$.\\
Hence, by linearity of the integral and the definition of $F_p$, we have 
\begin{align*}
\bar A^*\bigg(\int_{\cp^{m-1}}F_p^*\omega_{\cp^{m-2}}\, dp\bigg)&=\int_{\cp^{m-1}}\bar A^*F_p^*\omega_{\cp^{m-2}}\, dp\\
&=\int_{\cp^{m-1}}(F_p\circ\bar A)^*\omega_{\cp^{m-2}}\, dp\\
&=\int_{\cp^{m-1}}(\bar B\circ F_p)^*\omega_{\cp^{m-2}}\, dp\\
&=\int_{\cp^{m-1}}F_p^*\bar B^*\omega_{\cp^{m-2}}\, dp\\
&=\int_{\cp^{m-1}}F_p^*\omega_{\cp^{m-2}}\, dp
\end{align*}
and the statement follows by arbitrariness of $A\in U(m)$. 
\end{proof}
\end{lem}
\begin{lem}
\label{appendix projection lemma}
Let $m\ge 3$ and pick any two points $x,y\in\cp^{m-1}$. For every $j=1,...,m$, let $\tilde e_j:=\pi(e_j)\in\cp^{m-1}$, where $\{e_1,...,e_m\}$ denotes the standard complex euclidean basis of $\,\c^m$. Assume that
\begin{align}
\label{projections equality}
    F_{\tilde e_j}(x)=F_{\tilde e_j}(y), \qquad\forall\, j=1,...,m.
\end{align}
Then, $x=y$. 
\begin{proof}
Let $x=[x_1;...;x_{m}]$ and $y=[y_1;...;y_m]$. Fix any $j=1,...,m$. By definition of $F_{\tilde e_j}$, the condition $F_{\tilde e_j}(x)=F_{\tilde e_j}(y)$ implies that there exists $\lambda_j\in\c\smallsetminus\{0\}$ such that
\begin{align*}
    \lambda_jx_i=y_i, \qquad\forall\, i=1,...,m\,\mbox{ with }\, i\neq j.
\end{align*}
Hence, by enforcing \eqref{projections equality} we get that there exists $\lambda\in\c\smallsetminus\{0\}$ such that 
\begin{align*}
    \lambda x_i=y_i, \qquad\forall\, i=1,...,m.
\end{align*}
This implies $x=y$ in $\cp^{m-1}$ and the statement follows. 
\end{proof}
\end{lem}
\printbibliography
\end{document}